\documentstyle[amscd]{amsart}

\theoremstyle{plain}
\newtheorem{Thm}[subsubsection]{Theorem}
\newtheorem{Cor}[subsubsection]{Corollary}
\newtheorem{Lem}[subsubsection]{Lemma}
\newtheorem{Claim}[subsubsection]{Claim}
\newtheorem{Prop}[subsubsection]{Proposition}

\newtheorem{Sublem}[subsubsection]{Sublemma}
\newtheorem{Fact}[subsubsection]{Fact}

\theoremstyle{definition}
\newtheorem{Def}[subsubsection]{Definition}

\theoremstyle{remark}

\newtheorem{Rem}[subsubsection]{Remark}

\newcommand{\rtil}{{\tilde r}}

\newcommand{\End}{{\operatorname {End}}}
\renewcommand{\char}{{\operatorname {char}}}

\newcommand{\rank}{{\operatorname {rank}}}
\newcommand{\corank}{{\operatorname {corank}}}

\renewcommand{\Im}{{\operatorname {Im}}}

\newcommand{\Tbar}{\overline{T}}
\newcommand{\Pbar}{{P^-}}

\newcommand{\fed}{\operatorname{ def}}

\newcommand{\U}{{\frak U}}
\newcommand{\C}{{\frak C}}

\newcommand{\Zet}{{\Bbb Z}}
\newcommand{\Ce}{{\Bbb C}}
\newcommand{\Ze}{{\cal Z}}

\newcommand{\Q}{{\Bbb Q}}
\newcommand{\Heck}{{\cal H}}
\newcommand{\grH}{gr(H)}
\renewcommand{\H}{{\cal H}}
\newcommand{\Htil}{\tilde{\cal H}}
\newcommand{\dirlim}{\varinjlim}
    
\newcommand{\Hom}{{\operatorname {Hom}}}
\newcommand{\Aut}{{\operatorname {Aut}}}
\newcommand{\supp}{{\operatorname {supp}}}
\newcommand{\proof}{{\it Proof\ \ }}
\newcommand{\epf}{\square}

\newcommand{\lambdabar}{{\overline \lambda}}

\newcommand{\Gtil}{\tilde{G}}

\newcommand{\goth}{\frak}

\newcommand{\dsh}{{\underline V}}
\newcommand{\cHom}{{\underline {Hom} }}
\newcommand{\iso}{{\widetilde {\longrightarrow} }}

\newcommand{\V}{{\Bbb V}}
\newcommand{\A}{{\cal A}}
\newcommand{\B}{{\cal B}}
\newcommand{\cC}{{\cal C}}
\newcommand{\Ag}{{\frak A}}
\newcommand{\Agbar}{{\overline{\frak A}}}
\renewcommand{\a}{{\frak a}}
\newcommand{\Z}{{\Bbb Z}}
\newcommand{\F}{{\cal F}}
\newcommand{\J}{{\cal G}}
\renewcommand{\L}{{\cal L}}
\newcommand{\R}{{\cal R}}
\newcommand{\dR}{{\cal R}\galochka}
\newcommand{\tR}{{\tilde{\cal R}}}

\newcommand{\RE}{{\Bbb R}}

\renewcommand{\P}{{\frak P}}

\newcommand{\M}{{\cal M}}
\newcommand{\Sh}{{{\cal S}{\goth h}}}

\newcommand{\cS}{{\frak S}}
\newcommand{\cG}{{\cal G}}

\newcommand{\galochka}{\check{\;}}
\newcommand{\tM}{\widetilde{\cal M}}
\newcommand{\tdM}{\widetilde{\cal M}\galochka}
\newcommand{\dM}{{\cal M}\galochka}
\newcommand{\G}{{Geom}}

\newcommand{\Xbar}{{\overline X}}

\newcommand{\sm}{{\frak sm}}
\newcommand{\k}{\underline k}
\newcommand{\I}{{\goth I}}
\newcommand{\Deltabar}{\overline \Delta}

\newcommand{\cDL}{{\cal DL }}
\newcommand{\CDL}{\check{C}{\cal DL }}

\newcommand{\prf}{\proof}

\newcommand{\Mmod}{{\frak M}}

\newcommand{\prodl}{\prod\limits}
\newcommand{\suml}{\sum\limits}
\newcommand{\cupl}{\bigcup\limits}
\newcommand{\capl}{\bigcap\limits}
\newcommand{\oplusl}{\bigoplus\limits}

\newcommand{\<}{\langle}
\renewcommand{\>}{\rangle}
\renewcommand{\(}{\left(}
\renewcommand{\)}{\right)}
\newcommand{\imbed}{\hookrightarrow}
\newcommand{\K}{K}
\newcommand{\const}{\operatorname{const}}

\renewcommand{\leq}{\ourleq}
\renewcommand{\geq}{\ourgeq}

\renewcommand{\preceq}{\ourpreceq}
\renewcommand{\succeq}{\oursucceq}

\newsymbol\treugol 134A

\newsymbol\ourpreceq 1334
\newsymbol\oursucceq 133C

\newsymbol\ourleq 1336
\newsymbol\ourgeq 133E
\newsymbol\semidir 226E
\newsymbol\boxtimes 1202

\newcommand{\mvert}{\mid}

\newcommand{\la}{{\underset{^\alpha}{\leq}}}
 
\newcommand{\leqi}{{\underset{^{i+1}}{\leq}}}

\newcommand{\beq}{\begin{equation}}
\newcommand{\eeqn}{\end{equation}}

\newcommand{\codim}{{\operatorname{\codim}\,}}

\def\square{\hbox{\vrule\vbox{\hrule\phantom{o}\hrule}\vrule}}

\title[]{Homological properties of representations of p-adic
groups related to geometry of the group at infinity}
\author{Roman Bezrukavnikov}
\begin{document}
\maketitle

\centerline{\bf Thesis submitted for the degree 
``Doctor of Philosophy''}

\bigskip

\centerline{\bf Submitted to the senate of Tel-Aviv University
September 1998.}

\bigskip

\centerline{\bf This work was carried out under the supervision of}
\centerline{\bf Prof. Joseph~Bernstein}  

\bigskip
    
\centerline{Abstract.}
                   
Let $G$ be a reductive group with compact center over a nonarchimedian
local field, and let $X$ be its Bruhat-Tits building. 

The central result of the first part of the work establishes a connection
between the direct
image of an equivariant sheaf on $X$
to the Borel-Serre compactification  with Jacquet
functor. More precisely, we prove that there exists an exact functor $ L$
on a category of representations of the group, such that the following holds.
Let $F$ be an equivariant simplicial sheaf on $X$; then the restriction
of the direct image of $F$ to the boundary of the Borel-Serre compactification
is canonically isomorphic to $L(\Gamma (F))$ (here $\Gamma $ stands for global
sections). Furthermore, the restriction of $L$ to the category of admissible
representations is described explicitly  in terms of Jacquet functor.

As an application we obtain  a canonical isomorphism between the homological
duality on the derived category of smooth $G$-modules, and the composition
of  the Deligne-Lusztig functor with the Grothendieck-Serre duality
 (the latter 
comes from the action of the Bernstein center on every smooth $G$-module). 
(This isomorphism was formulated by J.~Bernstein; another proof was recently
obtained by P.~Schneider). 

\bigskip

In the second part I prove a conjecture by D.~Kazhdan which expresses
the Euler characteristic of Yoneda Ext's between two admissible representations
as the integral of the product of their characters
over the set of elliptic conjugacy classes (another proof of this result was 
recently obtained by P.~Schneider).
 
\newpage

\centerline{\bf Acknowledgments.}

\bigskip

\noindent 
I am glad to express my deepest gratitude 
to my scientific
advisor Prof. Joseph Bernstein for his guidance and support,
and to my fellow students A.~Braverman,
D.~Gaitsgory and E.~Sayag for various discussions and explanations.
All I know about representations of $p$-adic groups I know from them.

\bigskip

\noindent I wish to thank  Prof. R.~Kottwitz and Prof. M.-F.~Vigneras
 for the interest in the work they have expressed,
 and Prof. P.~Deligne for helpful
critical comments (in particular, for pointing out a gap in the 
preliminary version of the arguments).

\bigskip

\noindent Last but not least, I thank  M.~Finkelberg and A.~Stoyanovskii;
besides all the mathematics they have taught me during the past ten
years, I owe them taking the pain of reading the manuscript and pointing
out some misprints and expositional lapses.

\bigskip

\noindent A part of this work was written during the author's stay
 at the Institute for Advanced Study  in Princeton.
 I wish to
 thank this institution for the  excellent working conditions.

\newpage

\section*{Introduction}

\noindent 1.
Let $F$ be a non-archimedian local field, and $G=\underline G(F)$ be the 
group of $F$-points of a reductive algebraic group $\underline G$ defined 
over $F$. We assume for convenience that $G$ has compact center.
 
 We fix an algebraically closed 
 coefficient field $k$ of characteristic $0$.
 $\H$ will denote the Hecke algebra of
$k$-valued locally constant compactly supported measures on $G$.

The main object of our study is the category $\M$ of finitely generated
smooth representations of $G$ over $k$ \cite{BZ}.

\hfill

\noindent 2. In chapter 1 the category $\M$ is studied using equivariant
 sheaves on the Bruhat-Tits building $X$ of $G$.

For a $G$-space $Y$ we will denote by $Sh_G(Y)$ the category of $G$-equivariant
sheaves on $X$. 
 Let $\Sh\subset Sh_G(X)$ be the category of simplicial
$G$-equivariant sheaves of finite-dimensional vector-spaces 
on $X$. We have the derived functor of sections
with compact support $R\Gamma_c:D^b(\Sh)\to D^b(\M)$, and Verdier dual
functor $R\Gamma:D^b(\Sh)\to D^b(\M^{opp})$ (where $\M^{opp}$ is the opposite
category to $\M$).

Consider the imbedding $j:X\to \Xbar $ of $X$ into its Borel-Serre
compactification $\Xbar$ \cite{BSer}, and the equivariant direct image functor
$j_*^G:D^b(\Sh)\to D^b(Sh_G(\Xbar))$ (discussed below).
 Let $i$ stand for the imbedding
$i:\Xbar-X \imbed \Xbar$. The main technical result of 
chapter 1 is existence of a functor
$\L:\M^{opp}\to Sh_G(\Xbar-X)$ such that:

\smallskip

i)  We have $i^*\circ j_*^G = \L \circ \Gamma$ canonically.

ii) $\L$ is exact.

iii)  Let $\R\subset \M$ be the subcategory of admissible representations.
Then $\L|_{\R^{opp}\cong \R}$ is a well-known combinatorially defined functor:
for $M\in \R$ and $y\in \Xbar-X$ the stalk of $\L(M)$ at $y$ is the Jacquet
functor $r_{P_y}(M)$ where 
$P_y=Stab_G(y)$ is a parabolic subgroup in $G$.
(The equivalence between $\R$ and $\R^{opp}$ used above sends a representation
$\rho\in \R$ to its contragradient $\rho \galochka$).

We remark that construction of $\L$ and proof of i) and iii) are rather 
straightforward, while the proof of ii) required some effort.

\hfill

To apply this result to representation theory one needs to be able to
``localize'' representations to sheaves on the building.
An ideal result in this direction would be a one-sided inverse to the
functor $R\Gamma_c$. Unfortunately such functor is unavailable at present.
By the results of \cite{Schn}, \cite{Schn1}
 the full subcategory 
 $D^b(\R)\subset D^b(\M)$ lies in the image of 
$R\Gamma_c$. (And also $D^b(\M)$ lies in the image of the functor
$R\Gamma^c$ from $G$-equivariant simplicial sheaves of 
possibly  infinite
dimensional vector spaces to 
$ D^b(\tM)$, where $\tM$ is the larger category of 
{\it all} (not necessary finitely generated) smooth $G$-modules).

The result of \cite{Schn1} is sufficient for our application.
By a ``general nonsense'' argument we get  another localization type
claim sufficient for us
(but much less explicit, and
thus probably  less useful than the Schneider-Stuhler's result). It
 says the following: for any object 
$B\in D^b(\M)$ and $i\in \Zet$
 one can find an object $A\in D^b(\Sh)$ and a morphism
$R\Gamma_c(A)\to B$ which induces isomorphism on $H^j$ for $j>i$.

Using the connection between the Verdier duality
on sheaves with the  homological duality on representations
(also observed in Schneider's work \cite{Schndual}) 
one derives the following statement (which belongs essentially to
J.~Bernstein and was  the starting point for this part of the work):
there exists a  canonical isomorphism of functors:
$$D_h\cong DL\circ D_{Gr}$$
where $DL$ is the Deligne-Lusztig functor \cite{DL};
 $D_h$ is the homological duality
$D_h:M\to RHom_G(M, C_c^\infty(G))$; and $D_{Gr}(M)$ is the Grothendieck
dual of $M$ as a module over the Bernstein center \cite{BDKV} (for
admissible $M$ we have $D_{Gr}(M)=M\galochka$ where $M\galochka$ is the
contragradient module). (See the main text for  detailed definitions).

A proof of the latter isomorphism
using (co)sheaves on the Borel-Serre
 compactification of the Bruhat-Tits building
appears also  in \cite{Schn1}.

\hfill

\noindent 3. The starting point for  Chapter 2 was an attempt to find an
 ``a'priori'' proof for  
 the following conjecture by Kazhdan (proved recently  in \cite{Schn1} by
 making
use of an   explicit resolution of a $G$-module
constructed in  \cite{Schn1}, \cite{Schn}). 

\medskip

Let $Ell$ be the set of regular semisimple elliptic conjugacy classes
in $G$. Then $Ell$ carries a canonical measure (an analogue of the Weyl
integration measure for compact Lie groups,
 see section \ref{ellsect} of chapter 2
  for a recollection), which we
denote by $d\mu$. For an admissible $G$-module $\rho$ let $\chi_\rho$
denote its character, and let us use the same notation for the
corresponding locally constant function on $Ell$.
By \cite{K} we have $\chi_\rho \in L^2(Ell,d\mu)$ provided the base local
field $F$ has characteristic 0. Let now $\rho_1$, $\rho_2$ be two
admissible representations. Then the conjecture (= Theorem \ref{Kazhdan}
of chapter 2) says that
 \begin{equation}\label{Eulchar}
\sum (-1)^i \dim
Ext^i(\rho_1,\rho_2)=
\int_{Ell} \chi_{\rho_1}(g^{-1})\chi_{\rho_2}(g)d\mu(g)
\end{equation}
 
%The plan of our proof of \eqref{Eulchar} is as follows.
\smallskip

Recall the following construction (see e.g. \cite{Br}, \S IX.2
or \cite{Vigneras} \S 2 and references therein). 
Let $\A$ be a Noetherian associative
algebra 
with 1, having finite homological dimension.
 Then to a pair $(M,E)$ where $M$ is a finitely generated 
$\A$-module, and $E\in End(M)$ one can assosiate an element of
$\A/[\A,\A]$, denoted by $Tr_{H-St}(M,E)$, and called  the {\it
Hattori-Stallings trace} of $(M,E)$. This map enjoys the properties:
 
a) $Tr_{H-St}(M,E_1\circ E_2)=Tr_{H-St}(N,E_2\circ E_1)$ for two modules $M,N$
as above, and $E_1\in Hom(N,M)$, $E_2\in Hom (M,N)$;

b) $Tr(\A,R(a))=a \mod
[\A,\A]$ for $a\in \A$ where $R(a)$ is the right multiplication by $a\in A$
considered as an endomorphism of the free module;

c)  additivity on short
exact sequences. 

It is characterized by these properties uniquely.

Applying this construction
 (with slight modifications  necessary
because $\H$ is not unital, see \ref{HaS} of chapter 2)
 to our situation we obtain a map $Tr_{H-St}:
\{(\Mmod,E)\ |\ \Mmod \in \M, E\in End(\Mmod)\}\to \H/[\H,\H]$ satisfying above
properties a), b), c). For $\Mmod \in \M$ the element $Tr_{H-St}(\Mmod,Id)$
is called a pseudo-coefficient or rank function of $\Mmod$;
 we will denote it by $\< \Mmod \>$.

Recall that for $h\in \H$ and a regular elliptic element $g\in G$
the {\em orbital integral} $O_g(h)\in k$ is defined by
$$O_x(h):=
\int\limits_{y\in G} \frac{h}{\nu}
(yxy^{-1})d\nu(y),$$
where $\nu$ is a Haar measure on $G$; the resulting
expression does not depend on the choice of $\nu$.

A formal argument shows (see section 
\ref{ellsect} of chapter 2) that \eqref{Eulchar} is equivalent
to the equality
\begin{equation}\label{Ochi}
O_{g^{-1}}(\<\rho\>)= \chi_\rho(g)
\end{equation}
 being true for any admissible representation $\rho$
and regular  elliptic element
$g\in G$.

Notice that the latter statement 
makes sense (and will be proved) without the restriction $\char(F)=0$ or 
$k=\Ce$.

If $\rho$ is cuspidal, and $m_\rho \in \H$ is its matrix
coefficient normalized by $\int m_\rho\cdot \chi_\rho=1$, then 
$\<\rho\>= m_\rho \mod [\H,\H]$. Thus \eqref{Ochi} generalizes
a well-known equality going back to Harish-Chandra
$$
\chi_\rho(g)=O_{g^{-1}}(m_\rho).
$$

To prove \eqref{Ochi} we extend the functional $Tr(g,\rho)=\chi_\rho(g)$
from admissible to arbitrary smooth representations $\rho$, and then show 
(Theorem \ref{sumEii}) that 
for a projective $\rho\cong \H^{\oplus n} E$, where $E\in Mat_n (\H)$
is an idempotent we have 
\begin{equation}\label{sumeii}
Tr(g,\rho)=O_{g^{-1}}(\sum E_{ii}).
\end{equation} 

\medskip

The intuitive idea behind the definition of $Tr(g,\Mmod)$ is that for a small
open compact subgroup $K\subset G$ one can find a finite dimensional 
$g$-invariant subspace 
$M_0\subset \Mmod^K$ such that the quotient $\Mmod^K/M_0$ ``comes from the
parabolic induction''; thus it is reasonable to set $Tr(g,\Mmod)= tr(g,M_0)$.
Realization of this idea exploits certain natural (multi)filtration
on the Hecke algebra connected with the action of $\H$ on the principal 
series (see section 2 of chapter 2).

\medskip

 To prove \eqref{sumeii} we relate our filtration to the one defined
``geometrically'' in terms of support of a distribution $h\in \H$. 
Here some geometry of $G$ and its asymptotic cones $(G/U\times G/U^-)/L$
at infinity is used.

\newpage

\begin{section}{Sheaves on the building and representations}\label{1}

\subsection*{Preliminaries and notations}
In both chapters we will abuse the terminology by saying ``a parabolic/Levi
subgroup'' 
instead of ``the group of $F$-points of an $F$-rational parabolic
subgroup/Levi subgroup of a proper parabolic subgroup''. 
By rank of a reductive group defined over $F$ we mean its split rank. 

\subsubsection{} \label{prum}
$k$ is a fixed coefficient field of
characteristic 0; for a  set $S$ we denote by $k[S]$ the space $\oplusl_S
k$.
%of $k$-valued functions on $S$ with finite support.
 For a topological space $Z$ we write
$C^\infty(Z)$ and $C^\infty_c(Z)$ for the space of locally constant
functions on $Z$ and   locally constant
functions with compact support on $Z$. 
If $Z$ is an open subset of a group, we denote by $\Heck(Z)$ the space of 
locally constant measures on $Z$ with compact support; if
$Z$ is biinvariant under the action of a compact open subgroup $K\subset
G$, then $\Heck(Z,K)\subset \Heck(Z)$ is the subspace of $K$-biinvariant
measures. If $Z$ is a semigroup, then $\H(Z)$ and $\H(Z,K)$ carry the
structure of an associative algebra provided by convolution. 

By $\underline W$ or $\underline W_Z$
we will denote the constant sheaf on $Z$ with stalk $W$;
 if $\iota :Z\to X$ is a closed
embedding then we will also write $\underline W_Z$ instead of $\iota_*(
\underline W_Z)$.
By $\cHom$ we denote internal $Hom$ in the category of ($G$-equivariant)
sheaves.

For a $G$-space $X$ let $Sh_G(X)$ denote the category of $G$-equivariant
sheaves on $X$. Recall 
that for arbitrary topological group 
$G$ and a $G$-space $X$ 
a $G$-equivariant sheaf  is the same as a sheaf $\F$ on $X$
equipped with a $G$-action such that  each section of $\F$ is locally
constant. Here a section   $s\in \Gamma (U,\F)$ is called locally
constant if for  $x\in U$
there exist a neigborhood $U'$ of $x$ and a neigborhood $V$ of identity
$e\in G$ such that $g(s)|_{U'\cap g(U')}=s|_{U'\cap g(U')}$ for $g\in V$.
 
For a morphism of $G$-spaces $p:X\to Y$ we have the functor of inverse image
$p^*:Sh_G(Y)\to Sh_G(X)$ and the adjoint functor $p_*^G:Sh_G(X)\to Sh_G(Y)$.
Note that the inverse image commutes with the
 forgetful functor from equivariant
sheaves to sheaves while the direct image does not:
 for $\F \in Sh_G(X)$
one can describe $p^G_*(\F)$ as the sheaf of locally constant sections
of  $p_*(\F)$ (see \cite{Schnequiv} for a detailed discussion).

 From now on let $X$ denote the (semi-simple) building of $G$, unless stated
otherwise. 
Let $\Xbar$ be the Borel-Serre compactification of $X$  [BorSe];
put $Y=\Xbar- X$, and denote  the imbeddings by
$j:X \hookrightarrow \Xbar$ and $i:Y\hookrightarrow \Xbar$.
Recall that $Y$ is the spherical building of $G$. 
 For a parabolic $P\subset G$ we let $\Delta_P\subset Y$ denote the
corresponding simplex.

We will often say polysimplex meaning a polysimplex of the canonical
polytriangulation of $X$ (facette in the terminology of \cite{BT1}).

Let $\Sh$ be the full subcategory of $Sh_G(X)$ consisting of sheaves
of finite dimensional vector spaces which are constant on every (poly)simplex.

\subsubsection{} 
By a representation we will mean a representation  in a vector space
over a field $k$; the field  will be endowed with discrete topology.

Let $\tM$ be the category of smooth $G$-modules; let $\tdM$ be the dual
category. It will be convenient to use an explicit realization of $\tdM$.
Namely $\tdM$ is identified with the  category of
complete  topological vector spaces $V$ having a basis of neighbourhoods of 0
consisting of vector subspaces of finite codimension, together with 
a continuous homomorphism $\pi: G\to {\operatorname {Aut}}(V)$. (The  group
$ {\operatorname {Aut}}(V)$ is endowed with the topology of pointwise
convergence.) 
The antiequivalence  $*:\tM \to \tdM $ sends a $G$-module $M$
to its full linear dual equipped with the pointwise-convergence topology
(= the inverse limit topology on $M^*=\varprojlim V^*$, where $V$ runs over 
the set of finite dimensional subspaces in $M$).
 The inverse antiequivalence (also denoted by $*$)
sends a topological $G$-module $L$ to the space of
continuous linear functionals $L \to k$.

Let $\M \subset \tM$ (respectively $\R \subset \M$) be the full 
subcategory
of finitely generated smooth (respectively admissible) $G$-modules.
Let $\dM \subset \tdM$ and $\dR \subset \dM$ be the
dual categories.

We have an exact functor $Sm: \tdM \to \tM$ sending a module $M
\in \tdM$ to the module of smooth vectors in $M$.
The functor $Sm$ induces an equivalence $\dR \widetilde {\to} \R$.
(Of course  it does {\it not} send the category $\dM$ to $\M$.) 
 For $M\in \M$ we let $\M\galochka =Sm(*(M))\in \tM$
 be its contragradient.

 We have left-exact functors $\Gamma _c :Sh\to \M$ and $\Gamma : Sh \to
\dM$, where $\Gamma$ stands for global sections, and $\Gamma _c$ stands
for sections with compact support. The linear topology on the space
$\Gamma (X,{\cal F})$ is the inverse limit topology on
$\Gamma (X,{\cal F})= \lim\limits_{\longleftarrow K}   \Gamma (K, \F)$
where  $K$ runs over all compact subsets of $X$.

As we will soon recall,  $\Sh$ has enough injectives (and projectives).
 Hence the right
derived functors $R\Gamma:D^*(\Sh)\to D^*(\tM)$
and  $R\Gamma_c:D^*(\Sh)\to D^*(\M)$ where $*=b,+$
are defined. 

\subsection{Verdier duality for equivariant sheaves and homological
 duality}
Here we recall some basic generalities on homological algebra of $\Sh$,
 including Verdier duality and its connection with homological duality on
representations.
This section does not contain new results:
Proposition \ref{ve}a) follows from \cite{Schnequiv}, \S 3 and remark
before \cite{Schndual}, Proposition 2;
our \ref{ve}b) is Proposition 4 in \cite{Schndual}. We sketch the proof for
completeness.

\subsubsection{}\label{Verdier}
Let $Sh(X)$ be the category of sheaves of finite dimensional
$k$-vector   spaces on $X$ constant on any polysimplex of the canonical
polytriangulation. Then $Sh(X)$ has enough injectives; the injective
generators are the constant sheaves $\k _{\overline \Delta}$ on the closure
of a simplex. 
From the fact that the stabilizer of any point of $X$ is compact, and hence
its smooth representations form a semisimple category, it is easy to deduce
 that a sheaf $\F\in \Sh$ is injective iff it is  injective as
an object of $Sh(X)$.  For a polysimplex $\Delta \subset X$ and a smooth
finite dimensional  
representation $R$ of the stabilizer of $\Delta$ let $\F_{\Delta,R}\in \Sh$
be the $G$-equivariant sheaf
 $\F_{\Delta,R}:=\prod \limits _{g\in G/Stab(\Delta)}\underline k
_{g\overline \Delta} \otimes C^\infty\(g\cdot Stab(\Delta)\)\bigotimes \limits
_{Stab(\Delta)} R$. Obviously $\F_{\Delta,R}$ is injective.

Let $\I$ be the category of injective objects in $Sh(X)$, and  $\I_G$
be the category of injective objects in $\Sh$.

Let us call a module $M\in  \dM$ (respectively $M\in \M$) standard injective
(respectively projective) if it is a direct sum of modules induced
(respectively compactly induced) from smooth finite dimensional
representations of open  compact subgroups. 

\begin{Lem}\label{obv} a)  $\Sh$  is of finite homological dimension and
has enough injectives. Any injective object in $\Sh$ is a direct sum
of sheaves of the form $\F_{\Delta, R}$.

b) $\Gamma:\Sh \to \dM$ sends injectives to injectives and $\Gamma_c:\Sh
\to \M$ sends injectives to projectives.

c) An injective (respectively projective) object of $\dM$ (respectively
$\M$) is isomorphic to $\Gamma (\F)$  (respectively
$\Gamma_c(\F)$) for  $\F \in \I_G$ iff it is a standard injective
(projective). 

d) If $\F,\J \in \I_G$  then $\cHom(\F,\J)$ and $\F\otimes \J$
are also injective.

\end{Lem}
\proof Clear.
\epf 

We see that any left-exact functor on $\Sh$ admits the
(right) derived functor.

\begin{Prop}\label{ve} a)
There exists an object $\dsh \in D^b(\Sh)$ 
equipped with an isomorphism of functors on $D^b(\Sh) $
 \begin{equation}\label{verdier}
R\Gamma \circ  R\cHom (\underline {\ \ }, \dsh)\cong * \circ R\Gamma_c
\end{equation}

Isomorphism \eqref{verdier} defines $\dsh$  uniquely  up to the unique
isomorphism. 
 $\dsh$ is called the equivariant dualizing sheaf, and the functor
  $\V := R\cHom (\underline {\ \ }, \dsh)$ is called the 
equivariant Verdier duality. 

b) Let $D_h:D^b(\M)\to D^b(\M)$ be the homological duality defined by
$D_h(M)=RHom(M, \Heck)$. We have a canonical isomorphism:
\begin{equation}\label{verdier_homolo}
R\Gamma _c \circ \V \cong D_h \circ R\Gamma _c
\end{equation}
\end{Prop}
\proof a)  Obviously $Hom_{\Sh} (\F, \J)= (\Gamma (\cHom (\F,\J)))^G$.
 Using Lemma \ref{obv}b),d) we see that $RHom (\F,\J)=RI \circ R\Gamma (R\cHom
(\F,\J))$ where $I$ stands for invariants $I:\dM \to Vect$.  
So \eqref{verdier} yields an isomorphism of functors $Hom(\underline{\ \ },
\dsh) \cong H^0(RI\circ *\circ R\Gamma_c )$; uniqueness of $\dsh$
follows. 

To show existence we adapt a standard construction of the dualizing sheaf
to our situation. For a polysimplex $\iota:\Delta\hookrightarrow X$ let
$C_\Delta$ be the standard  resolution of $\iota_!(\underline
k[\dim(\Delta)])$ (so we have $C_\Delta^i = \bigoplus \limits_{\Delta'\subset
\overline \Delta; \dim (\Delta') = -i} \underline k _{\overline \Delta
'}$).

We set $\V(\k _{\Deltabar})=C_\Delta$; this obviously extends to a
contravariant functor
 $Kom(\I)\to Kom(\I)$ and defines  a contravariant
 functor $Kom(\I_G)\to Kom(\I_G)$. (We have $\V(\F_{\Delta,R})=\prod
\limits _{g\in G/Stab(\Delta)} C_
{g \Delta} \otimes C^\infty \(g\cdot Stab(\Delta)\)\bigotimes \limits
_{Stab(\Delta)} R^*$.)
 The induced functor  $D^b(\Sh)\to D^b(\Sh)$
is also denoted by $\V$.
One can easily define  isomorphism of complexes
 $\V(\F \otimes \J)\cong \cHom(\F, \V (\J))$ for $\F,\J\in \I_G$;
the isomorphism of bifunctors on $D^b(\Sh)$ follows. Hence $\V$
is corepresented by $\dsh = \V (\k)$.

For $\F\in \I_G$ we have a functorial map of complexes $\Gamma(\V(\F))\to
*(\Gamma_c(\F))$ which is a quasiisomorphism; here the pairing $\langle\
,\ \rangle: \Gamma_c (\F)
\times \Gamma (\V(\F))\to k$ is the standard Verdier pairing. This yields
isomorphism \eqref{verdier}.

b) The pairing  $\langle\ ,\ \rangle$ gives a map 
$\phi:\Gamma (\V(\F))\to Hom_G(
\Gamma_c(\F), C^\infty(G))$; here $\(\phi (y)(x)\)(g)=\langle x, g(y)\rangle$. 
It is clear that $\phi(\Gamma_c(\V(\F)))\subset Hom _G(
\Gamma_c(\F), C^\infty_c(G))\subset Hom_G(
\Gamma_c(\F), C^\infty(G))= \Gamma (\V(\F))$; more precisely for $\F\in
\I_G$ we get a functorial quasiisomorphism $\Gamma_c(\V(\F)))\to Hom _G(
\Gamma_c(\F), C^\infty_c(G))$ which yields \eqref{verdier_homolo}.
(In fact for $\F \in \I_G$ isomorphisms \eqref{verdier}, \eqref{verdier_homolo}
boil down to well-known equalities $$*(ind_K^G(R))=Ind_K^G(R^*),$$
$$Hom_G(ind_K^G(R), C^\infty_c(G))=ind_K^G(R^*).)$$

Proposition is proved.
\epf

\begin{Rem} In \cite{Schndual} it is proved that $\dsh$ is actually
concentrated in one homological degree; we will not use this fact here. 
\end{Rem}

\subsection{A ``localization'' type theorem}
 Let us call a morphism in $D^b(Sh)$  {\it
h-quasiisomorphism} (respectively {\it hc-quasiisomorphism}) if it induces
an isomorphism on $R\Gamma$ (respectively $R\Gamma_c$).

 Let $D^0(\dM)$, $D^0(\M)$ be the full subcategories in
$D^b(\dM)$, $D^b(\M)$ respectively
consisting of objects quasiisomorphic to a finite complex
of standard injectives (respectively projectives).
 $LD^b(\Sh)$, $L^0D^b(\Sh)$ be the localizations of  $D^b(\Sh)$
respectively with respect to h-quasiisomorphism and hc-quasiisomorphisms.

\begin{Prop}\label{loc} $R\Gamma$ provides an equivalence between
$LD^b(\Sh)$
and $D^0(\dM)$.

$R\Gamma_c$  provides an equivalence between
$L^0D^b(\Sh)$
and $D^0(\M)$.
\end{Prop}

For the proof of the Proposition we need the following

\begin{Lem}\label{blow} Assume that $\F, \J \in \Sh$ and $\F$ is injective;
let a morphism $\phi: \Gamma (\F) \to \Gamma (\J)$ be given.
Then there exist a  sheaf $\F' \in \Sh$ and morphisms 
$s: \F ' \to \F$, $\phi ' :\F ' \to \J$  such that $s$ is an
 h-quasiisomorphism,
 and $\Gamma (\phi ')= \phi \circ \Gamma (s)$.
\end{Lem}

{\sf Proof:} First notice that $X\times X$ is a union of poly-simplicial
subspaces $X_r$ such that $pr_i|_{X_r}$ is proper  for $i=1,2$
and  $pr_1|_{X_r}$ has geodesically contractible fibers.
 (To see this
let us choose an apartment $\Ag\subset X$, an open simplex $\Delta\subset
\Ag$,
and exaust $\Ag$ by  simplicial convex
neigborhoods $\C_r$ of $\overline\Delta$. Then take
$X_r=G(\overline \Delta\times \C_r)$. The desired properties of
$X_r$ are clear. We have $\cupl_r X_r=X\times X$  because
$G(\Deltabar)=X$ and  any two 
points of $X$ are $G$-conjugate to two points lying in $\Ag$.)
  Let $p_i=pr_i|_{X_r}$, $i=1,2$. 
 Set $\F _r =p_{2*} p_1^* \F$.

 We claim that for large enough $r$ we can take $\F ' = \F _r$.

It is clear that if $\F$ is injective then $R^ip_{2*}p_1^*\F=0$
 for $i\not =
0$. Hence $R\Gamma (\F_r)= R\Gamma ( p_1^* \F)=R\Gamma (\F)$ since $p_1$
has contractible fibers. Applying $p_{2*}$ to the morphism $p_1^* \F \to
\delta _* (\F)$ where $\delta : X\to X\times X$ is the diagonal imbedding 
we get a map $\F_r \to \F$ which induces isomorphism on cohomology. 

It remains to construct $\phi': \F_r \to \J$ for large $r$.

For any $x\in X$ consider the composition  $\Gamma(\F)\to \Gamma(\J)\to
\J_x$ where $\J_x$ is the stalk of $\J$ at $x$.
It is continuous and $\J_x$ is discrete;
 hence it factors through a map  $\Gamma (\F)\to \Gamma (\C, \F)\to
\J_x$ for some compact $\C=\C(x)\subset X$.
Then clearly the same is true for any $y$ lying in one polysimplex
with $x$ with the same compact set $\C(y)=\C(x)$. 

We can take a finite number of polysimplicies $\Delta_i$ such that $G(\cup
\Delta_i) =X$; and then find $r$ such that $X_r\supset \C(x_i)\times
\Delta_i$ for $x_i \in \Delta_i$. Then we see that the map
$\Gamma(\F_r)=\Gamma (\F) \to \J_x$ factors through $\Gamma(\F_r)\to (\F_r)_x$
for all $x\in X$, i.e. the composition $\k \otimes \Gamma (\F_r)\to
\k \otimes \Gamma (\J)\to \J$ factors through the canonical surjection
 $\k \otimes
\Gamma (F_r)\to \F_r$. The Lemma is proved.
\epf

\prf  of the Proposition.
  First note that image of $R\Gamma$ is contained in $D^0(\dM)$
since an injective sheaf goes to a standard injective module. 

Lemma \ref{blow} implies that for $\F,\J\in \I_G$
the map
$\Hom_{LD^b(\Sh)}(\F,\cG)\to \Hom_\dM(R\Gamma \F, R\Gamma \J)=
\Hom_\dM(\Gamma\F,\Gamma\J)$ is surjective. 

It is also injective for the following reason. Consider the commutative square
$$
\begin{CD}
{\cal F}           @>{s}>> \J       \\
@AAA                         @AAA         \\
\underline{\Gamma(\F)}   @>{\underline{\Gamma(s)} }>>
\underline{\Gamma(\J)}
\end{CD}
$$
where $\Gamma(s):\Gamma(\F)\to \Gamma(\J) $ is the map on global sections
induced by $s$, and $\underline{\Gamma(s)}$ is the corresponding map of
constant sheaves. 
The left vertical arrow in this diagram
 is surjective for injective $\F$, so we have $\Gamma(s)=0,\, \F\in \I_G
\Rightarrow s=0$.

Thus $R\Gamma$ induces an isomorphism
$\Hom_{LD^b(\Sh)}(\F,\cG)\iso \Hom_\dM(R\Gamma \F, R\Gamma \J)$ for $\F,\J
\in \I_G$,
and hence for any $\F,\, \J \in D^b(\Sh)$.

It remains to check that $R\Gamma:LD^b(\Sh)\to D^0(\dM)$ is surjective
on isomorphism classes of objects.
Let $C$
be a finite complex of standard injectives $0\to I_a\to...\to I_b\to 0$. We
then have some injective 
sheaves $\F_i \in \Sh$ with $\Gamma (\F_i)\cong I_i$.
We now use inductively
Lemma \ref{blow} to construct a complex of sheaves $\tilde C =
0\to \F_a'\to..
\to \F_b'\to 0$, where $\F_i'=(\F_i)_{r_i}$, $H^j(F_i')=0$ for $j\not =0$
 and $\Gamma (\tilde C)\cong C$.

 The first statement is proved, and the second one follows by Verdier
duality. 
\epf

\begin{Rem} I do not know whether in general 
$D^0(\M) = D^b(\M)$, i.e. whether  a
localization type theorem stronger (but  less explicit) than the one
obtained by Schneider and Stuhler in \cite{Schn},\cite{Schn1} holds.
 The positive answer to this question is equivalent to
 the claim that  the Grothendieck group of $\M$  is generated by the
 classes of standard projectives.

 Notice that the following unpublished result of J.~Bernstein agrees with the
positive answer to this question: the Grothendieck group of $\M$
tensored with $\Bbb Q$ is generated by the classes of standard projectives.
\end{Rem}

\begin{Rem} Our method also does not answer the question whether $D^+(\dM)$
or $D^+(\M)$ is localization of $D^+(\Sh)$. (Any module has infinite right
resolution by standard injectives; our procedure of constructing a complex
of sheaves from such a resolution was inductive
``moving from  right to the left''.)
\end{Rem}

\subsection{Main theorem}\label{centr}
 We start the construction of the functor $\L$ promised in the introduction. 
 
Let $Z$ be a closed subset of $\overline X -X$.

Consider the subspace $\Gamma _c (\overline X -Z,
j_*({\cal F}))
\subset \Gamma (X,{\cal F})$ and the subspace $\Gamma _{c} (\overline X -Z,
j_*^G{\cal F})
\subset \Gamma (\overline X,j^G_*{\cal F})=Sm(\Gamma (X,{\cal F}))$.
(The latter equality readily follows from compactness of $\Xbar$.)

We now explain how to recover these subspaces from an object
$\Gamma (X, {\cal F}) \in \dM$.

  Fix a point $p\in X$. Let
 $a: G\times X \to X$ be the action; set $a_p (g)=a(g,p)$.
 Let ${\goth F} _Z$ denote the filter of subsets in $G$
generated by the sets of the form $a_p ^{-1} (U)$, where $U$ is 
 a neighborhood of  $Z$.

For a module $L \in \dM$ define $L_Z = \{m\in L\; | \;{\lim}_{\goth F}
(g^{-1}(m))=0\}$, and $L_Z^\sm=Sm(L)\cap L_Z$. (Here  $ {\lim}_{\goth F}$
stands for the limit with respect to the filter $\goth F$.)

\begin{Lem}\label{premain} a) The filter $\goth F$ is invariant under right
 translations on $G$. It does not depend on the choice of a point $p\in X$.

b) We have:

  $\Gamma _c (\overline X -Z,
j_*({\cal F})) = (\Gamma (X,{\cal F}))_Z$

 $\Gamma _{c} (\overline X -Z,
j_*^G({\cal F})) = ( \Gamma (X,{\cal F}))_Z^\sm$
\end{Lem}

{\sf Proof}
To prove a) we need the following well-known property
of $\overline X$.

\begin{Fact} Assume that a sequence of points $x_n \in X$ 
has a limit $x\in \overline X-X$. Let $y_n \in X$ be another
sequence, such that $d(x_n,y_n)$ is bounded. Then $y_n \to x$
as well. 
\end{Fact}
\proof is included since no reference was found. We use Remark 5.5(2) in
\cite{BSer}. It says the following.
Fix an open polysimplex $\Delta_0\subset X$, and for any apartment $\Ag
\subset X$ containing $\Delta_0$ let $\Xi_\Ag:X\to \Ag$ be the contraction
of $X$ to $\Delta_0$ \cite{BT1} 2.3.4, 7.4.19. Recall that an apartment
$\Ag$
is an affine space, and let $\Agbar$ be its standard semisphere
compactification (see below); denote by $j_\Ag$ the imbedding $\Ag\imbed
\Agbar$. Then $\Xbar$ is the closure of the image of 
the map $\prodl_\Ag (j_\Ag \circ\Xi_\Ag) :X\to \prodl_\Ag \Agbar$. (The
product is taken with the Tikhonov topology.) In
particular two sequences $x_n$, $y_n$ have the same limit in $\Xbar$ if and
only if $\lim \Xi_\Ag(x_n)=\lim \Xi_\Ag(y_n)$ for any $\Ag\supset
\Delta_0$.
However, by \cite {BT1} 7.4.20(ii) the map $\Xi_\Ag$ does not increase
distances. Hence the statement is reduced to the analogous statement about
 sequences in the  compactified affine space $\Agbar$, which is obvious. \epf

The statement a) follows from the Fact directly.
 b) is clear. 
 \square

\subsubsection{} For $\Mmod\in \dM$ we now define a presheaf
 $\widetilde\L(\Mmod)$ 
 on $Y$ by $\Gamma(U,  \widetilde\L(\Mmod)):=\lim\limits_{\longleftarrow Z}
 Sm(\Mmod)/(\Mmod_Z \cap
Sm(\Mmod))$ where $Z$ runs over all closed subsets $Z\subset U$. We define
$\L(\Mmod)$ to be the associated sheaf.

From the definition it is clear that the stalk of $\L(\Mmod)$ at a point $z\in Y$ equals
$Sm(\Mmod)/(\Mmod_{\{z\} } \cap Sm(\Mmod))$.

The desired properties of the functor $\L$ follow from the next

 \begin{Thm}\label{main}  a) Let $z\in \overline X - X$ be a point.
Then the functor ${\goth M} \to {\goth M}_{\{z\}}^{\sm}$ on the category $\dM$  is 
exact.

b) Let $P$ denote the stabilizer of the point $z$, then $P$ is a parabolic
subgroup in $G$. Let $r_P$ be  the corresponding Jacquet
functor (that is coinvariants with respect to the unipotent
radical of $P$).  For any ${\goth M}\in \dM$ the subspace 
${\goth M}_{\{z\}}^{\sm}\subset Sm({\goth M})$ contains the kernel of the
projection $ Sm({\goth M}) \to r_P(Sm({\goth M}))$.
Furthermore, if ${\goth M}$ lies in the category $\dR$ then the sequence
$$0 \to {\goth M}_{\{z\}}^{\sm} \to Sm ({\goth M}) \to r_P(Sm({\goth M}))\to
0$$ is exact. 
\end{Thm}

{\sf Proof} occupies the rest of this section.

\subsubsection{}
We abbreviate ${\goth F} = {\goth F} _{\{z\}}$. Our first aim is to 
get a more tractable description of $\goth F$.

Fix a Levi decomposition 
 $P=L\cdot N$, and a maximal split torus  $T\subset L$.
 Let $A$ be the center of $L$,
and $A^s$ be its maximal split subtorus.

We have vector spaces 
 ${\goth a} = {\Hom}(\underline {F^\times},\underline
A)\otimes \RE
= (A/A^c)\otimes \RE = (L/L^c)\otimes \RE$ where $A^c\subset A$
is the  maximal compact subgroup and $L^c\subset L$ is the subgroup
generated by all compact subgroups;
and  ${\goth a} _T={\Hom}(\underline{F^\times},\underline
 T)\otimes \RE =
(T/T^c)\otimes \RE$ where $T^c\subset T$ is the maximal compact
subgroup. The  imbedding ${\goth a} \hookrightarrow {\goth a} _T$
is obtained from the imbedding ${\Hom}(\underline{F^\times},\underline
 A)={\Hom}(\underline{F^\times}, \underline{A_s})
\hookrightarrow {\Hom}(\underline{F^\times}, \underline T)$. 

 (Here we use the notation $\underline
H$ for an algebraic group over $F$, and $H$ for its group of $F$-points.) 
We identify Hom$(\underline{F^\times},\underline A)$
 and $L/L^c$ with subgroups in $\goth a$
and the groups   Hom$(\underline L, \underline{F^\times})$, Hom$(
\underline A,\underline{F^\times})$ with subgroups in the dual
space ${\goth a} ^*$.

We  fix a uniformizer $\goth p$ in  $F^\times$, and get  cocompact 
imbeddings $\Hom (\underline{F^\times} , \underline T)
\hookrightarrow T$,
$\Hom (\underline{F^\times} , \underline A)
\hookrightarrow A$, given by  $\psi \mapsto \psi({\goth p})$.
We denote the image of latter imbedding by $\Lambda$. 

Let $\pi$ be the projection $L\to L/L^c$.

By $\Lambda^+$ we denote the set of $\lambda \in \Lambda$
such that $\pi(\lambda)\in {\goth a} ^+$ where ${\goth a}^+$
is the positive chamber,  ${\goth a} ^+ = \{
a\in {\goth a} | r(a)<0$ if $r$ is a root of $A$ in $N\}$.

\smallskip

We fix notations and recall some basic facts about the geometry of $X$
(see \cite{BT1}, \cite {BT2}, \cite{Landv} for more details). 

 Recall that the set of maximal
split tori is in bijection with apartments in $X$. Let $\Ag$ be the
apartment corresponding to $T$, and let $\Agbar\subset \Xbar$
be the closure of $\Ag$.

Then $\Ag$  is  an affine space with underlying
vector space ${\goth a} _T$. The torus $T$ acts on $\goth A$ through
$T/T^c$ acting on the affine space by translations.

The closure $\overline {\goth A }$ of $\goth A$ in $\overline X$ 
is the standard semisphere compactification of the affine space
$\goth A$. In particular $\overline {\goth A }-\goth A $ 
is  the set of rays $({\goth a} _T - \{0\}) / \RE _+^*$. 
By our assumption 
the point $z$ lies in  $\overline {\goth A }-\goth A $ and corresponds
to a ray $\rho _z \subset {\goth a} \subset {\goth a} _T$. 
We have $\rho _z \subset \a^+$.

By a cone in a vector space we will mean a convex centered
cone. 

\begin{Def} a) A $z$-good cone is an open 
 cone $C\subset {\goth a}$
such that $\rho _z \subset C $.

For a $z$-good cone set $\Lambda_C = \Lambda\cap \pi^{-1}
( \overline C  \cap(L/L^c))$ and  $\Lambda_C^+ = \Lambda\cap \pi^{-1}
( C  \cap(L/L^c)).$

b) A $z$-good semigroup is an open compactly generated semigroup 
 $M\subset L$ such that 

$\;$ i) $M$ contains $\Lambda_C$ for a $z$-good cone $C$.

$\;$ ii) We have $L=\Lambda\cdot M$.
% The projection $L \to L/\Lambda $ restricted to $M$ is surjective.
\end{Def}

It is easy to see that an intersection of two  $z$-good semigroups
is again a $z$-good semigroup.

 \begin{Prop}\label{filterF} The filter  $\goth F$ is generated by the sets
$$ K' \cdot M \cdot a \cdot K_p,$$
where $K'$ is an open compact subgroup in $G$,
 $M$ is a $z$-good semigroup,     $a\in \Lambda $ and
   $p \in \goth A$.
\end{Prop}

We need some more notations and basic facts about $\Xbar$.

Recall that  $d(\cdot,\cdot)$ denotes the canonical metric on $X$. 

For $x \in X$ denote by $K_x$ the stabilizer of $x$.
It is an open compact
 subgroup in $G$.

\begin{Lem}\label{c_to_z}
 a) For any points $c\in X, \; z \in \overline {X} - X$
there exists a unique geodesic $[c,z)$ connecting $c$ and $z$.
(That is  $[c,z)$ is the unique geodesic ray starting at $c$ and
having $z$ as its  limit point in $\overline {X} - X$.)
For $c\in {\goth A}$ the ray $[c,z)$ coincides with $\rho_z +c$.

b) Assume that $c_1,c_2 \in X,\;z\in  \overline X - X$
and let $c_i ^t$ be the unique point on  $[c_i,z)$ at  the distance
$t$ from $c_i$ for $i=1,2$. Then we have
 \begin{equation}
\label{ineq}
d(c_1^t,c_2^t)\leq d(c_1,c_2)
  \end{equation}
\end{Lem}

\proof a) follows from \cite{Schn1}, Proposition on p. 166. More precisely,
 in {\it
loc. cit.} it is proven that there exists an apartment $\Ag$ containing
$c$ and $z$ in its closure, and that the ray $[c,z)_\Ag \subset \Ag$ of
 direction given by $z$ with 
vertex at $c$ does not depend on the choice of $\Ag$.

This gives existence; to show uniqueness it remains only to check that
any geodesic ray $\rho$ lies in some apartment. Let $\Ag$ be any apartment
containing the vertex $c$ of the ray $\rho$.
By \cite{BT1} 2.3.4 the set $\C_t:=\{g\in K_c \, |g(c'_t)\in \Ag\}$
is nonempty for all $t$. Each of the sets $\C_t$ is  closed (and  open)
since it is  
a finite union of cosets of an open compact subgroup $K_c\cap K_{c'_t}$; also 
$\C_{t_1}\subset \C_{t_2}$ if $t_1>t_2$ because
$\Ag$ contains the geodesic connecting any two points of $\Ag$,
 \cite{BT1} 7.4.20(iii).
 Since $K_c$ is compact we have $\capl_t \C_t\ne \emptyset$. Now we can take
$\Ag':=g^{-1}(\Ag)$ for any $g\in \capl_t \C_t$.

Let us prove b). Let $x_1=c_1, x_2,\dots,x_n=c_2$ be a sequence of points 
on $[c_1,c_2]$ such that $[c_1,c_2)=\cup [x_i,x_{i+1})$, and $x_i,x_{i+1}$
lie in the closure of an open polysimplex $\Delta_i$. 

By \cite{BT1} 7.4.18(ii) there exists an apartment $\Ag_i$ such that 
$\Ag_i\supset \Delta_i$, and $\Agbar_i\owns z$. Thus 
$[x_i,z)$, $[x_{i+1},z)$ are parallel rays in the affine space $\Ag_i$,
so $d(x_i^t,x_{i+1}^t)=d(x_i,x_{i+1})$ for all $t$. We finally get
$$d(c_1^t,c_2^t)\leq \suml_{j=1}^{n-1} d(x_j^t,x_{j+1}^t)=
\suml_{j=1}^{n-1} d(x_j,x_{j+1})=d(c_1,c_2).\ \epf$$

\subsubsection{}\label{mu} 
For any $x\in \goth A$ the action map $\mu : K_x \times \overline {\goth A}
\to \overline X$ is proper, 
surjective and the topology on $\overline X$ is the 
quotient topology with respect to this map \cite{BSer}, 5.4.1.
Hence the filter of neighborhoods of $z\in \Xbar-X$ is generated by images
of neighborhoods of $\mu^{-1}(z)$ in $K_x\times \Agbar$.

Fix a point $x\in \Ag$ which lies inside an open polysimplex $\Delta_0\subset
\Ag$. Then $K_x(z)\cap \Agbar=\{z\}$. Indeed, if, on the contrary,
$z\ne k(z)\in \Agbar$ for some $k\in K_x$, then $k\([x,z)\)= [x,k(z'))$ is a 
ray in $\Ag$ different from $[x,z)$; however, since $k$ fixes the
neighborhood $\Delta_0$ of $x$ in $\Ag$, the intersection of the  two rays
with this 
neighborhood coincide, which is impossible. 

Thus $\mu^{-1}(z)=Stab_{K_x}(z)\times\{z\}=(P\cap K_x)\times\{z\}$.

Let   $P^-=L\cdot N^-$ be the  parabolic opposite to $P$.
Since $x$ lies in an open polysimplex of $\Ag$ we have the following
 triangular decomposition:
\begin{equation}
\label{triang}
K_x = K_x^-\cdot K_x^0\cdot K_x^+  
\end{equation}
where $K_x^-=K_x \cap N^-,\;K_x^0=  K_x \cap L,\;  K_x^+=
 K_x \cap N^+$ (for $P$ minimal this follows from \cite{BT1} 7.1.4(2),
and the general case is a corollary). 

The fundamental system of neighborhoods of $K_x\cap P= K_x^0\cdot K_x^+$
 is formed by  subgroups of the form $O=O^- \cdot K_x ^0 \cdot K_x^+$
where $O^-$ is a small open compact subgroup in $N^-$.

The filter of neighborhoods of 
$z$ in $\Agbar$ is generated by closures of the cones
  $\widetilde C+t$ where $\widetilde 
C\subset {\goth a} _T$ is an open cone containing $\rho _z$, and
$t\in \goth A$ is a point.

Since $K_x$ and $\Agbar$ are compact it follows that the filter of
neighborhoods of $(K_x\cap P)\times\{z\}$ is generated by the sets 
$O\times (\widetilde C+t)$ in the above notations. Thus 
the sets $O(\widetilde C+t) $ generate the filter of neighborhoods of $z$
in $\Xbar$. 

Further, it is enough to take only such $\widetilde C$, $t$ that 
$\alpha|_{\widetilde C}\leq 0$ if $\alpha$ is a root of $T$ in $N$,
and $t\in x+\widetilde C$. Then we have $\widetilde C+t\subset 
\widetilde C+x\subset \cupl_{z',P_{z'}\subseteq P} [x,z')$.
Since $K_x^+\subset K_x\cap P'$ for any $P'\subset P$ we see that $K_x^+$
fixes $\widetilde C+t$ pointwise. 

Thus the system of neighborhoods of $z$ in $X$ is generated
by $O(\widetilde C+t)=O^-\cdot K_x^0 \cdot K_x^+(\widetilde C+t)=
O^-\cdot K_x^0(\widetilde C+t)$ for $O,\widetilde C,t$ as above;
we will call $O(\widetilde C+t)$ a {\it standard} neighborhood of $z$.

\subsubsection{} \proof of Proposition \ref{filterF}.
Fix  a $z$-good semigroup $M$, an element $a\in \Lambda_C$ for a $z$-good
cone $C$, and an open subgroup $K'\subset G$.
First we want to find a neighborhood $U$ of $z$ in $X$ such that
$a_p^{-1}(U)\subset K' \cdot M\cdot a \cdot K_p$, i.e. 
$U\subset  K' \cdot M\cdot a (p)$.

\smallskip

\noindent {\bf Claim} {\it Let $\pi_T $ denote  the projection
from $T$ to $T/T^c \subset {\goth a} _T$. Then  the semigroup
 $\pi _T(M\cap T)$ contains
$(T/T^c)\cap \widetilde C$ for 
some open 
cone $\widetilde C\subset {\goth a} _T$ containing $\rho _z$.}

\prf 
Notice that $\pi_T(M\cap T)+A/A^c=T/T^c$. 
 Moreover, $(\pi _T(M\cap T))\cap {\goth a}$ contains
$(T/T^c)\cap C = (A/A^c)\cap C$ for a $z$-good cone $C\subset {\goth a}$.
This yields the claim  by the following elementary argument.

Induction in $\dim(\a_T/\a)$  reduces the statement to the
 situation when $\dim(\a_T/\a)=1$; in this case one can argue as follows.
By the second sentence of the proof $C$ contains a
 basis $v_1,\dots ,v_r$ of $A/A^c$ such that $\sum a_i v_i\in \rho$ for
some $a_i>0$. Let $v$ be an element of $\pi_T(M\cap T)$
which does not lie in $\a$, and take $v'\in (-v+(A/A^c))
 \cap \pi_T(M\cap T)$; denote the sublattice generated by $v,v_1,\dots,v_r$
 by $\Lambda'$. 

The subset $\{a_0v+\sum a_iv_i\,|\,a_i \geq 0\}\cup
\{a_0v'+\sum a_iv_i\,|\,a_i \geq 0\}$ obviously contains an open cone
$C'$ such that 
$C'\supset \rho$. Any point of $\Lambda'\cap C'$ has nonnegative integral
coordinates with respect to one of the bases $(v,v_1,\dots,v_r)$ or $(v',
v_1,\dots,v_r)$, hence lies in $\pi_T(M\cap T)$. 

On the other hand, $\Lambda'$ is a sublattice of finite index in
$T/T^c$; let $x_1,\dots ,x_n\in T/T^c$ be a set of representatives for
the $\Lambda'$ cosets. By the first sentence of the proof we can assume
that $x_i\in \pi_T(M\cap T)$.  

It is clear that for a narrow enough subcone $C\subset C'$, $C\supset
\rho$, we have $y\in C\cap \pi_T(M\cap T),
 y\not \in \rho \Rightarrow y-x_i\in C'$ for all $i$. 
For some $i$ we have $y-x_i\in \Lambda'$,
  so $y-x_i\in \Lambda'\cap C'\subset 
\pi_T(M\cap T)$, and hence $y=x_i+(y-x_i)\in \pi_T(M\cap T)$.
If $y\in (T/T^c)\cap \rho$, then $y\in \pi_T(M\cap T)$ by the second
sentence of the proof. 
\epf

\smallskip

We can assume that $\alpha|_{\widetilde C}\leq 0$
if $\alpha$ is a root of $T$ in $N$. %%%%%%%%%% \subset {\goth a} _T ^+$.

Also we can find $b\in \Lambda_C$ such that $b+a+p\in x+\widetilde C$, so that 
 $b+a+p+\widetilde C$ is point-wise fixed by $K_x^+$. We can furthermore
 assume that $M\supset b\cdot K_x^0$. 

We now take $O=O^- \cdot K_x ^0 \cdot K_x^+$ as above such that
$O^- \subset K'$. Then $K'\cdot M \cdot a (p) \supset O^- \cdot K_x^0
 \cdot ab \cdot (\widetilde C\cap (T/T^c)) (p)$. So $U= O^- \cdot K_x^0
 (a+b+ \widetilde C +p)=O (a+b+ \widetilde C +p)$ is the desired neighborhood.

It remains to show that for any neighborhood $U$ of $z$ in
$\overline {X}$ there exist $M,a,K'$ as above such that
$K'\cdot M\cdot a \cdot K_p \subset a_p^{-1}(U)$ i.e.
$K'\cdot M\cdot a (p) \subset U$. We will do it with the help of another
auxilary geometric statement.

For 
$x \in X,\; r\in \RE$ write $B(c,r)=\{y\in X|d(y,c)<r\}$
for the ball of radius $r$ centered at $c$.
 For a point $y\in \goth A$ and $\lambda \in \RE^{>0}$ we define 
a {\it sector} $S(y,\lambda)$ by $S(y,\lambda)= \cup_{t\in \RE ^{>0}}
B([y,z)_t , \lambda t)$ (notations of \ref{c_to_z}b)).

\begin{Lem} a) For each neighborhood $U$ of $z$ in $X$ 
there exists a sector $S=S(y,\lambda)$ such that $S\subset U$.

b)Let $S\subset X$ be a sector.
 For any compact set $\C\subset L$ and a $z$-good
cone $C$ there exists $a \in \Lambda_C$ such that $ag(S)\subseteq S$
 for all $g\in \C$.
\end{Lem}

{\sf Proof} a) We can assume that $U=O(C+t)$ is a standard neighborhood.

It is obvious that $S(q,\lambda)\cap {\goth A} \subset U$ for $q\in U\cap \goth A$
and $\lambda $ small enough. It  will be convenient to take $q$ inside
of an open polysimplex $\Delta_q$ of the canonical polytriangulation.

By \cite{BT1} 2.3.4, 7.4.19 any point
$s\in X$ is $K_q$-conjugate to a unique point in $\goth A$. Moreover,
by \cite{BT1} 7.4.20 ii),
for any  $k \in K_q$ and $x_1, x_2 \in \goth A$ we have $d(k(x_1), x_2)\geq
d(x_1,x_2)$. Hence
$S=S(t,\lambda) \subset K_q(S\cap {\goth A})$.

We can insert $q$ instead of $x$ in the
 triangular decomposition  ~\eqref {triang}.
Now for  $\lambda $ small enough the cone 
we have $S(q,\lambda)\cap {\goth A} \subset a_T^+ + q$, where $a_T^+=
\{\lambda \in \a_T| (\alpha,\lambda)\leq 0$ if $\alpha$ is a root of $T$ in
 $N$, which implies that  
$S(q,\lambda)\cap{\goth A}$ is point-wise fixed
by $K_q ^+$ (see above). So $S(q,\lambda) \subset K_q(S\cap {\goth A})
= K_q^- \cdot K_q ^0 (S\cap {\goth A})$. It remains to  notice that
if we take e.g. $q = a(x)$ for appropriate  $a \in \Lambda_C$ 
($a$ should be ``deep enough'' in ${\goth a} _T^+$) 
then, besides $S(q,\lambda)\cap \Ag\subset U$,      we will get $K_q^- =
a(K_x^-)a^{-1}\subset O^-$,  $K_q^0 = K_x ^0$.
 Then
 $ S(q,\lambda) \subset  K_q^- \cdot K_q ^0 (S(q,\lambda)\cap {\goth A})
\subset O^- \cdot K_x ^0(S\cap {\goth A})\subset U$.
 a) is proven.

b) Let $\C\subset L$ be a compact subset. Then 
 $d(x,g(x))$ is bounded by a  constant $d_0$ for $g\in \C$.

Set $d_1 =\sup_{t>0} \inf _{a\in \Lambda_C} d(a(x),[x,z)_t)$.
It is clear that $d_1<\infty$.

Take  $s\in S, g \in \C$; let $t_0\in \RE^{>0}$, $a\in \Lambda_C$
 be such that $d(s,[x,z)_{t_0})<\lambda t_0 $, $d([x,z)_{t_0},a(x))
\leq d_1$.
Then we have
\begin{multline*}
%d(s,[x,z)_{t_0})<\lambda t_0 \Longrightarrow 
d(ag(s),[x,z)_{t+t_0})\leq 
 d(ag(s),ag([x,z)_{t_0}))+ d(ag([x,z)_{t_0}),[x,z)_{t+t_0} )< \\
\lambda t_0 +  d(ag(x),[x,z)_t).
\end{multline*}

Here in the second step we used the inequality
$d(ag([x,z)_{t_0}),[x,z)_{t_0+t})$
$\leq  d(ag(x),[x,z)_t)$ which is just ~\eqref {ineq}  applied to the case
$c_1=ag(x)$, 
$c_2=[x,z)_t$. Note that $ag(z)=z$ because $a,g \in L$, hence
$ag([x,z))=[ag(x),z)$.

But $$d(ag(x),[x,z)_t)\leq d(ag(x),a(x))+d(a(x),[x,z)_t)\leq d_0+d_1.$$

Thus we get
$$d(ag(s),[x,z)_{t+t_0})\leq \lambda t_0 +d_0+d_1< \lambda (t+t_0),$$
provided $\lambda t_0 >d_0+d_1$. So $ag(S)\subset S$
and the Lemma is proven. 
 \square

\subsubsection{} We are now ready to finish the proof of Proposition \ref{filterF}.
Assume a neighborhood $U$ of $z$ in $\overline X$ is given.
We are looking for $M,a,K'$ as above such that $K'\cdot M\cdot a(p)
\subset U$.

Let $S=S(y,\lambda)$ be a sector contained in $U$.
There exists
a $z$-good cone $C\subset {\goth a}$ and  $a\in \Lambda_C$ such that
$a+C+p\subset S\cap {\goth A}$.

Let $\C$ be a compact set generating $L$ as a semigroup. By the previous Lemma
we can find $b\in A$ such that $b\cdot \C$ preserves $S$.

We set $M$ to be generated by $\Lambda_C$, $K_y\cap L$ and $b\cdot \C$; the
element $a$ is chosen so that $a+C+p\subset S\cap {\goth A}$.

We can assume $U\cap X=O(C+t)$ is a standard neighborhood,  and take $K'=O$.

Let us check that $K'\cdot M\cdot a(p)\subset U$. 
By the construction $\Lambda_C\cdot a(p)\subset S\cap \Ag$.
Since $L$ commutes with $\Lambda$, and $b\cdot \C$, $K_y\cap L$ preserve
$S$ we have $M\cdot a(p)\subset S$. Hence $ K'\cdot M\cdot a(p)\subset
K'(S)=O(S)\subset U$.

The Proposition is proven.
 \square

\subsubsection {}
We return to the proof of the Theorem.

 For an   open 
compact subset $\C\subset G$
we denote by $e_\C \in \Heck$ the constant measure on $\C$
 of total volume 1. We write $\cdot $ for the convolution of
compactly supported measures on a group.

We can   fix  a small open compact subgroup
$K\subset  K_p$ which is normal in $K_p$
and satisfies the triangular decomposition
\eqref{triang} (see \ref{2.1} of part \ref{2} below for a stronger
statement).   

Set: $L^+=\{g\in L | g (K^-)g^{-1}
\subset K^-\}.$ It is readily seen that the semigroup $L^+$
is $z$-good.

We have a map $\phi : \Heck (L^+, K^0)\to \Heck (G,K)$ sending $h$ to
$e_K \cdot i_*(h) \cdot e_K$ (here $i:L^+\hookrightarrow G$
is the embedding, and $i_*$ denotes the direct image of a measure).
We abbreviate $ \Heck (L^+, K^0)$ to $ \Heck _+$.

\begin{Lem}\label{homo}
 $\phi$ is a homomorphism of algebras.
\end{Lem}

{\sf Proof} see e.g.  \cite{Al}, Theorem 1 (to apply {\it loc. cit.} as
stated $L^+$ has to consist of semisimple elements, i.e. 
to be a minimal Levi, and $N$, $N^-$
must be maximal unipotent subgroups. However these more restrictive 
 conditions are not used in the proof of Theorem 1).   \square

Let ${\goth F} _{\Heck _+}$ be the filter of subspaces in $\Heck _+$
 generated by  $\Heck (a\cdot M,K^0)$ where $M\subset L^+$ is a 
$z$-good semigroup containing $K^0$ and $a\in  \Lambda_C\subset M$.

The Hecke algebra $\Heck (G)$ acts on every smooth $G$-module;
this action is denoted by $h:m \mapsto h(m)$. For any $h\in \Heck (G)$
we denote by $h^*$ the image of $h$ under the involution 
$g\to g^{-1}$.

  \begin{Prop}  \label{algb}
a) Fix ${\goth M}\in \dM$. Assume that
$*({\goth M})$ is generated by its $K$-fixed vectors.
  Then for  $f \in {\goth M}^K$ the equality
$ \lim _ {\goth F } g^{-1}(f)=0$ is equivalent to
$ \lim _ {{\goth F} _{ \Heck _+}} \phi(h)^*(f)=0$.

b) The functor $H_z$ on the category of $ \Heck _+$-modules sending
a module ${\goth M}$  to the space $\{ f\in {\goth M}^*\;|\; \lim _
{{\goth F} _{ \Heck _+}} \langle h(m), f \rangle =0$ for any $m\in {\goth M}\}$
 (where ${\goth M}^*$ is the full linear dual, and the basic field is equipped 
with 
discrete topology)
 is exact. 
\end{Prop}

{\sf Proof}
a) Assume that $ \lim _ {{\goth F} _{ \Heck _+}} \phi(h)^*(f)=0$.
Then for any $m\in *({\goth M})^K$ we can find a $K$-invariant
$z$-good semigroup $M\subset L^+$ and $a\in \Lambda _C$
 such that
\begin{equation}
\label{condition}
h\in \Heck (a\cdot M, K^0)\Longrightarrow \langle m,
\phi (h)^*(f) \rangle =0 
\end{equation}

Hence for any finite set $m_1,..,m_l\in  *({\goth M})^K$ there is
some $M,a$ as above such that ~\eqref{condition} is satisfied
for $m\in \{m_1,..,m_l\}$. Fix $m\in  *({\goth M})^K$ and take
the finite set to be $K_p(m)$ (recall that we assumed that 
$K_p$ normalizes $K$, hence preserves $ *({\goth M})^K$).
Then we get $M,a$ such that 

$ h\in \Heck (a\cdot M)\Longrightarrow \langle k(m),
\phi (h)^*(f) \rangle =0$ for $k\in K_p$. This means
that   $ \langle m,
g^{-1}(f) \rangle =0$ for $g\in K\cdot M \cdot a \cdot K_p$
hence $ \lim _ {\goth F }\langle m,
g^{-1}(f) \rangle =0$.

It remains to notice that
 $\{  m\in *({\goth M})\,|\;\lim_{\goth F} \langle g(m),f \rangle =0\}$
is a $G$-submodule in $*(\goth M)$ (because ${\goth F}$ is right-translation
invariant). Thus if  $ *({\goth M})$
is generated by its $K$-fixed vectors, then 
  $ \lim _ {\goth F }\langle m,
g^{-1}(f) \rangle =0$ for any $m \in *({\goth M})$, i.e.
 $ \lim _ {\goth F } g^{-1}(f) =0$.

Conversely, if $ \lim _ {\goth F } g^{-1}(f)=0$ then for
any $m$ there is  a $z$-good
semigroup $M\subset L^+$ and $a\in \Lambda \cap L^+$ such that 
 $ \langle m,
g^{-1}(f) \rangle =0$ for $g\in  M \cdot a $. We can assume that
 $a \in \Lambda_C$. Any  $z$-good
semigroup contains one normalized by $K^0$, so we can assume also
that
$M$ is  normalized by $K^0$. Then $M\cdot K^0$ is again a  $z$-good
semigroup and $ \langle m,
\phi (h)^*(f) \rangle =0$ for $ h\in \Heck (a\cdot M\cdot K_0)$
provided $m$ is $K$-invariant.  a) is proved.

\subsubsection {} To prove b) we need an auxilary definition.

The homomorphism $\pi|_{L^+}: L^+ \to
 L/L^c$ defines an $L/L^c$-grading on $\Heck _+$.

\begin{Def} \label{zgoodalg} A $z$-good subalgebra in $\Heck_+$ is a graded
subalgebra ${\cal A}=\oplus_{\lambda\in L/L^c} \A_\lambda$ in $\Heck ^+$
 such that 

i) For some $z$-good cone $C$ all elements of the form $ e_{K^0\cdot a}$
belong to $\A$ for $a\in \Lambda_C$.

ii) The map $pr_*:\A\to \H(L/\Lambda, K^0)$ is surjective.
\end{Def}

\begin{Sublem}\label{sublem} Let $\A$ be a $z$-good semigroup, $C$ be a
$z$-good cone 
such that $e_{\eta\cdot K^0}\in \A$ for $\eta\in \Lambda_C$. Then for any
$h\in \H(L,K^0)$ there exists $\nu\in \Lambda_C$ such that
$\nu\cdot h \in \A$. 
\end{Sublem}
\proof We can assume that $h$ is supported in one $L^c$-coset $\lambda$.
 By condition ii) of the definition we have $pr_*(h)=\sum pr_*(h_i)$,
where $h_i\in \A_{\lambda_i}$. Since $\H(L/\Lambda)$ is graded by the
 finite group  $L/(L^c\cdot A)$, we can further assume that 
$\lambda_i=\lambda \mod \Lambda$. Then we get: $h=\sum
 (\lambda-\lambda_i)(h_i)$. Clearly for ``large enough'' $\nu \in
 \Lambda_C$
we have $\nu+\lambda-\lambda_i\in \Lambda_C$, so $(\nu)h=
\sum (\nu+\lambda-\lambda_i)(h_i)\in \A$. \epf

\begin{Lem}\label{A_i} a) The filter ${\goth F}_{\H_+}$ is generated
by the subspaces $a\cdot {\cal A}$, where $\cal A$ runs over the set
of $z$-good subalgebras and $a$ runs over $\Lambda$.

b) There exist a set of $z$-good subalgebras $\A _i$ and elements
$a_i \in \Lambda $, $e_{a_i\cdot K^0 }\in \A _i$ such that:

$\;$ i) ${\goth    F}_{\Heck ^+}$ is generated by  $a_i\cdot \A_i$

$\;$ ii) Each of 
the algebras $\A_i$ is finite over its center, which has finite type.
\end{Lem}

{\sf Proof} a) It is obvious that if $M$ is a $z$-good semigroup, then
$\H(M,K^0)$  is a $z$-good subalgebra. 

On the other hand, let $\C$ be a
compact set generating $L^c$ as a semigroup. 
Fix a $z$-good subalgebra
$\A$. Since $\Heck(\C,K^0)$ is finite dimensional Sublemma \ref{sublem}
gives a coweight $\nu\in \Lambda_C$ such that $(\nu)\cdot \H(\C,K^0)
\subset \A$.

Let  $M$ be the semigroup generated by $K^0$, $\C\cdot
 \nu$  and $\Lambda_C$. Obviously $M$ is a $z$-good semigroup, and
 $\H(M,K^0)\subset  \A$. This proves a).

The proof  of  b) consists in constructing 
the Noetherian algebras $\A_i$ explicitly.

 It is known \cite{BDKV} that  $\Heck (L, K^0)$
is finite over its center $Z_L$; and $Z_L$ is a finitely
generated commutative algebra. Let $z_1,..,z_n$ be generators
of $Z_L$, and $h_1,..,h_l$ be the generators of $\Heck (L,K^0)$
as $Z_L$-module. We have: $h_i\cdot h_j = \sum z_{ij}^l h_l$
for $   z_{ij}^l\in Z_L$. 

Fix a rational (with respect to the lattice
$L/L^c \subset {\goth a}$) $z$-good cone $C\subset {\goth a}^+$.
Then the semigroup $\Lambda _C$ is finitely generated.

 We can find
$a \in \Lambda _C^+$ such that $a \cdot z_i$, 
 $a \cdot z_{ij}^l$ and  $a \cdot  h_l$ are supported inside
$L^+$. 

For such  $a,C$ define $\A_{C,a}$ to be the subalgebra generated by
$e_{\eta \cdot K^0}$ for $\eta \in \Lambda _C$
and  $a \cdot z_i$, 
 $a \cdot z_{ij}^l$,  $a \cdot  h_l$.

Then  the subalgebra  ${\cal Z}_\A
\subset \A_{C,a}$ generated by
$e_{\eta \cdot K^0}$ for $\eta \in \Lambda _C$
and by $a \cdot z_i$,  $a \cdot z_{ij}^l$ is central. Also  $\A_{C,a}$
is finite  over ${\cal Z}_\A$; namely, it is clearly generated by
  $a \cdot  h_l$ as a ${\cal Z}_\A$-module. Since ${\cal Z}_\A$
is a finitely generated commutative algebra, it is 
 Noetherian, so the center of $\A_{C,a}$ is a finite ${\cal Z}_\A$-module,
hence itself is a finitely generated commutative algebra.

To finish the proof it is enough, in view of a),  to check that every $z$-good
subalgebra contains a subalgebra of the form $\A_{C,a}$.
For this let $\A$  be a $z$-good algebra, and $C$ be a $z$-good cone
such that $e_{\lambda \cdot K^0}\in \A$ for $\lambda \in \Lambda_C$. 
Then Sublemma \ref{sublem} guarantees  existence of $a\in \Lambda_C$ such that 
 $a \cdot z_i$, 
 $a \cdot z_{ij}^l$,  $a \cdot  h_l\in \A$, i.e. such that $\A_{C,a}
\subset \A$. \epf

Let us finish the proof of the Proposition.

Fix $\A\in\{\A_i\}$. We can pick $\A_i$, $a_i$ in the Lemma so that
 $\A_i\subset \A$.

The filter of subspaces in $\A$ generated by $a_i\A_i$
is invariant under both left and right multiplication by an element
$h\in \A$, in the sense that 
for any $i$ and $h\in \A$ there exists $j$ such that
$h\cdot a_i\A_i\subset a_j\A_j$,
$ a_i\A_i\cdot h\subset a_j\A_j$.
Thus the space $\dirlim(\A/a_i\A_i)^*$  carries
a natural $\A$-module structure. Moreover,
for any $\Mmod\in \tM$ we have
                                                \begin{multline}
H_z(\Mmod)=
\{f\in \Mmod^*\ |\ \forall m\in \Mmod\; \dirlim_{\F_{\H_+}}\<h(m),f\>
=0\}=                                                             \\
\{f\in \Mmod^*\ |\ \forall m\in \Mmod\  \exists \ z\operatorname{-good}
\ \A,\,a\in
\Lambda \ \forall h\in a\A :\ \<h(m),f\>=0\
\}                =\\
\{f\in \Mmod^*\ |\ \forall m\in \Mmod\  \exists {i}\ \forall h\in a_i\A_i :
  \<h(m),f\>=0\
\}
= \\
\Hom(\Mmod, \dirlim (\A/a_i\A_i)^*)
                                                 \end{multline}

So to prove the Proposition we have only to show that $\dirlim(\A/a_i\A_i)^*$
is an injective $\A$-module.

It is not hard to see that for fixed $i$ the filter of subspaces in $\A$
 generated
by $a_i^n\A_i$, $n\in \Zet^{>0}$ is also invariant under left and right
multiplication by $h \in \A$, thus we have an $\A$-module $\A_{z,i}^*
\overset{\fed}{=} \dirlim_n(\A/a_i^n\A_i)^*$.

Notice that $\dirlim(\A/a_i\A_i)^*=\dirlim _i(\dirlim_n (\A/a_i^n\A_i)^*)
=\dirlim _i \A_{z,i}^*$
(in fact we can assume that $\forall i,\,n\, \exists j \ |\ \A_i=\A_j,\, a_j=a_i^n$).

Since $\A$ is Noetherian the direct limit of injective $\A$-modules 
is injective.
Hence we have only to show that $\A_{z,i}^*$ is an injective $\A$-module.
For any $\A$-module $M$ we have
\begin{multline}
  \Hom (M,\A_{z,i}^*)=\{f\in M^*\ |\ \forall m\in M\; \exists n\in \Zet^{>0}
:\ \<f,a_i^n \A_i(m)\>=0\}\\
=\Hom_{\A_i}(M,\dirlim _n(\A_i/a_i^n\A_i)^*).
\end{multline}

So we need only to show that $\dirlim _n  ( \A_ i/a_i^n\A_i)^*$
is an injective $\A_i$-module.

 Since $\A_i$ is Noetherian, this is equivalent to saying that the functor
$\Hom_{ \A_i}(\,\underline {\;\;}\, ,\dirlim _n  (\A_i /a^n_i\A_i)^*)$
 is exact on the category of {\it
finitely generated}  $\A_i$-modules.

On this category the functor coincides with the functor
$M\to \dirlim _n (M/a^n)^*$, and
 its exactness is
guaranteed by the standard Artin-Rees lemma because of \ref{A_i}b), 
property ii).

The Proposition is proved.  \square

\subsubsection{} We are ready at last to finish the proof of Theorem
\ref{main}. 

Proof of a): The statement is equivalent to saying that 
the functor $M\to (*(M))_{\{z\}}^\sm$ on the category $\M$ is exact.
Let $0\to M_1 \to M_2 \to M_3 \to 0$ be a short exact sequence in $\M$.
Let $K$ be as above (a compact open subgroup with triangular
decomposition ~\eqref{triang}), and assume that $M_i$ is generated
by its $K$-fixed vectors for $i=1,2,3$. Then by Proposition ~\ref{algb}a):
$$ (*(M_i))_{\{z\}}^K=\{f\in (M_i^K)^*|\lim_{{\goth F}_{\Heck _+}}
\phi (h)^*(f)=0\}.$$
 Hence the sequence
$$0\to  (*(M_3))_{\{z\}}^K\to  (*(M_2))_{\{z\}}^K \to
 (*(M_1))_{\{z\}}^K \to 0$$
is exact. The statement a) of \ref{main} follows.

\hfill

To prove b) it is convenient to use the following

\begin{Lem} Let $M\in \tM$ be a module, $K$ as above. For an element 
$m\in M^K$ the following are equivalent:

$\;$ i) $m$ lies in the kernel of the projection $M\to r_P(M)$

$\;$ ii) $e_{Ka^{-1}K}(m)=0$ for some $a \in \Lambda ^+$.
\end{Lem}

{\sf Proof} is well-known; see e.g. \cite{BZ}, Lemma 3.22 for the case 
$G=GL(n)$, or \cite{Be}, Claim 15.4 for the (absolutely analogous) general
case.   \square 

Now for ${\goth M} \in \dM$ assume
 that $m\in ({\goth M})^K$ lies in the kernel of the projection
$Sm({\goth M})\to r_P (Sm({\goth M}))$. Then for
some $a \in \Lambda ^+$ we have: 

$\phi (e_{aK^0})^*(m)=
e_{Ka^{-1}K}(m)=0\Longrightarrow \phi (e_{aK^0}\cdot \Heck_+)^*(m)=0$

i.e. $m\in {\goth M}_{\{z\}}$ by Proposition ~\ref{algb}a).

Assume now that ${\goth M} \in \dR$, and $m \in \Mmod^K_{\{z\}}$.
For $v \in *(\Mmod)^K$ we have $\lim _{{\goth F} _{\Heck _+}}
\langle \phi(h)^*(m), v \rangle =0$ (by Proposition ~\ref{algb}a))
hence there exists $a\in \Lambda ^+$ such that 
$\langle e_{Ka^{-1}K}(m), v \rangle =\langle \phi(e_{aK^0})^*(m), v \rangle 
=0$. Since $*({\goth M})^K$ is finite dimensional we can find
one $a\in \Lambda ^+$ such that $\langle e_{Ka^{-1}K}(m), v \rangle
=0$ for all $v \in *({\goth M})^K$. But this means that $ e_{Ka^{-1}K}(m)
=0$. By the last Lemma this finishes the proof of Theorem \ref{main}.
 \square

\begin{Thm}\label{cormain}
 For $M\in \M$ let $L(M)\in Sh_G(Y)$ be the unique quotient
of the constant sheaf $\underline{M}$ such that for a parabolic $P$
and $p\in Y$ with  $Stab(p)=P$ (i.e. $p\in \Delta_P$)
 the stalk of $L(M)$ at $p$ is the Jacquet functor $r_P(M)$. 

 The functor $\L:\dM \to Sh_G(Y)$ satisfies:

i) $\L \circ \Gamma \cong i^*\circ j^G_*$

ii) $\L$ is exact

iii) We have a functorial surjection $L(Sm({ M}))\to \L({ M})$;
it is an isomorphism provided ${ M} \in \dR$.
\end{Thm}
\prf Property i)  follows from the  Lemma \ref{premain}.
Properties ii) and iii) follow from  Theorem \ref{main}.
\epf

\subsection{Application to representation theory: dualities}
Recall the definition of {\it Deligne-Lusztig} duality (originally defined
in the context of  groups over finite field in \cite{DL}). 
For  two parabolic subgroups $P_2\subset P_1 \subset G$ we have the canonical
morphism $d_{P_1,P_2}:
i_{P_1}\circ r_{P_1} \to i_{P_2}\circ r_{P_2}$. 
For any $M\in \M$ define the complex $DL(M)$ as follows. Set $(DL(M))^i= 
\bigoplus \limits _{P\in S_i}i_{P}\circ r_{P}(M) $ where $S_i$ is  the set of
standard parabolics of corank $i$; the differential $\displaystyle
d_i:=\sum
\limits_{P_1\in S_i, P_2\in S_{i+1}} \pm d_{P_1,P_2}$ (definition of
signs is straightforward and can be found in \cite{DL}). Obviously
$DL$ extends to an exact functor $Kom(\M)\to Kom(\M)$ and defines a functor
on the derived category.

We will need the following standard Lemma

\begin{Lem}\label{subcat}
$D^+(\R)$ (respectively $D^b(\R)$)
 is equivalent to the full subcategory of $D^+(\tM)$ (respectively $D^b(\tM)$)
consisting of complexes with admissible cohomology.
\end{Lem}  

\prf Let $\tR\subset \tM$ be the full subcategory consisting of direct
limits of admissible modules. 
 By \cite{BDKV}  $\Heck$ is finite over its center $\Ze$ which is locally of
finite type. Let $I\subset \Ze$ be an ideal of finite codimension. Then 
the module $\H/I^n$ is admissible for every $n$, and by the Artin-Rees
Lemma the functor $M\to \Hom_\H(M,\dirlim_n (\H/I^n)\galochka)=\Hom_\Ze(M,
\dirlim_n (\Ze/I^n)\galochka)$ is exact, i.e. the module $\check{\H}_I:=\dirlim_n
(\H/I^n)\galochka$ is injective. It is clear that any module in $\tR$ injects
into $\check{\H}_I$ for some $I$. Hence by \cite{H}, Proposition I.4.8 the
tautological functors $D^b(\tR)\to D^b(\tM)$, $D^+(\tR)\to D^+(\tM)$
induce equivalences with the full subcategories of objects, whose
cohomology lie in $\tR$.

To finish the argument notice that for any complex $C \in Kom^b(\tR)$
(respectively  $C^\cdot \in Kom^+(\tR)$) it is easy to represent $C$ as a
direct limit of its quasiisomorphic subcomplexes $C_i\in Kom^b(\R)$
(respectively $C_i\in Kom^+(\R)$). It follows that $D^{b,+}(\R)$ is
equivalent to the full subcategory of $D^{b,+}(\tR)$ whose objects have
admissible cohomology. \epf

We are now in the position to prove
 
\begin{Thm}\label{dual} 
a) There exists a  canonical morphism  of functors 
from $D^b(\M)$ to $D^b(\tM)$
\begin{equation}\label{d0}
 DL\circ \galochka\to D_h.
\end{equation}

b) Restriction of \eqref{d0} to the full subcategory 
 $D^b(\R)\subset D^b(\tM)$ gives an isomorphism of functors from $D^b(\R)$
to itself
\begin{equation}\label{d}
D_h\cong DL\circ \galochka
\end{equation}
\end{Thm}

\begin{Rem} Isomorphism \eqref{d}, as well as its
generalization \eqref{d1}, 
is due to J.~Bernstein (unpublished). An attempt to interpret
it through the geometry of $\Xbar$ was the starting point
for this part of the work.

Recall that another construction of \eqref{d} using compactified
building appears in \cite{Schn1}. 
\end{Rem}

\begin{Rem} Theorem \ref{dual} is the main step in the proof of
Zelevinsky's conjecture \cite{Zel},  which says that the map on the
Grothendieck 
group of the category of admissible representations induced by $DL$  
sends the class of an irreducible representation to the class of an
irreducible representation (up to a sign). The conjecture is proved in
\cite{Be2},  \cite{Schn1}.
\end{Rem}

 \proof of the Theorem.
 For  $M\in \tM$ denote by $\cDL(M)\in Sh_G(\Xbar)$ the
kernel of the canonical surjection 
 $\underline{M}_{\Xbar}\to 
L(M)$.
%%%%%%%%%%, where  $\underline{M}$ is the constant sheaf on
%%%%%%%%%%$\Xbar$ with fiber $M$, carrying  the obvious equivariant
%%%%%%%%%%structure. 
 Then $\cDL$ is an exact functor from $\tM$ to
$Sh_G(\Xbar)$. 

For an object $I\in \I_G$ we have
the canonical surjections ${\sf pr}_I:\underline{Sm(\Gamma(I))}_{\Xbar}\to
j_*^G(I)$ and 
$L(\Gamma(I))\to \L(\Gamma(I))=i^*j^G_*(I)$ (\ref{cormain}(iii)), 
which yield a morphism of functors
$s: \cDL\circ Sm \circ \Gamma (I)\to Ker(j_*^G(I)\to i^*j_*^G(I))=j_!(I)$. 

Clearly $\cDL$, $s$ extend to a functor 
  (respectively, a morphism of functors) on bounded derived categories; we
  will denote these extensions by the same symbols.  

According to \cite{Schnequiv}, Lemma 1 for any  locally compact $G$-space
 $Z$ the category $Sh_G(Z)$ has enough injectives, so the derived
functor of a left-exact functor on  
$Sh_G(Z)$ is defined. Suppose further that $Z$ is compact; then one
checks immediately that for $\F\in Sh_G(Z)$ the $G$-module $\Gamma(\F)$ is
smooth. Thus in this case we have the derived functor $R\Gamma:D^+(Sh_G(Z))
\to \tM$, which commutes with the forgetful functors $For_{sh}:Sh_G(Z)\to
 Sh(Z)$,  
$For_{mod}:\tM\to Vect_k$. By Corollary 3 in \cite{Schnequiv} the forgetful
 functor 
$Sh_G(X)\to Sh(X)$ sends injective equivariant sheaves into $c$-soft
 sheaves; if $Z$ is 
compact they are  adapted to $\Gamma$,
hence we have $For_{mod}\circ R\Gamma \cong R\Gamma\circ For_{sh}$.
%%%%%%%%%% In particular we see that $R\Gamma$ has
%%%%%%%%%% finite  on equivariant sheaves 
%%%%%%%%%%homological dimension, if it does on the usual sheaves. 

We will be interested in the case $Z=\Xbar$; then all of the above
properties are satisfied.

\begin{Lem} \label{cDL}  
$R\Gamma(\cDL(M))$ is canonically isomorphic to $DL(R\Gamma(M))$
for any $M\in D^b(\tM)$.

\end{Lem}

\proof Recall that $\Xbar-X$ is the spherical building of $G$,  that is the
realization of the simplicial topological space, whose space of
nondegenerate $n$-simplices is the disjoint union of partial flag
varieties 
$G/P$, where $\corank(P)=n+1$; the structure maps come from  the canonical
projections $\pi_{P_1,P_2}:G/P_1\to G/P_2$ for $P_1\subset G/P_2$.
The sheaf $L(M)$ is the realization of a simplicial sheaf on this space;
such a sheaf has a canonical resolution by acyclic sheaves.

More precisely, for a set $s$ of simple roots (of the system of
nonmultipliable roots of $G$) let $\P(s)$
denote the 
corresponding conjugacy class of parabolics (so $s$ is the set of simple
roots which appear in the nilradical of $P\in \P(s)$).
 Let $S(s)=\cupl_{P\in\P(s)}$
be the corresponding stratum of $\Xbar-X$, and $U(s)$ be the union of the
stars 
of all simplices $\Delta_P$, $P\in s$. Then we have $U(s_1\cup
s_2)=U(s_1) \cap U(s_2)$; thus $\{U_s\,|\, \#(s)=1\}$ form an open
covering of $\Xbar-X$, and the corresponding  \v{C}ech resolution of $L(M)$
has the form 
\begin{equation}\label{chech1}
0 \to \oplusl_{|s|=1} j_s^*j_{s_*}(L(M))
\to \dots \to  \oplusl_{|s|=\rank(G)} j_s^*j_{s_*}(L(M)) \to 0,
\end{equation}
where $j_s$ denotes the imbedding $U_s\imbed \Xbar$. 
Hence for any $M\in \tM$ the complex 
\begin{equation}\label{chech2}
0\to \underline{M}_{\Xbar}
 \to \oplusl_{|s|=1}  j_s^*j_{s_*}(L(M))\to
 \dots \to  \oplusl_{|s|=\rank(G)}  j_s^*j_{s_*}(L(M)) \to 0
\end{equation}
is a resolution for $\cDL(M)$. Moreover, there exists a canonical
retraction of $U_s$ on $S(s)$, hence
$H^i(U_s,j_s^*(L(M)))=H^i(S(s),L(M)|_{S(s)} )=0$ for $i>0$, and 
$H^0(U_s,j_s^*(L(M)))=H^0(S(s),L(M)|_{S(s)} )=i_P^G\circ r_P^G(M)$.
Of course, $H^i(\underline{M}_{\Xbar})=0$ for $i>0$ and
$H^0(\underline{M}_{\Xbar})= M$ since $\Xbar$ is contractible.
Thus the terms of the complex of global section of \eqref{chech2} are
identified with that of $DL(M)$; it is easy to see that the differentials
are also the same. 

It is also evident how to extend the assignment $M\mapsto $ \eqref{chech2}
to the definition of an exact functor $\CDL$
from $Kom(\tM)$ to the complexes of
$\Gamma$-adapted objects of $Sh_G(\Xbar)$,
 together with a canonical
quasiisomorphism $\cDL\to \CDL$, and a canonical isomorphism 
$\Gamma\circ \CDL\cong DL$. 
\epf

We can now construct the morphism \eqref{d0}.
For $\F\in D^b(\Sh)$ we have functorial morphisms 

\begin{multline}\label{morph}
DL(R\Gamma_c(\F)\galochka) \overset{\ref{ve}a)}{\cong} DL( Sm
\circ R\Gamma(\V\F)) 
\overset{\ref{cDL}}{\cong} R\Gamma(\cDL( Sm \circ R\Gamma(\V\F)))\\
\overset{R\Gamma(s)}{\longrightarrow} R\Gamma(j_!(\V\F))=
R\Gamma_c(\V\F)\overset{\ref{ve}b)}{\cong} D_h
(R\Gamma_c(\F)).
\end{multline}

Hence by the second statement of Proposition \ref{loc} we get the desired
morphism \eqref{d0} on the full subcategory $D^0(\M)$. a) of the Theorem
now follows from the evident

\begin{Claim} \label{claim}
 a) Let $\A$, $\B$ be additive categories, and 
$\cC$ be a full subcategory of $\A$. Let $F_1$,
$F_2:\A\to \B$ be additive functors, and $s:F_1|_\cC\to F_2|_\cC$ be a
transformation. Assume that every object of $\A$ is isomorphic to a direct
summand in an object of $\cC$. Then $s$ extends in a unique way to a
transformation of functors $F_1\to F_2$.

 b) Any object of $D^b(\M)$ is a direct summand in an object
of $D^0(\M)$.
\end{Claim}

\proof Take $A, A'\in Ob(\A)$ such that $A\oplus A'\in Ob(\cC)$. By the
definition of a morphism of functors we have $s_{A'\oplus A''}\circ
(id\oplus 0)=id\oplus 0$, $s_{A'\oplus A''}\circ
(0\oplus id)=0\oplus id$, hence $s$ preserves the direct sum decomposition,
i.e. $s_{A\oplus A'}=s_A\oplus s_A'$. If $A''\in Ob(\A)$ is another object
such that $A\oplus A''\in Ob(\cC)$, then considering the morphism $id\oplus
0:A\oplus A'\to A\oplus A''$ we see that $s_A$ does not depend on the
choice of $A'$. This proves a). 

To prove b) recall that $\M$ has enough projectives and finite
homological dimension (\cite{Be}, Theorem 38, see also e.g.
\cite{Vigneras}, Proposition  37),  
so any object of $D^b(\M)$ is represented by a finite complex of projective
modules. Since standard projectives form a set of projective generators in
$\M$, we
can add to such a complex a complex with zero differential to get a complex
of standard projectives. \epf

For part b) of the Theorem we need another

\begin{Lem}\label{acycl} 
 For $I\in \I_G$ the cohomology  $H^i(\Xbar, j^G_*(I))$ vanish
for $i>0$.
\end{Lem}

\begin{Sublem} Let $\Delta\subset X$ be a polysimplex of the canonical
polytriangulation, $x\in \Delta$ be a point, and $R$ be a smooth  finite
dimensional  representation of $Stab(R)=K_x$.  Let $\iota_x$ be the imbedding
 $\iota_x:G(x)\imbed X$, where $G(x)$ is the $G$-orbit of $x$, and
set $\F_{x,R}=\iota_{x*}\iota_x^*(\F_{\Delta,R})$ (notations of
\ref{Verdier}).  Consider the canonical arrow $\F_{\Delta,R}\to \F_{x,R}$,
and the induced 
morphism $f:j_*^G( \F_{\Delta,R})\to j_*^G(\F_{x,R})$. Then $i^*(f)$
and $R\Gamma(f)$ are isomorphisms.

\end{Sublem}
                                  
\proof of the Sublemma. Let  $z\in \Xbar-X$ be a point.
Let $K_i\subset G$ be compact open subgroups forming  a fundamental system
of neighborhoods of $1$; then there exists a basis  $\{U_i\}$ of the filter of
neigborhoods of $z$ such that $U_i\cap X$ is a $K_i$-invariant closed
simplicial 
subset. Indeed, if $\Ag$ is an apartment containing $z$ in its closure,
then it is clear that one can find a basis $C_i$ of the filter of
neigborhoods of 
$z$ in $\Agbar$ such that $C_i\cap \Ag$ is closed and simplicial; then by
\ref{mu} 
we get the desired basis of the filter setting $U_i=K_i(C_i)$.

Now we see that $\Gamma(U_i,\F_{\Delta,R})\cong \prodl_{\Delta'\in
G(\Delta), \Delta'\subset U_i } R\iso \Gamma(U_i,\F_{x,R})\cong \prodl
_{G(x)\cap U_i} R$. Hence $j^G_*(\F_{\Delta,R})|_z=\dirlim_{i} \Gamma(U_i\cap
X,\F_{\Delta,R})^{K_i}\iso j^G_*(\F_{x,R})=\dirlim_{i} \Gamma(U_i\cap
X,\F_{x,R})^{K_i}|_z$ (here we write $?|_z$ for a stalk of a sheaf $?$ at a
point $z$). This proves the first statement.  

It  implies that $f$ is a surjection of sheaves, and the kernel of $f$ is
(noncanonically) isomorphic to $\underline{R}_{G(\Deltabar-x)}$ extended by
0 to $\Xbar$. Since $H_c^i(\Deltabar-x)=0$ for all $i$ the latter sheaf is
acyclic, so we get the second statement. \epf

\proof of the Lemma. $\F_{x,R}$ is an injective object in $Sh_G(X)$; hence
$j_*^G(\F_{x,R})$ 
is an injective object of $Sh_G(\Xbar)$. Hence by \cite{Schnequiv},
Corollary 3 it is a $c$-soft sheaf;
since  $\Xbar$ is compact this implies
$H^{>0}(\Xbar,j_*^G(\F_{x,R}))=0$. Thus Lemma follows from the Sublemma. \epf

\subsubsection{} \proof of Theorem \ref{dual}b). 
To prove b) we must check the following. Suppose that $\F\in
D^b(\Sh)$ is such that $R\Gamma_c(\F)\cong M_1\oplus M_2$ where $M_1\in
D^b(\R)\subset D^b(\M)$, $M_2\in D^b(\M)$. Then we have
$R\Gamma(s)=S_1\oplus S_2$, where $S_i:DL(M_i\galochka)\to D_h(M_i)$ for
$i=1,2$, and $S_1$ is an isomorphism. We have $R\Gamma(s)=S_1\oplus S_2$
for some $S_1$, $S_2$ by Claim \ref{claim}. Now recall
that for $I\in \I_G$ the morphism $s:\cDL(I)\to j_*^G(I)$ includes in a
functorial morphism of exact sequences 
$$\begin{CD}
0 @>>>    \cDL  @>>>  \underline{Sm(\Gamma(I))} @>>> L(Sm(\Gamma(I)))
\longrightarrow 0 \\
@.     @VV{s}V           @VVV                         @VVV    \\
0 @>>>   j_!(I) @>>>     j_*^G(I)         @>>>
i^*j_*^G(I)\overset{\ref{cormain}}{=} \L(Sm\circ
\Gamma(I)) \to            0 
\end{CD}
$$
Thus we have a functorial morphism of distinguished triangles in $D^b(\tM)$

$$\begin{CD}
DL(Sm \circ R\Gamma(\F))  @>>> Sm\circ R\Gamma(\F) @>>> R\Gamma(L(Sm\circ
R\Gamma(\F)))  \\
@VV{R\Gamma(s)}V    @VVV         @VVV  \\
R\Gamma( j_!(\F))=R\Gamma_c(\F)        @>>> R\Gamma( Rj_*^G(\F))    @>>>
R\Gamma (\L 
(Sm\circ R\Gamma(\F)) )  
\end{CD}
$$
The second vertical map is an isomorphism by Lemma \ref{acycl}.

Also it is clear that decomposition $R\Gamma_c(\F)\cong M_1\oplus M_2$
splits the whole diagram into the direct sum of two. The condition $M_1\in
D^b(\R)$ implies that the right vertical arrow in the first  summand is an
isomorphism by Theorem \ref{cormain} ii), iii).  
Hence $S_1$ is also an isomorphism. By Lemma \ref{subcat} the
quasiisomorphism \eqref{d} in the category $D^b(\tM)$ yields the
corresponding quasiisomorphism in $D^b(\R)$.
 \epf

\begin{Rem} We could have used the Schneider-Stuhler ``localization'' Theorem
\cite{Schn1} instead of our Proposition \ref{loc} in order to deduce the
isomorphism \eqref{d} from Theorem \ref{cormain}.
\end{Rem}

\subsubsection{} We finish the chapter with a  generalization of the last
Theorem which connects  homological and Deligne-Lusztig dualities on $D^b(\M)$
(not only on complexes with admissible cohomology). 

Recall that the Bernstein center $\Ze$ is by the definition the algebra of
 endomorphisms of identity functor on $\M$; so it is a commutative algebra 
acting on every module $M\in \M$ commuting with the $G$-action. It is known
\cite{BDKV} that $\Ze$ is an infinite product of finitely generated commutative
 algebras; in fact each factor in this decomposition is isomorphic
to algebra of functions on an algebraic torus of dimension $\leq$ rank$(G)$
invariant under the action of a finite group.

 We also recall that the strong 
admissibility Theorem \cite{BDKV} asserts that for any $M\in \M$ and any open
compact subgroup $K\subset G$ the space of $K$-invariants 
$M^K$ is a finite $\Ze$-module  supported on a finite number of components of
 $Spec(\Ze)$.
If $M$ is a $G$-module, and $V$ is a $\Ze$-module
then the  space $Hom_{\Ze}(M,V)$ is naturally a $G$-module; moreover, from
the strong admissiblity Theorem it follows that $Sm(Hom_{\Ze}(M,V))\in \M$ if
$M\in \M$ and $V$ is locally finitely generated.
Since $\M$ has enough projectives 
 we can  derive the bifunctor $(M,V)\mapsto Sm(Hom_{\Ze}(M,V))$
 in both arguments to get a
 bifunctor $D^-(\M)^{opp} \times D^+(\Ze-mod) \to D^+(\M) $, where
$\Ze-mod$ is the category of locally finitely generated $\Ze$-modules. Let
 $D\in D^b(\Ze-mod)$  be the Grothendieck dualizing complex. 

We denote by $D_{Gr}$ the Grothendieck-Serre duality
 functor $RHom(\underline{\ \ }, D)$; we have $D_{Gr}:D^b(\M)\to D^+(\M)$

 Note that for a finite dimensional 
module $V$ over a commutative algebra we have $D_{Gr}(V)\cong V^*$
canonically. 
It follows that $D_{Gr}|_{D^b(\R)}\cong \galochka$. 

\begin{Rem} The 
quotient of a smooth variety by an action of a finite group is known to be
 Cohen-Macaulay, hence the dualizing complex $D$ is concentrated in one
homological 
 dimension. Also, for such a variety the Grothendieck-Serre duality
 functor sends the bounded derived category
$D^b(Coh)$ into itself, where $Coh$ stands for coherent
 sheaves. It follows that $D_{Gr}$ sends $D^b(\M)$ to itself.
 We will not use these facts below; the latter one follows also from the
next Theorem. 
\end{Rem}

\begin{Thm}\label{main_dual}
 We have a canonical isomorphism of functors 
on $D^b(\M)$
\begin{equation}\label{d1}
D_h\cong DL\circ D_{Gr}
\end{equation}
\end{Thm}

In order to prove the Theorem
we need to extend the statement \ref{dual} to the following set-up.
Let $B$ be a commutative unital $k$-algebra of finite type.
 Let  $\tM_B$  be the category of smooth $G\times B$ modules,
i.e. $B$-modules equipped with a smooth $G$-action;
in other words $\tM_B$ is a category of modules
over the algebra $\Heck_B:= \Heck \otimes B$ . The category
 $\M_B \subset \tM_B$ consists of all finitely generated modules.
We also define $\tdM_B$ to be  the  category of complete topological
$G\times B$ modules  having a basis of neighbourhoods of 0
consisting of $B$-submodules such that the  quotient is a finite $B$-module;
the subcategory
 $\dM_B \subset \tdM_B$ consists of modules which are topologically
 finitely generated. 

We call a module $M\in \M_B $ or $M\in \dM_B$ admissible if
the space of $K$-invariants in $M$ is a finite $B$-module for any open compact
subgroup $K\subset G$. We denote by
 $\R_B \subset \M_B$ and $\dR _B\subset \dM_B$ the (full)
subcategories of admissible modules.

\begin{Lem}  The category $D^+(\R_B)$ (respectively  $D^b(\R_B)$) is
equivalent to the full subcategory of   $D^+(\tM_B)$ (respectively
$D^b(\tM_B)$) 
consisting of complexes with admissible cohomology.
\end{Lem}
\prf is parallel to that of \ref{subcat}. \epf

For $V\in B-mod$ we have a contravariant functor $Hom_B(\underline{\  \ }, V):
\M_B \to \dM_B$. We can derive it in both arguments to get a bifunctor
$D^-(\M_B)\times D^+(B-mod)\to D^+(\tdM_B)$.

We have Deligne-Lusztig functors $DL:D^{b,\pm}(\M_B)\to D^{b,\pm}(\M_B)$;
 the definition remains the same. 

\begin{Prop}\label{admze}
 For any  $V\in D^b(B-mod)$ there exists a canonical isomorphism
of functors $D^b(\R_B)\to D^+(\R_B)$
$$RHom_{G\times B} (\underline{\ \ }, \Heck \otimes V)|_{\R_B} \cong DL\circ
 RHom_B(\underline{\ \ }, V)$$
\end{Prop}

\prf is parallel to the proof of \ref{dual}. The only difference
is that one should consider $G$-equivariant sheaves of $B$-modules
instead of  $G$-equivariant sheaves of $k$-vector spaces, and replace
the Verdier dualizing sheaf by the tensor product of Verdier dualizing sheaf
(with coefficients in $k$) and $V\in D^+(B-mod)$.\epf

\subsubsection{}
By the definition we have imbedding of a full subcategory $I:\M \hookrightarrow
  \M_{\Ze}$.
By the strong admissibility Theorem its image actually lies in
$\R_{\Ze}$. So to finish the proof of the Theorem it is enough to prove
the following standard fact.

\begin{Prop}\label{ipr}
  We have a functorial isomorphism 
 $$RHom_{G\times \Ze} (I\otimes Id(M), N \otimes D)\cong I(RHom_G (M , N))$$
where $M\in D^-(\M)$; $N$ is a bounded below complex of finitely
generated smooth $G$-bimodules;
$D\in D^b(\Ze-mod)$ is the dualizing complex. Here both sides lie
in $D^+(\tM_\Ze)$.
\end{Prop}

\prf We have an isomorphism
$$Hom_{G\times \Ze} (I(M),?)=Hom_{G}(M,Hom_{G\times \Ze}(I(\Heck), ?)),$$ 
where
the action of $G$ on $Hom_{G\times \Ze}(I(\Heck), ?)$ comes from 
the right $G$-action on $\Heck$; and also:
$$Hom_{G\times \Ze}(I(\Heck), ?)=Hom_{\Ze\otimes \Ze}(\Ze, ?)$$
where the first of the two $\Ze$-actions on $?$ comes from the 
$G$-action on $?$, and $G$ acts on $Hom_{\Ze\otimes \Ze}(\Ze, ?)$ via its 
action on $?$. Since $Hom_{\Ze\otimes \Ze}(\Ze, ?)$
is injective provided $?$ is injective,
we get the corresponding isomorphisms of derived functors.

We have the standard quasiisomorphism  in the derived category of $\Ze$-modules
$$ RHom_{\Ze\otimes \Ze} (\Ze, N\otimes D)=i^!\circ pr_1^! (N)\cong
 ( pr_1\circ i)^! N \cong N$$
where  $i$ is the diagonal embedding $i:Spec(\Ze)\to
  Spec(\Ze\otimes \Ze)$, $pr_1: Spec(\Ze\otimes \Ze)\to Spec(\Ze)$ is the
 projection, and we do not distinguish between a module over a commutative
algebra and the corresponding quasicoherent sheaf over its spectrum.

It is easy to lift this quasiisomorphism to a quasiisomorphism in the category
of $G$-bimodules. More precisely, if $N$ lies in the heart of 
the derived category (i.e. $H^i(N)=0$ for $i\not = 0$;  the only case
we need $N=\Heck$ satisfies this assumption) then 
$ RHom_{\Ze\otimes \Ze} (\Ze, N\otimes D)$ has only non-trivial cohomology
 in degree 0  isomorphic to $N$, hence in this case
we have a canonical quasiisomorphism
$N\leftarrow \tau_{\leq 0}(  RHom_{\Ze\otimes \Ze} (\Ze, N\otimes D))\to
  RHom_{\Ze\otimes \Ze} (\Ze, N\otimes D)$. If $N^\cdot$ is any complex 
(of injective modules) and $D^\cdot$ is any complex of injective $\Ze$-modules
quasiisomorphic to the dualizing sheaf, then $RHom_{\Ze\otimes \Ze}
 (\Ze, N\otimes D)$ is represented by the bicomplex $Hom _{\Ze\otimes \Ze} (\Ze, 
N^\cdot \otimes D^\cdot) $; we can apply the above
canonical quasiisomorphism in each column of this bicomplex to get the desired 
quasiisomorphism for any $N$. 
\epf

\subsubsection{} \prf of Theorem \ref{main_dual}. Since $I(M)\in \R_{\Ze}$
for $M\in \M$ we see  using \ref{ipr} that
$$D_h(M)=RHom_G(M,\Heck)\cong RHom_{G\times \Ze}(I(M), \Heck \otimes D)$$
(more precisely we apply the forgetful functor $G \times \Ze-mod\to
G-mod$ to the isomorphism \ref{ipr} to get the last isomorphism).
By \ref{admze} we have:
$$RHom_{G\times \Ze}(I(M), \Heck \otimes D)\cong DL(RHom_{\Ze}(M,D))$$  

The Theorem is proved. \epf

\end{section}

\newpage

\begin{section}{Elliptic pairing and filtrations on the Hecke
algebra}\label{2} 

\subsection*{Preliminaries and notations} In this chapter we prove

\begin{Thm}\label{Main} The equality \eqref{Ochi}:
$$O_{g^{-1}}(\<\rho\>)= \chi_\rho(g)$$
holds for any $\rho\in\R$, and any elliptic element $g\in G$.
\end{Thm}
and deduce from it

\begin{Thm}\label{Kazhdan} Assume $\char(F)=0$ and $k=\Ce$.
The equality \eqref{Eulchar}:
$$ \sum (-1)^i \dim
Ext^i(\rho_1,\rho_2)=
\int_{Ell} \chi_{\rho_1}(g^{-1})\chi_{\rho_2}(g)d\mu(g)
$$
is valid for all $\rho_1,\,\rho_2\in \R$.
\end{Thm}

It will be convenient to allow  
$G$ to be an arbitrary reductive group (not requiring that
$G$ has  compact center). Notice however that a group with non-compact
 center does not contain elliptic elements (in the terminology we adopt),
 so the  RHS of \eqref{Eulchar} vanishes; it is well-known
 (see  Claim \ref{noncomp}
below) that  so does the LHS. The statement of Theorem \ref{Main} 
 is vacuous for a group with noncompact center. 

\subsubsection{Hattori-Stallings trace for the Hecke algebra}\label{HaS}
 Since the Hecke algebra $\H$ is not unital, 
a comment on the definition of the Hattori-Stallings trace appearing in
the formulation of Theorem \ref{Main} is required.

Recall from \cite{BDKV}
that for a {\it nice} open compact subgroup $K\subset G$
(see \cite{BDKV} 2.1b) for the definition) the full subcategory $\M ^K\subset 
\M$ consisting of modules generated by $K$-fixed vectors is a direct summand in
$\M$. Moreover, the functors $\Mmod \mapsto \Mmod^K$ and $M\to M\otimes 
_{\H(G,K)}\H$ are mutually inverse equivalences between $\M^K$ and the category
of $\H(G,K)$-modules (here $\H(G,K)$ is the
 subalgebra of $K$-biinvariant measures in $\H$).

Hence $\H(G,K)$ is a Noetherian unital algebra of finite homological dimension,
 so the usual construction of the Hattori-Stallings trace applies to it.
Now for any $\Mmod \in \M^K$, $E\in End_G(\Mmod)$ 
 and  $K'\subset K$ we have   
$\Mmod^{K'}=\Mmod^K\otimes _{\H(G,K)}\H(G,K')$, hence
%%%%by ??
 the natural map $\H(G,K)/[,] \to \H(G,K')/[,]$ induced by the imbedding
   $\H(G,K) \imbed \H(G,K')$ sends $
Tr_{H-St}(\Mmod^K, E|_{\Mmod^K})$ to  $
Tr_{H-St}(\Mmod^{K'}, E|_{\Mmod^{K'}})$. Thus for small enough $K'$
the image of
 $Tr_{H-St}(\Mmod^{K'}, E|_{\Mmod^{K'}})$ in $\H/[\H,\H]$ does not depend 
on $K'$; this is by the definition the Hattori-Stallings trace of $(\Mmod ,E)$.

\subsubsection{}\label{Kellipt}
 A regular elliptic element $g\in G$ is fixed throughout the argument.

 Since $g$ is elliptic it is also compact; we fix a maximal  compact subgroup
$K_0$ containing  $g$.

Let $K\subset K_0$ be a normal subgroup; since $g$ normalizes $K$ it acts
on the set of double cosets $K\backslash G/H$ for any subgroup $H\subset G$.
We say that $g$ is {\it $K$-elliptic} if it acts  on $K\backslash G/P$
without fixed points  for any proper parabolic $P\subset G$. 
Obviously $g$ is $K$-elliptic for small enough normal open $K\subset K_0$.

We say that $K$ is {\it nice} if it is nice with respect to some minimal Levi
subgroup $L_0$, which is in good position with $K_0$ (see \ref {nota}, 
\ref{2.1} below). 

 We  fix a  normal open subgroup
$\K\subset K_0$, such that
 $g$ is $K$-elliptic.

{\it From now on we will denote by $\Heck$  the subalgebra
$\Heck(G,\K)\subset \Heck(G)$ of $\K$-biinvariant distributions.}

\medskip

We will in fact work with the following form of Theorem \ref{Main}

\begin{Thm}\label{iii}
 One can define for any $g\in Ell$, and any finitely generated
$\H$-module $M$ a number
$Tr(g,M)$ $\in k$ so that

\ i) $Tr(g,M)$ is additive on short exact sequences in $M$.

Fix $\Mmod \in\M$. 

\ ii) If $\Mmod$ is admissible, and its character $\chi_\Mmod$ 
is constant on the coset $K\cdot g$,
 then $Tr(g,\Mmod^K)=\chi_\Mmod(g)$
 is the character value.

\ 
iii) Let $M$ be a projective $\H$-module.
 Then there exists an idempotent $E\in
Mat_m(\H)$ such that $M$ is isomorphic to the image of $E$ acting on
$\H^{\oplus m}$ on the right. 

Suppose that $K$ is nice, and the
 locally constant 
invariant function $x\to O_{x^{-1}}\left(\< \Mmod\> \right)$ 
defined on the regular
elliptic set is (defined and) constant on the coset $\K\cdot g$. 
Then we have
$$
Tr(g,M)=O_{g^{-1}}(\sum E_{ii}).
$$
\end{Thm}

It is immediate to deduce
\ref{Main} from \ref{iii} (see \ref{proofMain}).

The next three sections are devoted to the  construction of
$Tr(g,M)$ and verification of the properties i)--iii).

\subsubsection{Notations}\label{nota}
 Let us choose a minimal parabolic with Levi decomposition
$P_0=L_0\cdot U_0\subset G$ so that $L_0$ is in good relative position
with $K_0$, i.e. the point $p$ of the Bruhat-Tits building fixed by $K_0$
lies in the apartment normalized by $L_0$.

 A standard parabolic is the one containing $P_0$;
a standard Levi is the Levi subgroup of a standard parabolic containing $L_0$.

The letters $P$, $L$, $U$ with indices will be reserved 
for a standard parabolic, its standard Levi, and its unipotent radical
respectively, unless  stated otherwise.

 Let $W=Nm(L_0)/L_0$, $W_L=(Nm(L_0)\cap L)/L_0$
  and
$W_{aff}=Nm(L_0)/L_0^c$ be respectively the Weyl group, the Weyl group of $L$
and the affine Weyl group.

For an $L$-module $\rho$ we will sometimes write $i_L^G(\rho)$ instead of
$i_P^G(\rho)$  for the parabolically induced module (of course
there is no ambiguity here, because $P$ and $L$ determine each other
uniquely). 

Let $A_L$ be the center of $L$, $X   _L$ be the lattice of coweights of $A_L$.
Thus $X   :=X   _{L_0}$ is the lattice of abstract coweights of $G$.
Set $\a=X   \otimes \RE$; let $\a_+\subset \a$ (respectively $\a^+$)
 be the dominant Weyl chamber (respectively its closure).
For a subset $?\subset \a$ we denote the intersection $?\cap \a^+$
by $?^+$.
 
For a standard Levi $L$ let $\a_L\subset  \a$ denote $X_L\otimes \RE$.

Fixing a uniformizer ${\goth p}\in F^\times$
we obtain a canonical 
imbedding $\iota_L:X   _L\imbed A_L$, $\chi
\mapsto \chi({\goth p})$, inducing an imbedding of finite index
$X   _L\imbed L/L^c$.

For a topological group $H$ we denote by $H^c$ the subgroup generated
by all compact subgroups.

 Consider the composition
$L/L^c\imbed (L/L^c)
\otimes \RE=\a_L\imbed \a$. We identify $L/L^c$ with its 
image under this map.

\medskip 

By $\P_L$ we denote the orthogonal projection $\a \to \a_L$.
 
By a root we will mean a nonmultipliable root, and by a coroot the
corresponding coweight of $A_{L_0}$. Thus coroots (respectively, roots)
 form a reduced root
system in $\a$ (respectively $\a^*$) (see \cite{BoT}, Theorem 7.2).
The choice of a minimal parabolic $P_0$ determines the basis of simple
(co)roots, and the dual basis of fundamental (co)weights.

 Let $\Sigma^+$ be the set of positive coroots; for a Levi
$L$ let  $\Sigma^+_L\subset \Sigma^+$ be the set of positive coroots of 
the derived group $L^{[,]}$.

If $\alpha$ is a simple coroot, $\omega_\alpha$ is the dual fundamental
weight,
 and $\lambda_1,\,\lambda_2 \in X   $ then we
 will write $\lambda_1 \la \lambda_2$  meaning $(\lambda_1,\omega_\alpha )
\leq  (\lambda_2,\omega_\alpha )$. We also write $\lambda \la n$ meaning
$ (\lambda,\omega_\alpha) \leq n$.

 The standard partial order 
 $\preceq$ on $\a$ is defined by $\lambda_1 \preceq \lambda_2$
if $\lambda_1 \la \lambda_2$ for all $\alpha$; as usual 
we write $\lambda_1 \prec
\lambda_2$ instead of $\lambda_1 \preceq \lambda_2\ \& \ \lambda_1 \not =
 \lambda_2$ etc.

Recall the notation  $k[S]:=\bigoplus\limits_{s\in S}k$.
When $S$ is a (semi)group this is the (semi)group algebra;
for $s\in S$ let  $[s]\in k[S]$ be the delta-function of $s$.

If $(\Lambda, \treugol)$ is a (partially) ordered semigroup,
 we use the notation  
 $\Lambda _{\treugol \lambda}:= \{\mu \in X    \mid \mu\treugol
\lambda \}$ for $\lambda\in \Lambda$.

We have a $(\Lambda, \treugol)$-filtration on the semigroup algebra
$k[\Lambda]$, defined by $ k[\Lambda]_{\treugol \lambda}:=
k[\Lambda_{\treugol \lambda}]$.

Analogous notations apply when $\treugol$ is a filtration on $\Lambda$
indexed by another (partially) ordered semigroup. 

 We use the same symbol for a partial order on $\a$
 and its restriction to  a sublattice.

We will say that $\lambda \in X    _L$ is dominant (or large) enough
if $( \lambda, r )\gg 0$ when  $r$ is a root of $U$.
(The latter condition will be abbreviated as $\< \lambda,
 \Sigma^+-\Sigma^+_L\> \gg 0$.

\subsection{``Spectral'' filtration: definition of the trace functional}
\label{1spar}
\subsubsection{}\label{ograni1} 
Let 
 $\rho$ be an admissible representation of $L$. We can form a
 (non-admissible) $G$-module  $\Pi_\rho:=i_L^G(\rho\otimes
 k[L/L^c])=ind_{L^c\cdot U}^G(\rho|_{L^c})$.
Then $\Pi_{\rho}$ has also a $
 k[L/L^c]$-action commuting with the $G$-action, which comes  from
the  action of $k[L/L^c]$ on the second multiple of $\rho\otimes
 k[L/L^c]$:
\begin{equation}\label{canact}
l(f)(g)=(Id\otimes l)(f(g))
\end{equation}
for $l\in L/L^c$, $f:G\to \rho\otimes
 k[L/L^c]$, $f\in \Pi_\rho$.
For $\lambda\in L/L^c$ let $[\lambda]_\rho\in \End(\Pi_\rho)$ denote 
action \eqref{canact}.

For an unramified character $\psi:L/L^c\to k^\times$ we have a canonical 
isomorphism of $G$-modules  (but not of $L/L^c$-modules) $I_\psi:
\Pi_{\rho}\iso \Pi_{\rho\otimes \psi}$, satisfying 
\begin{equation}\label{commrel}
I_\psi(l(m))=\psi(l)
l(I_\psi(m))\end{equation}
 for $l\in L/L^c$, $m\in \Pi_\rho$.

From \eqref{commrel} it follows that $E\in \End(\Pi_\rho)$ commutes
with $\lambda_\rho$ for all $\lambda\in L/L^c$ iff $I_\psi\circ E\circ
I_\psi^{-1} $ commutes with all $[\lambda]_{\rho\otimes \psi}$, thus we
have a canonical isomorphism
\begin{equation}\label{isoend}
\End_{L/L^c}(\Pi_\rho)\cong \End_{L/L^c}(\Pi_{\rho\otimes\psi}).
\end{equation} 

For any standard Levi $L$ let us choose a set $Cusp_L$ of cuspidal irreducible
representations of $L$, in such a way that any  cuspidal irreducible
representation of $L$ is isomorphic to a unique representation $\rho
\otimes \psi$, where $\rho\in Cusp_L$, and $\psi:L/L^c\to k$
 is an unramified character of $L$.
We can (and will) make this choice in such a way that $\rho|_{\iota_L(X_L)}$ is
trivial for $\rho\in Cusp_L$.

 Let  $Cusp^K_L\subset Cusp_L$ be the set of
such $\rho\in Cusp_L$ that $i_L^G(\rho)^K\ne 0$;
set $Cusp:=\cupl_{L} Cusp_L$ and
 $Cusp^K:=\cupl_{L} Cusp_L^K$.

\medskip

Set $\Htil:= \bigoplus\limits_{\rho\in Cusp^\K}
\End_{k[L/L^c]}(\Pi_\rho^{\K})$.

We have $[\lambda]_\rho\in \Htil$ for $\rho\in Cusp_L^K$.
For $\lambda\in X$ define $[\lambda]\in\Htil$ by 
\begin{equation}\label{[lam]}
\lambda =\suml_{L}\suml_{\rho\in Cusp^K_L}
[\lambda]_\rho   ,
\end{equation}
 where $L$ runs over the set of standard Levi subgroups, such that $\lambda
\in X_L$.

We have an imbedding $I=\suml_{\rho\in Cusp^K}I_\rho :\H \imbed \Htil,$
 where $I_\rho:\H \to \End(\Pi_\rho^K)$.

\subsubsection{}
We are going to define an $\a$-multifiltration on $\Htil$.
 To clarify the idea we first carry out the construction under an 
 additional simplifying assumption
 that the maximal compact subgroup $K_0$
containing $g$ is good, i.e. $G=K_0\cdot P_0$
(the general case will be treated in the next subsection).

If $K_0$ is good, then the $k[L/L^c]$-module $\Pi_\rho$ can be canonically
trivialized in the following fashion.

By the definition the space of the induced representation $\Pi_\rho$
 is the space of functions on $G$ with values in 
 $\rho\otimes k[L/L^c]$ which transform accordingly under the right
 action of $P$. Let us say that $f\in \Pi_\rho$ is constant if 
$f|_{K_0}$ takes values in $\rho \otimes 1\subset
 \rho\otimes k[L/L^c]$. Then the space  of constant
elements is a $K_0$-submodule, and  generates $\Pi_\rho$ freely
as a $k[L/L^c]$-module.
We denote the space of constant
 elements in $\Pi_\rho$ by $(\Pi_{\rho})_0$.
For an unramified character $\psi$ of $L$ the isomorphism $I_\psi$ sends
 $(\Pi_{\rho})_0$ to $(\Pi_{\rho\otimes \psi})_0$,
 and  we have the canonical 
isomorphism $(\Pi_{P,\rho})_0\iso i_P^G(\psi \otimes \rho)
=\Pi_\rho \otimes_{k[L/L^c]}\psi$.

We abbreviate  $(\Pi_{\rho})_0^{\K}=i_0 (\rho)$.
So we have  $\Htil = \oplusl_{Cusp}  \End_{k[L/L^c]}(\Pi_\rho^{\K})=
 \oplusl_{Cusp^K}\End(i_0(\rho))\otimes k[L/L^c]$.

\smallskip

For a (partial) order $\treugol$ on $\a$ we have an $(\a,
\treugol)$-filtration on $\Htil$ given by $\Htil_{\treugol \lambda}: =
 \bigoplus\limits_{\rho\in
Cusp^\K}\End(i_0(\rho))\otimes k[L/L^c]_{\treugol \lambda}$.

\subsubsection{} Consider now the case when $K_0$ is an arbitrary maximal 
compact subgroup. 

Recall that $p$ is the point of the Bruhat-Tits building $X$ fixed by 
$K_0$; let $\Ag\owns p$ be the  apartment corresponding to  $L_0$. 
Let  $W_p=K_0\cap Nm(L_0)/L_0^c\subset
W_{aff}$ be the stabilizer of $p$, and $W_{aff}^L\subset W_{aff}$ be the 
subgroup $(L^c\cap Nm(L_0))/L_0^c$. The groups $Nm(L_0)$, $L^c\cap
Nm(L_0) $ act on $\Ag$
respectively through $W_{aff}$, $W_{aff}^L$.

Notice that  $P^c=L^c\cdot U$.

\begin{Lem}\label{bijection}
The map $w\mapsto K_0 w P^c$ provides a bijection
\begin{equation}\label{KP^c}
 W_p\backslash W_{aff}/W_{aff}^L
\iso K_0\backslash G/P^c\ .
\end{equation}
\end{Lem}
\prf First notice that since $K_0$, $L_0$ are in good relative position we have
$K_0\supset L_0^c$, and obviously $L_0^c \subset  P^c$, so the
expression $ K_0 w P^c$ is meaningful. 

Let $W_{aff}(L)$ denote $Nm_{L}(L_0)/L_0^c$.
On the LHS of  \eqref{KP^c} the group
$W_{aff}(L)/W_{aff}^L$ acts from the right. From the Bruhat
decomposition $L=(P_0\cap L)\cdot Nm_L(L_0) \cdot (P_0\cap L)$ we see,
using the inclusion $P_0\cap L=  L_0\cdot (U\cap L)\subset L_0\cdot L^c$,
  that the map 
$W_{aff}(L)/W_{aff}^L=Nm_{L}(L_0)/Nm_{L^c}(L_0)\to L/L^c$
 is an isomorphism. 

On the RHS of \eqref{KP^c} the group $L/L^c=P/P^c$ acts from the right. The
map \eqref{KP^c} obviously agrees with the action of $L/L^c$. Moreover, we
claim that the action on both sides is free. Indeed, the right action of 
$L/L^c=P/P^c$ on $G/P^c$ is obviously free; since the map 
$G/P^c\to  K_0\backslash G/P^c$ is proper, the stabilizer of any point in $
K_0\backslash G/P^c $ is compact. However $L/L^c$ is a (discrete) free
abelian group, so it does not contain compact subgroups. The same argument
applies to the LHS (with ``finite-to-one'' replacing ``proper'').

Thus it is enough to prove that the map $ W_p\backslash W_{aff}/W_{aff}(L)
\to K_0\backslash G/P$ is an isomorphism. 
This can be reformulated via the Bruhat-Tits building $X$ as follows.

Any points $x\in X$, $z\in \Xbar-X$ lie in the closure of some apartment
$\Ag$. Moreover, 
if $x'\in \Ag$, $z'\in \Agbar-\Ag$ is another pair of points, such that
$x'=g(x)$, $z'=g(z)$, then also $x'=w(x)$, $z'=w(z)$ for some element $w\in
W_\Ag$, where $W_\Ag=Nm(\Ag)/Stab(\Ag)\cong W_{aff}$. 

The fact that $x$ and $z$ lie in one apartment follows from
\cite{BT1} 7.4.18(ii). To check the second statement
recall that by Lemma \ref{c_to_z} of part \ref{1} there exists a unique
geodesic ray connecting a point in $X$ and a point in $\Xbar-X$. In
particular  the rays $[x,z)$, $g\([x,z)\)=[x',z')$ lie in $\Ag$ since there
endpoints lie in $\Agbar$. Hence by \cite{BT1} 7.4.8 there exists $w\in
W_\Ag$ such that $g|_{[x,z)}=w|_{[x,z)}$. By continuity of the action of
$g$, $ w$ we get $g(z)=w(z)$. \epf

\subsubsection{} For $x\in W_p\backslash W_{aff}/W_{aff}^L$ let 
$(G/U)_x$ be the corresponding orbit of $K_0\times P^c$.

Let the subgroup $W_f^p\subset \Aut(\Ag) $ be generated by
 reflections in all hyperplanes parallel to root 
hyperplanes and passing through $p$; let
 $W_{aff}^p$ be  generated by 
$W_{aff}$ and $W_f^p$. Let $\Lambda$ be the 
intersection of $W_{aff}^{p}$ with the group of translations; it is 
obvious that $\Lambda$ is a lattice in $\a=X   \otimes \RE$ containing 
$X   $, and $W_{aff}^p=W_f^p\semidir \Lambda$.

Consider the projection $\pi_L:W_p\backslash W_{aff}/W_{aff}^L
\to W_f^p \backslash W_{aff}^p/W_{aff}^L = \Lambda/W_{aff}^L
\overset{\P_L}{\longrightarrow}
 \a_L\subset \a$; denote its image by
$\Lambda_L$.
For $\lambda \in \Lambda_L$ set
$(G/U)_\lambda=\cupl_{\pi_L(x)=\lambda} (G/U)_x$.
 
Now we get a $\Lambda_L$-grading on $\Pi_{\rho}$ for  $\rho\in Cusp_L$
 as follows
\begin{equation}\label{grad}
\Pi_{\rho}=C_c(G/U)^K\otimes_{\H(L^c)} \rho =
\oplusl_{\lambda\in \Lambda_L}C_c(G/U)_\lambda^K \otimes_{\H(L^c)} \rho.
\end{equation}

For $\lambda\in \a$ we now set 
\begin{equation}\label{gradHtil}
\Htil_\lambda  = \{ h\in \Htil\,
|\, h(\Pi_{\rho}^K)_\nu\subset (\Pi_{\rho}^K)_{\nu+\lambda}\},\end{equation}
\begin{equation}
\Htil _{\treugol \lambda} = \oplusl_{\mu \treugol \lambda}\Htil_\mu.
\end{equation}
 Comparison of the definitions shows that the latter
filtration coincides with the one introduced in \ref{ograni1}
 when $K_0$ is good.

\subsubsection{}\label{ograni3}
 For a (partial) order $\treugol$ on $\a$ we have
defined an $(\a, \treugol )$-filtration on $\Htil$.
Let $\H_{\treugol \lambda}$ be the induced filtration on $\H\subset \Htil$.

\medskip

The next Proposition is the key technical result needed
for our definition of the trace functional.

For a standard Levi $L_1$, $\rho \in Cusp_{L_1}$
 the space of $\Pi_ {\rho}$ is  the space of
$\rho\otimes k[L_1/L_1^c]$-valued functions on $G$ transforming
 accordingly under the right
action of $P_1$;
 for $L\supset L_1$ and $x\in K\backslash G/ P$ let
$m_x^{\rho}\in
\End(\Pi_{\rho}^K)$ be  multiplication by the
$\delta$-function of the corresponding right $K\times P$ double coset.

Set also  $m_x=\suml_{L'\subset L} \suml_{\rho\in Cusp_{L'}^K} m_x^{\rho}$.

\begin{Prop}\label{image} Fix a standard Levi $L\supset L_1$,
 $\rho \in Cusp _{L_1}^K$,
  $x\in  K\backslash G/ P$ and $a\in \RE_{>0}$.
Then for dominant enough $\lambda \in X   _{L}^+$ there exists an element
$h_{x,\lambda}^\rho \in \H$
such that $I(h_{x,\lambda}^\rho)=\suml_{L_2=w(L_1)\subset L}
 [\lambda]_{w(\rho)} \cdot m_x^{w(\rho)} +h'$, 
where $h'\in  \suml_{r\in\Sigma^+-\Sigma^+_L}\Htil_{\preceq \lambda - a\cdot
r}$.
\end{Prop}

\begin{Rem} The construction of the trace functional will use only the 
fact  $I(h_{x,\lambda}^\rho)-\suml_{L_2=w(L_1)\subset L}
 [\lambda]_{w(\rho)} \cdot m_x^{w(\rho)} \in \Htil_{\prec \lambda}$.
The stronger condition stated in the Proposition will play a role only
in the proof that the constructed functional satisfies property iii) of
\ref{iii}.
\end{Rem}

 We will deduce the Proposition from the {\it matrix Payley-Wiener
Theorem} which describes the image of the imbedding $I$ (see \cite{Be}, \S 21).
Let us first recall its contents. 

If  $M_1,M_2$ are two modules over a free abelian group 
 ${\Bbb Z}^n$
then by a rational morphism from $M_1$ to $M_2$ we will mean an
element of $Hom_{Frac(k[{\Bbb Z}^n])}\( M_1\otimes_{k[{\Bbb Z}^n]} Frac(k[{\Bbb
 Z}^n]),
M_2\otimes_{k[{\Bbb Z}^n]} Frac(k[{\Bbb Z}^n])\)$ where $Frac(?)$ is the
 field of 
 fractions of a commutative ring $?$.

We recall the necessary part of the theory of intertwining operators. 
Suppose that standard Levi $L_1$, $L_2$ are conjugate, so we have
$L_2=w(L_1)$ for some $w\in W$.

Then there exists a canonical rational isomorphism of $G$-modules 
$I_w=I_w^\rho:\Pi_\rho\to \Pi_{w(\rho)}$, such that $I_w\circ [\lambda]_\rho
=[w(\lambda)]_{w(\rho)}\circ I_w$ for
$\lambda\in L/L^c$ 
 satisfying:

i) Transitivity: if $L_2=w_1(L_1)$, $L_3=w_2(L_2)$ are standard Levi
 then $I_{w_2w_1}(m)=I_{w_2}\circ I_{w_1}(f_{w_1,w_2}m)$ for some rational
function $f_{w_1,w_2}\in Frac(k[L_1/L_1^c])$.

ii) Compatibility with induction: if $L\supset L_1, L_2$, $L_2=w(L_1)$
 then for a representation $\rho\in Cusp_{L_1}$ we have
$I_{w}^\rho=i_L^G(I^\rho_{w,L})$, where $I^\rho_{w,L}:
i_{L_1}^L(\rho\otimes k[L_1/L_1^c])\to i_{L_2}^L(\rho\otimes k[L_2/L_2^c])$
is the intertwining operator between the $L$-modules.

\smallskip

We define the rational morphisms $R_w^\rho:\End_{L_1/L_1^c}
( \Pi_\rho)\to \End _{L_2/L_2^c}
(\Pi_{w(\rho)}) $ by 
$R_{w}(E)=I_{w}\circ E \circ I_{w}^{-1} $. Notice that the operators
$R_w$ satisfy the stronger transitivity condition:
 $$R_{w_1w_2}=R_{w_1}\circ R_{w_2}.$$

Let $L_1$, $L_2$ be standard Levi subgroups such that $L_2=w(L_1)$ and
$\rho_1\in Cusp_{L_1}$,  $\rho_2\in Cusp_{L_2}$ be such that
$\rho_2=w(\rho_1)\otimes \psi$ for an unramified character
$\psi$. Identify $\End_{L_2/L_2^c}(\Pi_{w(\rho_1)})$ with
$\End_{L_2/L_2^c}(\Pi_{\rho_2})$ by means of \eqref{isoend}; then
$R_w^\rho$ provides  a rational morphism
$\End_{L_1/L_1^c}(\Pi_{\rho_1})\to \End_{L_2/L_2^c}(\Pi_{\rho_2})$.
We denote it again by $R_w$ or $R_w^\rho$, this abuse of notation 
hopefully will not lead to a confusion.

For future reference we record the following

\begin{Lem}\label{actonweights} For standard Levi subgroups $L_2=w(L_1)$
and $\rho_1\in Cusp_{L_1}$,  $\rho_2\in Cusp_{L_2}$ such that
$\rho_2=w(\rho_1)\otimes \psi$ as above we have
   $R_w([\lambda]_{\rho_1})=
[w(\lambda)]_{\rho_2 }$ provided $\lambda$ lies in $X_{L_1}\subset L_1/L_1^c$.
\end{Lem}

\proof Since $R_w^{\rho_1}([\lambda]_{\rho_1})=[w(\lambda)]_{w(\rho_1)}$
we have only to check that
$[w(\lambda)]_{w(\rho_1)}=[w(\lambda)]_{\rho_2}$ (more precisely, that
the two sides agree under \eqref{isoend}). By \eqref{commrel} this reduces to 
$\psi(\mu)=1$ for $\mu \in X   _{L_2}$.
 But this is clear because $X   _{L_2}$
acts trivially both on $w(\rho_1)$ and on $\rho_2
\in Cusp_{L_2}$.  \epf

There exists also 
another type of operators acting on $\oplusl_{Cusp} \Pi_{\rho}$.
Namely, for a given irreducible representation $\rho$ of $L$ there might
exist a finite number of unramified characters $\psi$ such that $\psi \otimes
 \rho\cong \rho$. For such $\psi$ the canonical isomorphism
  $I_\psi:\Pi_\rho\iso \Pi_{\rho\otimes \psi}$ induces an automorphism of 
$\Pi_\rho$. This automorphism is defined uniquely up to a constant
(we have to choose an isomorphism $J: \rho\otimes \psi\iso \rho$, and then
get $i_L^G(J)\circ I_\psi: \Pi_\rho \iso \Pi_\rho$); thus we have a uniquely
 defined
automorphism $T_\psi$ of $\End_{L/L^c}(\Pi_\rho)$,
$T_\psi:  E\mapsto (i_L^G(J)\circ I_\psi) \circ E\circ (i_L^G(J)
\circ I_\psi)^{-1}$.

\medskip

The matrix Payley-Wiener Theorem (\cite{Be}, Theorem 33) asserts that 

{\it An element $E=\sum\limits_{Cusp^\K} E_{\rho}\in \Htil$
lies in the image of $I$ iff

{\rm i)} For any conjugate standard Levi subgroups $L_2=w(L_1)$ and 
$\rho_1\in Cusp_{L_1}$, $\rho_2\in Cusp_{L_2}$ such that  $\rho_2\cong
w(\rho_1\otimes \psi)$ for an unramified character $\psi$ we have
$R_w^{\rho_1}(E_{\rho_1})=E_{\rho_2}$.  

{\rm ii)} $T_\psi(E_{\rho})=E_{\rho}$ whenever $\psi\otimes \rho\cong \rho$.}

\subsubsection{} We return to the proof of \ref{image}.   

For standard Levi subgroups $L'=w(L)$ and $\rho\in Cusp_L^K$
we will not distinguish between $\Pi_w(\rho)$ 
and  $\Pi_{w(\rho)\otimes
\psi_w(\rho)}$ when this is not likely to lead to a confusion. 
We will also abuse notations writing $R_w$ instead of $T_{\psi_{w(\rho)}}
\circ R_w$.
Thus $R_w:\End(\Pi_\rho)\to \End(\Pi_{w(\rho)})$ is a rational
endomorphism of $\Htil$.

\medskip

The plan is to start
with $m_x^\rho\cdot [\lambda]_\rho$ and then average over the  set of
intertwining operators; we have to take care of the poles of intertwining
operators, and also to ensure that we really get an element with 
the desired highest term.
It will be  not hard to see that each summand in the resulting sum
is  $T_\psi$-invariant.

\begin{Lem}\label{2.5b}  If $ L_2=w(L_1)$, $L\supset L_1,\,L_2$ are standard
Levi,
 $x\in K\backslash G/P   $,  $\lambda \in X   _{L   }\subset X   _{L_1},
\, X   _{L_2}$ then
 we have
 $R_w(m_x^{\rho}\cdot [\lambda]_\rho)=
m_x^{\rho}\cdot [\lambda]_\rho$ for any $\rho$.

\end{Lem}

\proof By \ref{actonweights}
 the action by intertwining operators induces the standard
 geometric action on coweights, thus 
 $R_{w}([\lambda]_\rho)=[\lambda]_{w(\rho)  }$   % \otimes \psi_{w(\rho)}}$
 for $\lambda \in X   _{L   }$ and $w\in W_L$.
%%koriavo

If we think of $\Pi_\rho =i_{L_1}^G(\rho\otimes
 k[L_1/L_1^c])=    i_{L   }^G\left( i_{L_1}^{L   }(\rho\otimes
 k[L_1/L_1^c]) \right)$ as sections of the corresponding  sheaf
on $G/P   $ then the action of  $m_x^\rho$ 
comes from an action on the sheaf; moreover the corresponding endomorphism
 of a stalk is either identity or 0. From $R_w(Id)=Id$
 we conclude by compatibility that it sends $m_x^\rho$ to $m_x^{w(\rho)}$. \epf

 For $h\in \Htil$, $h=\suml_{L,
\lambda\in \Lambda_L} h_\lambda$, where $h_\lambda \in \Htil_\lambda$,
 let the support of $h$ be
\begin{equation}\label{supp^spec}
\supp(h):=\{\lambda \in \a
 \mvert 
h_\lambda \not = 0\} .
\end{equation}

\subsubsection{}
Consider the following situation. Let $L_1\ne L_2$ be two conjugate standard 
Levi lying in a standard Levi $L$ with $\rank(L)=\rank(L_1)+1=\rank(L_2)+1$.
Then necessary $L_2=w_0^L(L_1)$, where $w_0^L\in W_L$
 is the longest element.
Let $\alpha$ be a nonzero element of the group $X   _{L_1}\cap (\a_L)^\perp$.
 Notice
that $\alpha$ is proportional to $\P_{L_1}(r)$ where $r$ is the only  simple
coroot of $L$ which does not lie in $\Sigma_{L_1}$. 

\begin{Lem}\label{hyperplane} a)
 Components of the divisor of poles of
$R_{w_0^L}$ have  the form  $[\alpha]_\rho =\const$; in other words
for some $f\in k[t]$ for any $h\in \End_{L_2/L_2^c}\Pi_{w_0^L(\rho)}^K$
 we have $R_{w_0^L}(f([\alpha]_\rho)h)\in 
 \End_{L_2/L_2^c}\Pi_{w_0^L(\rho)}^K\subset 
 \End_{L_2/L_2^c}\Pi_{w_0^L(\rho)}^K\otimes Frac(k[L_2/L_2^c])$. 

b) There exists $n_0\in \Z$ such that the following holds. Suppose that
$h\in \End_{L_1/L_1^c}\Pi_{\rho}^K$ is such that $R_{w_0^L}(h)$ is
regular, i.e. belongs to
 $ \End_{L_2/L_2^c}\Pi_{w_0^L(\rho)}^K$ rather than to
 $ \End_{L_2/L_2^c}\Pi_{w_0^L(\rho)}^K
\otimes Frac(k[L_2/L_2^c])$. 
Then $\nu \in  \supp (R_{w_0^L}(h))\Rightarrow
 (w_0^L)^{-1}(\nu)+n\alpha \in \supp(h)$ for some $n$ with
 $|n|\leq n_0$.

\end{Lem}
\proof a) We can decompose  $L$ as an almost
direct product of its center $A_L $ and derived group $L^{[,]}$,
 i.e. we have a homomorphism $\pi:A_L
\times L^{[,]} \to L$ with finite kernel and cokernel. 
Obviously for $j=1,2$ the pull-back of the  representation
 $i_{L_j}^{L}(k[L_j/L_j^c]\otimes \rho)$ under $\pi$ decomposes as the tensor 
product of the free module $k[A_L/A_L ^c]$
with the parabolically induced module 
$i_{L_j\cap L^{[,]}}^{L^{[,]}} (\oplusl _{L_j/L_j^c\cdot A_L\cdot (L^{[,]}\cap 
L_j)}
 k[(L_j\cap L^{[,]})/(L_j\cap L^{[,]})^c]\otimes \rho|_{L_j\cap L^{[,]}})$. The
 intertwining
operator $I^\rho_{w_0^L,L}$ also decomposes as 
$Id \otimes I_{w_0^L, L^{[,]}}^{\rho|_{L_1\cap L^{[,]}}}$, hence so does
$R^\rho_{w_0^L,L}$. Thus the image
of the divisor of poles of  $R^\rho_{w_0^L,L}$
under the finite covering $Spec (k[L/L^c])\to
Spec(k[A  _{L}/A_L^c])\times Spec(k[(L^{[,]}\cap L_1)/
(L^{[,]}\cap L_1)^c]) $ is the product of  $Spec(k[A  _{L}/A_L^c])$
with a finite subscheme of $ Spec(k[(L^{[,]}\cap L_1)/
(L^{[,]}\cap L_1)^c])$  (notice that  $ (L^{[,]}\cap L_1)/
(L^{[,]}\cap L_1)^c$
is a free cyclic group, so   $ Spec(k[(L^{[,]}\cap L_1)/
(L^{[,]}\cap L_1)^c])$ is the punctured affine line). This finite 
subscheme is necessary killed
by the function $\prod\limits_{j=1}^n ([\alpha]_\rho-q_j)$
 for some finite set of elements
$q_1,...,q_n\in k^\times$, hence the same holds for the divisor of poles of
 $R^\rho_{w_0^L,L}$. a) follows by compatibility of
 the intertwining operators with induction.

Let us prove b). For $j=1,2$ consider the projection $pr_j:
K\backslash G/P_j^c \to K\backslash G/P^c   $.
For $y\in K\backslash G/P^c   $
denote $(\Pi^K _\rho)_y:=\oplusl_{\pi_1(x)=y}(\Pi^K _\rho)_x$,
 $(\Pi^K _{w_0^L(\rho)})_y:=\oplusl_{\pi_2(x)=y}(\Pi^K _{w_0^L(\rho)})_x$.
Notice that $(\Pi^K _\rho)_y$ is a $k[\Z \alpha]$-submodule for any
$y\in K\backslash G/P^c  $. 

By a) the intertwining operator $I_{w_0^L}$ is a well-defined morphism
$\Pi^K _\rho \otimes _{k[\Z \alpha]} k([\alpha]_\rho)\to 
\Pi^K _{w_0^L(\rho)} \otimes _{k[\Z w_0^L(\alpha)]} k([w_0^L(\alpha)]_{w_0^L(\rho)
})$.
From compatibility with induction it follows
that $I_{w_0^L}(\Pi^K _\rho)_y \otimes k([\alpha]_\rho)\subset (\Pi^K _{w_0^L(\rho)})
_y\otimes 
%?k([\alpha]_\rho)$. 
 k([w_0^L(\alpha)]_{w_0^L(\rho)})$.

Now from the definition it is easy to see that for $j=1,2$
 the following diagram is commutative:
$$
\begin{CD}
K\backslash G/P^c_j @>{pr_j}>> K\backslash G/P^c  \\
@V{\pi_{L_j}}VV                           @V{\pi_{L   }}VV\\
\Lambda_{L_j}                @>\P_{L_j}>>  \Lambda_{L   }
\end{CD}
$$

Hence for $\mu \in \Lambda_{L   }$ we have $k[\Z\alpha]$-submodule
$(\Pi^K _\rho)_\mu:=\oplusl_{\lambda \in\P_{L_1}^{-1}(\mu)} (\Pi^K _\rho)_\lambda$,
 and
a $k[\Z\cdot w_0^L(\alpha)]$-submodule 
$(\Pi^K _{w_0^L(\rho)})_\mu:=\oplusl_{\lambda \in\P_{L_2}^{-1}(\mu)}
 (\Pi^K _{w_0^L(\rho)}
)_\lambda$, so that $I_{w_0^L}((\Pi^K _\rho)_\mu\otimes k([\alpha]_\rho))
\subset (\Pi^K _{w_0^L(\rho)})_\mu \otimes k([w_0^L(\alpha)]_{w_0^L(\rho)})$. 

It follows that $R_{w_0^L}( (\End_{L_1/L_1^c}\Pi^K _\rho)_\mu) \subset
(\End_{L_2/L_2^c} \Pi^K _{w_0^L(\rho)})_\mu \otimes k([w_0^L(\alpha)]_
{w_0^L(\rho)})$
where the $[\alpha]_\rho$ 
(respectively $[w_0^L(\alpha)]_{w_0^L(\rho)}$)-invariant
$\Lambda_{L   }$-grading on $\End_{L_1/L_1^c}\Pi^K _\rho$ 
(respectively on $\End_{L_2/L_2^c} \Pi^K _{w_0^L(\rho)}$)
 is induced by the grading on $\Pi^K _\rho$, $\Pi^K _{w_0^L(\rho)}$ as in 
\eqref{gradHtil}.

Recall  the following elementary fact. 
For a Laurent polynomial $P=\sum a_i t^i\in k[t,t^{-1}]$ 
set $\supp(P)=\{n\in \Zet\,|\,a_n\ne 0\}$. Let $f\in k(t)$ be a rational
function, and set $n_f=\max(\operatorname{ord}_0(f),
- \operatorname{ord}_\infty(f))$. Then for $P\in k[t,t^{-1}]$ we have
$fP\in  k[t,t^{-1}] \Rightarrow \supp(fP)\subset[\inf(\supp(P))-n_f,\sup(\supp
(P))+n_f]$.

Now notice that for $\mu\in \Lambda_{L   }$:
$\P_{L_2}^{-1}(\mu)=w_0^L(\P_{L_1}^{-1}(\mu))$, and recall that by Lemma
\ref{actonweights} $R_{w_0^L}
([\alpha]_\rho h)
=[w_0^L(\alpha)]_{w_0^L(\rho)}R_{w_0}^L(h)$.

It is clear that $(\End_{L_1/L_1^c}\Pi^K _\rho)_\mu$ is a finite graded
 $k[\Z\alpha]$-module for $\mu\in \Lambda_{L   }$. Hence the above
elementary fact implies  that at least for fixed $\mu$ there exists such
$n_0$ that the condition b) of the Lemma holds for $h\in (\End_{L_1/L_1^c}
\Pi^K _\rho)_\mu$. But if it holds for some $h$ it also holds 
(with the same $n_0$) for
$h'=l(h)$ where $l$ is any element of $ A   _{L   }/A_L^c$. 
Since a finite number of components
of the $\Lambda_{L   }$-grading generate 
$\End_{L_1/L_1^c}\Pi^K _\rho$ as an $A   _{L   }/A_L^c$-module, b) is proved.
 \epf

\subsubsection{} The next statement implies Proposition \ref{image}.

Let us choose such $m\in \Z^{>0}$ that $\P_L(mr)\in X   _L$
for any $L$ and any coroot $r$. In the notations of \ref{hyperplane} 
choose $\alpha$ to be $\P_{L_1}(mr)$.

Let $q_1,\dots, q_n\in k^\times$ be a set of elements
such that $\displaystyle R_{w_0^L}^\rho \cdot( \prod ([\alpha]_\rho -
q_i)h)$ is regular 
for any  $\rho$, $L_1$, $L_2$, $L$, $w_0^L$
 as in \ref{hyperplane},  $\alpha$
as above and $h\in \End_{L_1/L_1^c}(\Pi_\rho)$
(such a set  exists by \ref{hyperplane}a)).

For $\lambda\in L_1/L_1^c$ 
denote 
\begin{equation}\label{Mxlambda}
M_{x,\lambda}^\rho= m_x^\rho\cdot
[\lambda]_\rho \prod _{r \in \Sigma^+-\Sigma^+_{L}}\prodl_{i=1}^n 
([\P_{L_1}(mr)]_\rho - q_i).
\end{equation}

\begin{Prop}\label{q1..qn}  Let $L_2=w(L_1)\ne L_1$, $L\supset L_1$
be standard Levi subgroups, $\rho \in Cusp_{L_1}$, and
$x\in  K\backslash G/ P$.

a) For any $\lambda \in L_1/L_1^c$ the element  
$R_w(M_{x,\lambda}^\rho)$ is regular
(i.e. lies in $\End(\Pi_{w(\rho)})$ rather than in  $\End(\Pi_{w(\rho)})
\otimes Frac (k[L_2/L_2^c])$).

b) If $L_2\subset L$ then  $R_w\(
M_{x,\lambda}^\rho\)
= M_{x,w(\lambda)}^{w(\rho)}$.
Moreover, if the integer $m$ used in the definition of $M_{x,\lambda}^\rho$
is large enough, then
$ M_{x,w(\lambda)}^{w(\rho)}-m_x[\lambda]_{w(\rho)}\in 
 \suml_{r\in \Sigma^+-\Sigma+_L}\Htil
_{\preceq \lambda+2n\delta^L-ar }$.

c) If $\lambda \in X _L$
is dominant enough and $L_2\not \subset L$, then 
 $R_w\(M_{x,\lambda}^\rho
% m_x^\rho\cdot [\lambda] \prod _{i,r \in \Sigma^+-\Sigma^+_{L}} ([\P_{L_1}
%(mr)] - q_i)
\)\in  \suml_{r\in \Sigma^+-\Sigma+_L}\Htil_{\preceq \lambda
+2n\delta^L-ar }$, 
where $\delta^L={1 \over 2} \suml _{r \in
\Sigma^+-\Sigma^+_L}r$.
\end{Prop}

\proof b) The first statement
 follows directly from Lemma \ref{2.5b} and \ref{actonweights}. 
To see the second notice that 
\begin{multline*}
\label{supprost}
\supp(R_w(M_{x,\lambda}^\rho-m_x[\lambda]_{w(\rho)}))=
 \{ \lambda +\sum\limits_{r \in  \P_{L_1}
(\Sigma^+-\Sigma^+_L)} i_r m r
  \mvert 0\leq i_r \leq n \}-\\
\{ \lambda +\sum\limits_{r \in  \P_{L_1}
(\Sigma^+-\Sigma^+_L)} n m r\}\subset \cupl_{r \in  \P_{L_1}
(\Sigma^+-\Sigma^+_L)}\a_{\preceq \lambda+2n\delta^L-m\P_{L_1}(r)}
\end{multline*}

(Here we used that  $\lambda,\,2\delta^L\in X_L\subset X_{L_1}$,
so $\P_{L_1}(\lambda+2n\delta^L)=\lambda+2n\delta^L$).
It is well-known that $\P_{L_1}$ preserves $\preceq$   (see e.g. \cite{BW}
Lemma 6.4 on p. 139, statement 2),
so the statement is clear.

Let us prove a). Recall that the set of standard Levi conjugate to $L_1$
is in bijection with chambers into which the root hyperplanes divide 
$\a_{L_1}$; the bijection sends a chamber $C$ to $L_C=w_C(L_1)$,
where $w_C\in W$ is an element such that $w(C)\in \a^+$
(such an element
 is defined uniquely up to right multiplication by an element of $W_L$;
we also have $w_C(C)=\a_{L_C}^+$).
Let $C_+= \a_L^+$, $C$ be the chambers corresponding
respectively to $L_1$, $L_2$. 
On the set of chambers we have the length function
$\ell(C)=\#\{r\in \P_{L_1}(\Sigma^+)\,|\, \<r,\alpha\> <0$ for $\alpha \in C\}$
(the latter condition will be abbreviated as $ \<r,C\> <0$).
We can choose a sequence of chambers $C_0=C_+,C_1,\dots,C_n=C$ so that
$C_i$ and $C_{i+1}$ have a common face of codimension 1, and $\ell
(C_i)=\ell(C_{i-1})+1$. Let $r_i\in \P_{L_1}(\Sigma^+)$ be the (only)
element such that $\<r_i, C_{i-1}\> >0$, $\<r_i, C_{i}\> <0$, and 
$L_i=L_{C_i}$. 
The sequence  $C_1,\dots C_n$ can be chosen  so as to satisfy the following
additional property: 
\begin{equation}\label{Lsnacha}
r_i\in\P_{L_1}(\Sigma ^+_L) \Longrightarrow r_j
\in\P_{L_1}(\Sigma ^+_L)\,\& \, L_j\subset L\operatorname{ \ for\  all\ } j<i.
\end{equation}
Indeed,  such a
 sequence  can be  constructed inductively as follows. Suppose that 
$C_1,\dots,C_i$ are already chosen,
and that $L_1,\dots,L_i\subset L,\,  r_1,\dots, r_i
\in\P_{L_1}(\Sigma ^+_L)$. 
Assume first that there exists a coroot $r\in \Sigma_L$
  separating $C_i$ from $C$ (i.e. $\<r,\alpha\>\cdot \<r,\beta\> <0$
for $\alpha \in C_i$, $\beta \in C$).  Then
there exists also a coroot $r\in \Sigma_L^+$ which
is orthogonal to a codimension 1 face of $C_i$,
and separates $C_i$ from $C$, for otherwise all simple roots of $L$
are nonnegative on the cone $w_{C_i}(C)$, while some positive root of $L$
is negative
there, which is impossible. Any (co)root separating $C_i$ from $C$
separates also $C_+$ from $C$; since $r\in \Sigma_L^+$ we get $\<r, C_+
\> >0$, $\<r, C_i\> > 0$.
  We now take $r_{i+1}= \P_{L_1}(r)
\in \P_{L_1}(\Sigma^+_L)$, and $C_{i+1}$ to be the chamber
neighboring with $C_i$ and separated from $C_i$ by $r'$.
It is clear that $L_{i+1}\subset L$.

If none of the roots separating $C_i$ from $C$ lies
in $\Sigma_L$, we choose the remaining $C_{i+1},\dots,C_n=C$ arbitrary
(keeping only the requirements $\dim (C_i\cap C_{i+1})=\dim (C_i)-1$,
$\ell(C_{i+1})=\ell(C_i)+1$). Since $r_i$ separates $C_j$ from $C$ 
for $j<i$ we see that the resulting sequence of chambers satisfies
\eqref{Lsnacha}.

Fix $(C_1,\dots,C_n)$ satisfying \eqref{Lsnacha}; let $l$
be the smallest integer for which $r_j\not \in \Sigma_L^+$. We can write: 

\begin{multline}\label{regular}
R_w(M_{x,\lambda}^\rho)=( R_{w_C \cdot w^{-1}_{C_{n-1}}} ) \circ 
( R_{w_{C_{n-1}} \cdot w^{-1}_{C_{n-2}}} )\circ \dots
 ( R_{w_{C_{l+1}} \cdot w^{-1}_{C_{l} } } )\circ
R_{w_{C_{l}}}(M_{x,\lambda}^\rho)   = \\
( R_{w_{C} \cdot w^{-1}_{C_{n-1}}} ) \circ 
( R_{w_{C_{n-1}} \cdot w^{-1}_{C_{n-2}}} )\circ \dots \circ
 ( R_{w_{C_{l+1}} \cdot w^{-1}_{C_{l}}} )(M^{w_{C_l}(\rho)}_{x,w_{C_l} 
(\lambda)})=
\\
\( R_{ w_{C_n}
 \cdot w^{-1}_{C_{n-1}}}\prod ([w_{C_{n-1}}(m r_n)]_{w_{C_{n-1}}(\rho)}
-q_j)  \)
\circ \dots 
\\
\circ   
 \( R_{w_{C_{l+1}} \cdot w^{-1}_{C_{l}}}\prod ([w_{C_{l}}(mr_{l+1})]_
{w_{C_{l}}(\rho) }-q_j) \) 
\\
\([w_{C_l}(\lambda)]_{w_{C_l}(\rho)  }m_x^{w_{C_l}(\rho)}          
\prodl_{r\in \P_{L_1}(\Sigma^+-\Sigma_L^+), r\ne r_i} \prodl_{j=1}^n        
 ([mr]_{w_{C_l}(\rho)  }-q_j)\).
\end{multline}

Notice that $ L_{i-1}$, $L_i$
generate  a standard Levi $L^i$ of rank equal to $\rank(L_i)+1$,
$w_{C_i}w_{C_{i-1}}^{-1}\in w_0^{L^i}\cdot W_{L_{i-1}}$ and
$w_{C_{i-1}}(r_i)$ is the only simple root  of $L^i$ which is not
a root of $L_{i-1}$. Hence, by the choice of $q_1,\dots,q_n$,
each of the operators
$\displaystyle
 R_{w_{C_{i}} \cdot w^{-1}_{C_{i-1}}} \prod ([w_{C_{i-1}}(mr_{i})]_
{w_{C_{i-1} } (\rho)}-q_j)$
is regular, and thus the whole expression  \eqref{regular} is regular.

Let us prove c).
We show that for $C_1,\dots, C_n$, $l$ as above and $i>l$ 
$$\mu\in \supp (w_{C_i}(M^\rho_{x,\lambda})) \Rightarrow
\exists \nu \in  \supp (w_{C_{i-1}}(M^\rho_{x,\lambda})), \,r\in \Sigma^+
-\Sigma^+_L:\  \mu \preceq \nu-ar.$$
 
Recall that $w_{C_{i-1}}w_{C_i}^{-1}\in w_0^{L^i} \cdot W_{L_{i}}$;
also the vector $\alpha:= w_{C_{i-1}}(r_i) \in X_{L_i}$ is orthogonal to
$X_{L^i}$.
Thus by Lemma \ref{hyperplane}b) there exists $\nu\in
 \supp (w_{C_{i-1}}(M^\rho_{x,\lambda}))$ such that 
$\nu=  w_0^{L^i} (\mu)+N \alpha$, where $|N|<n_0$.
Decompose $\nu$ as $\nu=\P_{L^i}(\nu)+\frac{\<\nu,\alpha\>}{\<\alpha,\alpha\>}
\alpha$.
% Then we have $ (w_0^{L^i} )^{-1}(\nu)= w_0^{L^i} (\nu)
%= \P_{L^i}(\nu)+\frac{\<\nu,\alpha\>}{\<\alpha,\alpha\>} w_0^{L^i} (\alpha)$. 

Let $\rtil_i\in \Sigma^+$ be a coroot such that $\P_{L_1}(\rtil_i)=r_i$.
Then  $w_{C_{i-1}}(\rtil_i) \succ 0 \succ w_{C_i} (\rtil)$,
because  (the dual root to)
$w_{C_{i-1}}(\rtil_i)$ is positive on $w_{C_{i-1}}(C_{i-1})\subset
\a^+$, while (the dual root to) $w_{C_i} (\rtil)$ is negative on 
 $w_{C_i}(C_i)\subset \a^+$.

Since the orthogonal projection $\P_L$ preserves $\preceq$,
%%%%%%%%%%for \wp this is Lemma 6.13 of {BW}
we have:
$\alpha=\P_{L_{i-1}}(w_{C_{i-1}}(\rtil_i))\succ 0\succ 
w_0^{L^i}(\alpha) =(w_0^{L^i})^{-1}(\alpha) 
 \P_{L_{i}}(w_{C_{i}}(\rtil_i))$. Moreover, since $\alpha$ belongs to
 a fixed finite subset of $\a_{L_1}^+$, there exists $b>0$ and $r\in \Sigma_+-
\Sigma_+^L$ such that $\alpha\succeq br$. 

To finish the proof it is enough to show that $\<\nu,\alpha\> \gg 0$
if $\lambda \in X_L$ is very dominant; in fact it is enough to ensure
that 
\begin{equation}\label{ocenka}
\<\nu,\alpha\> > \left( \frac{a}{b}+n_0 \right) \<\alpha, \alpha \>,
\end{equation}
 for then we get $\mu = (w_0^{L^i})^{-1}(\nu-N\alpha)=
\P_{L^i}(\nu)+\( \frac{\<\nu,\alpha\>}{\<\alpha,\alpha\>}-N\)
 w_0^{L^i}(\alpha)\prec \P_{L^i}(\nu) \prec
 \P_{L^i}(\nu)+\frac{\<\nu,\alpha\>}{\<\alpha,\alpha\>}\alpha \preceq \nu
-ar$.

To check \eqref{ocenka} notice that from \ref{hyperplane}b)
it follows by induction on $i$ that for $w=w_{C_i}$ 
there exists a finite set $S$ such that $\supp(R_w(M_{x,\lambda}^\rho))$
is contained in the convex hull of $w(\supp(M_{x,\lambda}^\rho))+S$
for any $\lambda, \rho$
(where we used the notation $A+B=\{a+b|a\in A, b\in B \}$). Hence 
$$\min\limits_{\nu\in \supp(R_{w_{C_{i-1}}}(M_{x,\lambda}^\rho))}
 ( \<\nu,\alpha\> )\leq \min\limits_{\nu\in w_{C_{i-1}}
(\supp(M_{x,\lambda}^\rho))}
 ( \<\nu,\alpha\> )+\const .$$
Since $\alpha = w_{C_{i-1}}(r_i)$ we have
$ \min\limits_{\nu\in w_{C_{i-1}}(\supp(M_{x,\lambda}^\rho))}
  \<\nu,\alpha\> = \min\limits_{\nu\in\supp(M_{x,\lambda}^\rho)}
  \<\nu,r_i\> \gg 0$ if $\nu\in X_L$ is very dominant, because
 our assumption $i>l$
implies that $r_i\in \Sigma^+-\Sigma^+_L$.
This proves \eqref{ocenka}.
\epf

\subsubsection{}
To deduce Proposition \ref{image} from \ref{q1..qn} we take 
\begin{equation}\label{hxlamrho}
h_{x,\lambda}^\rho:= 
%\suml_{L_1\subset L, \rho\in Cusp_{L_1}^K}
\suml_{L_2=w(L_1)}
%\frac{1}{\#\{L'\subset L|L'\in W_L(L_1)\}}
 R_w (M_{x,\lambda}^\rho).
\end{equation}

From Proposition \ref{q1..qn} it follows that $h_{x,\lambda}^\rho
\in \Htil$ and that
 $h_{x,\lambda}^\rho-\suml_{L_2=w(L_1)\subset L} [\lambda]_{w(\rho)}
 m_x^{w(\rho)}
\in \suml_{r\in\Sigma^+-\Sigma^+_L}
\Htil _{\preceq \lambda+2n\delta^L-ar}$. It is obvious that $h_{x,\lambda}^\rho
$ is invariant under all the intertwining operators $R_w$ (more precisely
$R_w(h_{\rho'})=h_{w(\rho')}$ where $h_{x,\lambda}^\rho
=\sum_{Cusp^K} h_{\rho'}$, $h_{\rho'}\in
\End (\Pi_{\rho'})$).
 Also it is immediate to check that $T_\psi(m_x^\rho)=
m_x^\rho$, $R_w\circ T_\psi\circ R_w^{-1}=T_{w(\psi)}$ when $T_\psi$ is
defined. 
If $\rho\otimes\psi\cong\psi$ for some irreducible representation $\rho$ of
a standard Levi $L$, then $\psi|_{\iota_L(X_L)}=1$, hence  for $\lambda\in X_L$
we have $T_\psi([\lambda]_\rho)=[\lambda]_{\rho\otimes\psi}$.

 This implies that each summand in \eqref{hxlamrho} is
invariant under $T_\psi$ when the latter is defined, provided $\lambda\in X_L$.
 By the matrix
 Payley-Wiener Theorem we conclude that  $h_{x,\lambda}^\rho\in I(\H)$.
Since $2\delta^L\in X_L$ the coweight $\lambda+2n\delta^L$ runs over the set
of very dominant elements of $X_L$ when $\lambda$ does.
Proposition \ref{image} is proved. \epf

\begin{Cor}\label{hxlam} Set
$$h_{x}^\lambda:= 
\suml_{L_1\subset L, \rho\in Cusp_{L_1}^K}
\frac{1}{\#\{L'\subset L|L'\in W_L(L_1)\}}h_{x,\lambda}^\rho.$$
Then we have $$I(h_x^\lambda)-[\lambda]m_x\in
  \suml_{r\in\Sigma^+-\Sigma^+_L}\Htil_{\preceq \lambda - a\cdot
r}$$
provided $\lambda\in X_L^+$ is dominant enough, and the integer $m$
used in \eqref{Mxlambda} is large. \epf
\end{Cor}

\subsubsection{} We next want to consider associated graded objects, so we
have to pass to a complete order.

Let $\leq$ be any complete order on $X   $ strengthening $\preceq$.
It yields filtrations on $\H,\,\Htil$ as in \ref{ograni1}--\ref{ograni3}.

We write $gr\H$ for the associated graded algebra of $\H,\, \leq$,
and use the letter $I$ for the induced imbedding
 $gr\H \hookrightarrow gr\Htil \cong \Htil$.

We have an obvious map $K_0\to gr_0\H$.

\begin{Thm}\label{tr0} a) 
 $gr\H$ is Noetherian.

b) Let $M$ be any finitely generated graded module over
 $gr\H$. Then 
for all but finitely many $\lambda \in X   $ we have 
$Tr(g, M_\lambda)=0$.
\end{Thm}

\proof Let  $Z$ be the center of $\H$.

The explicit description of $Z$ (Theorem 2.13 of \cite{BDKV}) implies 
 that the following elements form a basis in
 $I(grZ)\subset gr\Htil\cong\Htil$.
For a standard Levi $L$, $\rho\in Cusp_L^K$,
  $\lambda\in (L/L^c)^+$  such that
$\psi(\lambda)=1$ whenever $\psi\otimes \rho\cong \rho$
denote 
$z_\lambda^\rho:= \suml _{w(\lambda)=\lambda} [\lambda]_{w(\rho)\otimes
\psi_{w(\rho)}}$, where
$w$ runs over the set of such $w\in W/W_L$ for which $w(L)$ is a standard
parabolic and $w(\lambda)=\lambda$.

It follows that for $\lambda\in X^+$ we have $[\lambda]
%=\suml_{L,X_L\owns\lambda}\suml _{\rho\in Cusp_L^K} \frac{1}{\#\{L'=w(L),
%w(\lambda)=\lambda\} }z^\rho_\lambda
 \in grZ$. 

It is clear that $grZ$ is  central in $grH$, and
(as we will soon see) it is of finite type.

Thus  to prove a) it suffices to show that $gr\H$  is a finite $grZ$ module.
The main step  is the next:

\begin{Lem}\label{ogr}
There exists $\lambda_0\in \a   $ such that $gr\H_\lambda \ne 0\Rightarrow
\lambda \in \lambda _0+   \a^+$.
\end{Lem}

\proof of the Lemma. Assume the contrary. Then for any $N>0$
we can find  $\lambda\in X   $ such that $\<\lambda, r \><-N$
for a simple coroot $r$, and $grH_\lambda \ne 0$. Let $\lambda$ be
like that, then there exists
 $h\in \H$ such that $\lambda \in \supp (h) \subset X   _{\leq
\lambda}$. Decompose $h=\sum h_\rho$, $h_\rho\in \End(\Pi_\rho)$,
and fix $\rho\in Cusp_L^K$ for which  $\supp(h_\rho)\owns\lambda$.
Since $\<r,\lambda\>\ne 0$ we have $r\not \in \Sigma_L$. 
Let $L'$ be the standard Levi such that  $L'\supset L$,
$\rank(L')=\rank(L)+1$ and $r\in \Sigma_{L'}$. Then writing $\lambda=\P_{L'}
(\lambda)+\frac{\<\lambda,r\>}{\< \P_L(r),\P_L(r)\>}\P_L(r)$
we get by \ref{hyperplane}b) an element 
$\nu \in \supp(R_{w^{L'}_0}(h))=\supp(h)$ such that
 $\nu=w^{L'}_0(\lambda) +N'w^{L'}_0(\P_L(r))=\P_{L'}
(\lambda)+\( \frac{\<\lambda,r\>}{\< \P_L(r),\P_L(r)\>} +N' \)
w^{L'}_0(\P_L(r))$, where $|N'|<n_0$. 
We have $\P_L(r)\succ 0$, $w^{L'}_0(\P_L(r))=\P_{w^{L'}_0(L)}(w^{L'}_0(r))\prec 
0$, hence
taking $N>n_0\cdot \< \P_L(r),\P_L(r)\>$ we obtain
$\nu \succ \P_{L'}(\lambda)\succ \lambda$
 $\Rightarrow \nu > \lambda$.
This contradicts the condition $\supp(h)\subset \a_{\leq \lambda}$.\epf

\begin{Lem}\label{dur} For any $L$ fix $\lambda_L\in X  _L^+$. The subalgebra
$grZ_L\subset grZ$ 
generated by $[\lambda]$ for $\lambda\in \lambda_{L'}+ X   _{L'}^+$, 
$L'\subseteq L$
 is of finite
 type. The algebra 
${\cal Z}_L:=\oplusl _{L', \rho\in Cusp_{L'}^K}
\oplusl _{ \lambda\in (L'/(L')^c)\cap \a_L^+}
  k\cdot [\lambda]_\rho\subset \Htil$,
 where $L'$ runs over the set of standard Levi,
is a finite $grZ_L$-module.
\end{Lem}
\proof follows from the next easy fact. Fix $n\in \Zet_{\geq 0}$ and consider
the subsemigroup in $\Zet_{\geq 0}^m$
 generated by $\Zet_{\geq n}^I$, where $I$ runs over the subsets of $[1,m]$.
This semigroup is finitely generated, and $\Zet_{\geq 0}^m$ is the union of 
a finite number  of cosets of this semigroup. 
\epf

It is  clear that $ \bigoplus\limits_{\lambda\in
\lambda_0+\a   _+} gr\Htil_\lambda$ is a finite module over
${\cal Z}_G$. Hence by \ref{dur} it is also finite over 
 $grZ_G$. By \ref{ogr} we have $gr\H\subset \bigoplus\limits_{\lambda\in
\lambda_0+\a   _+} gr\Htil_\lambda$ for some $\lambda_0$ and
by \ref{dur} $grZ_G$ is Noetherian, hence $\grH$ is finite over 
 $grZ_G$. a) is proved.
Notice that along the way we have shown also that $gr\H$ is finite
over $grZ\supset grZ_G$.
 
Proof of b).
Consider a subalgebra ${\cal S} _L \subset gr\H$ generated by the highest
terms of the elements $h_x^{\lambda}$ (see \ref{hxlam}),
with large $\lambda\in X^+_L$, i.e. by the elements 
$[\lambda] m_x$.

Notice that 
${\cal S}_L\supset  grZ_L$ (for appropriate choice of coweights
$\lambda_L$ in \ref{dur}),
 because $\suml_{x\in K\backslash G/P}
 [\lambda]m_x=[\lambda]$. Hence $gr\H_L$ is finite over ${\cal S}_L
$, and thus any finite $gr\H_L$-module is also finite over ${\cal S}_L$. 

We prove by induction on corank of $L$, that for
any graded finitely generated  $\cal S$-module $M$ we have $Tr(g,
M_\lambda)=0$ for all but finitely many $\lambda \in X   _L$.
The step of induction is as follows. 

Let $M$ be an $X_L$-graded finitely generated ${\cal S}_L$-module. 

By the definition   ${\cal S}_L$ contains for some $\lambda_L\in X_L^+$
the subalgebra
$  % \oplusl_{\lambda \in \lambda_0+X   _L^+}\oplusl_{x,\rho/\sim}
  % k\cdot s_{x,\lambda}^\rho \cong
 \A_L [\lambda_0 +X   _L^+]
,$
 where
\begin{equation}\label{AL}
 \A_L:= 
\oplusl _{x\in  K\backslash G/P} k\cdot m_x.
\end{equation}
 Hence the proof of the Hilbert
Basis Theorem implies
 that for any graded finitely generated ${\cal S}_L$-module 
 $M$
 there exists  
$\mu_0$ and a finite  $\A_L$-module $M_0$
 such that
$M_\mu\cong M_0$ for $\mu \in \mu_0+X   _L^+$,
this isomorphism being compatible with the action.
In other words we have an isomorphism of graded 
$\A_L[\lambda_0 +X   _L^+]$-modules
$\oplusl _{\mu\in \mu_0+X_L^+} M_\mu
\cong M_0 [\mu_0 +X   _L^+]$. 

Further, for $\lambda\in \lambda_0+X   _L^+$,
$M_\lambda\cong M_0$ is the direct sum of the images of idempotents 
 $m_x\in \A_L$, and we have $g(\Im (m_x))=
\Im(m_{g(x)})$. We have assumed that $g$ acts on $
K\backslash G/P$ without fixed points, hence $tr(g,M_\lambda)=0$.

%This way we get a finite number of $\lambda_i\in \Lambda_L$ such
%that $tr(g,M_\lambda)=0$ for $\lambda\in \lambda_0+X   _L^+$, and such
It is obvious that
 $\Lambda-\cup (\lambda_0+  X   _L^+)$ is a  finite union
of cosets $\mu_i+X   _{L_i}^+$, $L_i\supset L$.
Then 
$ \bigoplus\limits_{\lambda \in
\mu_i+X   _{L_i}^+} M_\lambda$ 
is a finitely generated graded
module
over ${\cal S} _{L_i}$, because it is contained
in $ \bigoplus\limits_{\lambda \in
\mu_i+\a   _{L_i}^+} M_\lambda$ which is finite over
$grZ_\lambda$. Hence the statement for $M$ follows
from the induction hypothesis. 
The Theorem is proved. \epf

Now for any  finitely generated
graded module $N$ over $gr\H$ we can define
$Tr(g,N):=\sum\limits_\lambda Tr(g, N_\lambda)$ the sum being finite
by the last Theorem.

Let  $M$ be a finitely generated $\H$-module.

\subsubsection{}
We can choose an   $(\a   , \leq)$-filtration $F$ on $M$ compatible with
$\leq$ so that associated graded module $gr_F M$ is a finite
$gr\H$-module (such a filtration is called a good filtration).

\subsubsection{}\label{ei}
{\it From now on we assume that the complete order $\leq$ is induced by the
lexicographic order under the isomorphism  $\a    \cong \RE^r$,
$\lambda \mapsto \( (e_1, \lambda), ...,  (e_r, \lambda)\)$,
where  $(e_1,\dots , e _r)$ is a basis in $(X\otimes{\Q})^*$, 
such that $(e_i,r)\geq 0$ for $r\in \Sigma^+$.}

\begin{Prop}\label{Tr} a) $Tr(g,gr_F M )$ does not
 depend on $F$.
 
We set $Tr(g,M):= Tr(g,gr_F M)$.
Then $Tr(g,M)$ satisfies:

b) \ref{iii}i), and
% $Tr(g,M)$ is additive on short exact sequences (in $M$).

c)  \ref{iii}ii).

%If $\Mmod$ is an admissible $G$-module, and $\K$ is so small that $\chi_\Mmod$
%is constant on the coset $g\cdot K$, then 
%$$Tr(g,\Mmod)= \chi_\Mmod(g).$$
\end{Prop}

\proof Recall the following:
\begin{Fact}\label{fact} 
 Let $H$ be an algebra with an increasing $\Z$-filtration 
bounded below (i.e. $F_{\leq n}=0$ for $n\ll 0$),
such that the associated graded algebra is Noetherian;
let $M$ be a finite $H$-module. 

  Let $\grH$-Mod be the category whose objects are finitely generated
graded $\grH$-modules
and morphisms are homomorphisms of modules preserving the degree up to a shift
(i.e.  $Mor(M,N)=\sum Mor(M,N)_\nu$ where $Mor(M,N)_\nu=
\{\phi:M\to N\ |\ \phi (M_\lambda )\subset N_{\lambda+\nu}$).
 Then  the  class  of $grM$ in the Grothendieck group  $K^0(\grH$-Mod)
does not depend on the choice of a good filtration on $M$. 

Thus $[M]\to [grM]$ yields a well-defined homomorphism from $K^0(H$-Mod)
to $K^0(grH$-Mod).
\end{Fact}

\proof: see e.g. \cite{Gin}  Corollary 2.3.19 on p. 79. 
(The authors of {\it loc. cit.} consider the Grothendieck 
group of the category of  finitely generated $grH$-modules without grading;
however their argument actually proves the above statement as well).    \epf

We apply \ref{fact} inductively. 

Let the sublattice $\Lambda_0\subset \a$ be generated by $\Lambda_L$ for all
$L$. It is obvious that the set of indices $\lambda$ for which $gr\H_\lambda
\ne 0$ is contained in $\Lambda_0$.

For $i\in [1,..,n]$ let $\Lambda   ^{(i)}$ be the image of $\Lambda_0   $ 
in $\Q^i$
under $pr_i:\lambda \to (( \lambda, e _1) ,\dots , ( \lambda, e
_i))$; it is convenient to choose $e_i$ so that $\Lambda^{(i)}\subset \Zet^i$.
 Equip $\Lambda^{(i)}$  with the lexicographic order 
$\leq$ induced from $\Zet^i$.
We have a $\Lambda    ^{(i)}$-filtration 
$F^{(i)}_{\leq \overline \lambda}(\H)=\bigcup\limits_{\lambda, pr_i(
\lambda)= \overline \lambda}  \H_ {\leq \lambda}$ on $\H$; notice that
$F^{(i)}$ is induced by the corresponding $\Lambda   ^{(i)}$-filtration on 
$\Htil$.

 Furthermore, 
consider the  one parameter filtration $\Htil$ given by
$\Htil_{\leqi n} =\oplusl _{(\lambda,e_i)\leq n}\Htil_\lambda$, and the induced
filtration on $\H$ denoted by the same symbol. 

It yields a filtration on $gr_{F^{(i)}}\H$ given by the rule:
  $gr_{F^{(i)}}\H_{\leqi n}= \bigoplus \limits
_{\lambdabar \in X   ^{(i)}}[gr_{F^{(i)}}\H_\lambdabar]_{\leqi n}$
where $[gr_{F^{(i)}}\H_\lambdabar]_{\leqi n}= F^{(i)}_{\leq \overline
\lambda}(\H) \cap \H_{\leqi n}/  F^{(i)}_{< \overline
\lambda}(\H) \cap \H_{\leqi n}$. 

Notice also that this filtration on $gr_{F^{(i)}}\H$ is induced
by the $\leqi$-filtration on $\Htil \cong gr_{F^{(i)}}(\Htil)$
under the imbedding $gr_{F^{(i)}}(I)$.

We see that $gr_{F^{(i+1)}}\H =
gr_{\leqi}[gr_{F^{(i)}}\H ] $.

 If now  a good
$X   $-filtration $F$ 
on an $\H$-module $M$ is given, then again we can form  $\Lambda   
^{(i)}$-filtrations $F^{(i)}(M)$, where
 $F^{(i)}_{\leq \overline \lambda}(M)=\bigcup\limits_{\lambda, pr_i(
\lambda)= \overline \lambda}  M_ {\leq \lambda}$.
Consider the associated graded
$gr_{F^{(i)}}M$; on $gr_{F^{(i)}}M$ we have a one parameter filtration
denoted by $\leqi$, and given by
 $[gr_{F^{(i)}}M_\lambdabar]_{\leqi n} = \bigcup \limits
_{pr_i(\lambda)=\lambdabar, (e_{i+1},\lambda) \leq n} M_{\leq \lambda}/
 \bigcup \limits
_{pr_i(\lambda)<\lambdabar, (e_{i+1},\lambda) \leq n} M_{\leq \lambda}$.

Obviously, $gr_{F^{(i+1)}} M=gr_{\leqi}[gr_{F^{(i)}}M]$.

By inverse induction in $i$ we see that $\leqi$ is good;
hence using \ref{fact} and induction in $i$ we get the following
statement: 

{\it The class of $gr_F^{(i)}M$ in the Grothendieck group of the 
category of $\Lambda    ^{(i)}$-graded
$gr_F^{(i)}\H$-modules, with morphisms preserving grading up to a shift,
is independent on $F$.}

 (Actually to make a step of induction we need
a variant of \ref{fact} applicable in the situation when $H$ and $M$
have some multi-grading respected by the filtration while $grH$ and $grM$
are considered as 
multi-graded objects with one additional index. One may say that what
we need is exactly the statement \ref{fact} but with algebras/modules 
being algebras/modules in the category of multigraded vector spaces
rather than in the category of vector spaces).
% anyway the standard proof of
%\ref{fact}  extends to this situation with no difficulty). 

a), b) follow from the last italicized statement directly. c) is obvious:
if $\Mmod$ is an admissible representation, and $M=\Mmod^K$ then for large
$\lambda$ we have:
\begin{multline*}
Tr(g,M)=tr(g, F_{\leq \lambda}(M))=tr(g,M)=tr(\delta_K\cdot g, \Mmod)
=\\
\int _{g'\in K} \chi_\Mmod(g'g)dg'=\chi_\Mmod(g)\end{multline*}
provided the integrand is constant. \epf

\subsection{Geometric filtration and asymptotic cones}
In \ref{Tr} we have defined a functional $Tr(g,M)$ satisfying property
i), ii) of Theorem
\ref{iii}; it remains to prove that it satisfies also property iii). This
will be accomplished in the next two sections by passing from the ``spectral''
filtration constructed above to a ``geometric'' one defined
in terms of support of a distribution $h\in \H$; the latter filtration
is directly related to orbital integrals. 

The main results of this section used in further arguments are Theorem
\ref{soglasov}, Proposition \ref{F}, and Theorem \ref{upada}.

\subsubsection{}\label{2.1} We can assume that $K$ is {\it
 nice} with respect to $L_0$
(nice in the terminology of \cite{BDKV} 2.1b), i.e. that for any pair 
of opposite parabolics $P=L\cdot U,P^-=L\cdot U^-$,
with $L\supset L_0$ standard we have $K=K^+\cdot K^0\cdot K^-$
where  as usual $K^+=U\cap K$, $K^-=U^-\cap K$, $K^0=L\cap K$.
 (More precisely, 
the condition that $K_0$, $L_0$ are in good relative position implies
that   there exists a
 $K_0$-stable
 lattice ${\frak g}\subset \operatorname{Lie}\, G$ such that 
${\frak g}=({\frak g} \cap \operatorname{Lie}\,U^-) \oplus
({\frak g} \cap \operatorname{Lie}\,L)\oplus
({\frak g} \cap \operatorname{Lie}\,U)$;  we can take $\frak g$ to be
 the Lie algebra of the $O$-algebraic group attached to the
 $K_0$-fixed
point $p$  of the building  by the Bruhat-Tits theory \cite{BT2}.
%%%% [see Moy-Prasad p. 399].
 We then obtain $K$ from ${\frak p}^N \frak g$ for large $N$
 as in \cite{BDKV}2.1b)).  

Fix opposite parabolics $P=L\cdot U$, $P^-=L\cdot U^-$
with $P,L$ standard, and  set $C_P=(G/U\times G/U^-)/L$. 
By $\H(C_P)$ we denote the space of $K$-biinvariant compactly supported 
measures on $C_P$.

Consider the imbedding $\psi:L\hookrightarrow C_P$ given by  $\psi:l\mapsto
 (l\cdot
U, 1\cdot 
U^-) \mod L$.
 
Let $T\subseteq L_0$ be the (unique) maximal split torus contained in $L_0$.
For $\lambda\in X$ let $T_\lambda\in T/T^c$ be the coset of $\lambda$.
(Recall that we have an imbedding $X \imbed T$, $\chi \mapsto \chi ({\goth p})
$).

\begin{Prop}\label{bije} Let $\C\subset G$ be a compact $K$-biinvariant
set. 
 For $\lambda \in X    ^+$ such that   $\< \lambda, \Sigma^+-\Sigma^+_L
\> \gg 0$,  and $t \in T_\lambda$,
 there exists a unique bijection 
$\Psi_{\C,t}: K\backslash \C\cdot t \cdot  \C 
/K \to K\backslash   \C\cdot \psi(t) \cdot  \C /K$ which satisfies
$\Psi_{\C,t}(K \cdot c_1 t c_2 \cdot K)= K\cdot c_1 \psi(t)c_2 \cdot K$
for $c_1,c_2\in \C$.
\end{Prop}
\proof We have only to  check that 
\begin{multline*}
K \cdot g_1 t g_2 \cdot K=K \cdot g_1' t g_2' \cdot K
\Leftrightarrow K \cdot g_1 \psi(t) g_2 \cdot K
=K \cdot g_1' \psi(t) g_2' \cdot K \\
\operatorname{for\ } g_1,g_2,g_1',g_2'
\in \C.
\end{multline*}

Set $\C'=\{g_1g_2^{-1}|g_1,g_2\in \C\}\cup \{g_2g_1^{-1}|g_1,g_2\in \C\}$,
 and $K':= \capl_{g\in \C'} gKg^{-1}$ for $g\in \C$. 
It is clearly enough to verify that 
\begin{equation}\label{svert1}
g_1 \psi(t) g_2=\psi(t),\ g_1,g_2\in \C'
 \Longrightarrow   g_1 t g_2 \in K'\cdot t \cdot K' 
\end{equation}
\begin{equation}\label{svert2}
g_1 t g_2=t,\ g_1,g_2\in \C' \Longrightarrow     
g_1 \psi(t) g_2 \in K'\cdot \psi(t) \cdot K' 
\end{equation} 
If equality in the LHS of \eqref{svert1} holds then 
 we have: $g_1=lu$,
$g_2=t^{-1}l t u^-$ where $u\in U$, $u^-\in U^-$, $l\in L$. Since
$g_1,g_2$ lie in the compact set $ \C'$, the elements $l,\,t^{-1}l t,\,
  u,\,u^-$ are in some bounded subsets
of $L,\,U,\, U^-$. Hence the condition 
 $\<\lambda,  \Sigma^+-\Sigma^+_L\> \gg 0$
implies that
 $(tg_2)^{-1}u (t g_2) \in K$, $(lt)u^-(lt)^{-1}\in K$,
which yields $g_1 t g_2^{-1}=[(lt)u_-(lt)^{-1}]^{-1}t[t^{-1}u t]\in K t K$.

If equality in the LHS of \eqref{svert2} holds, 
then 
$g_1=tg_2t^{-1}\in \C' \cap t \C' t^{-1}$. It is well-known
%%%%%%%%%%% \cite??
that for $t\in T_\mu$, $\<\mu,  \Sigma^+-\Sigma^+_L\> \gg 0$
 we  have $ \C' \cap t \C' t^{-1} \subset \C^- \cdot  \C^L
\cdot \C^+$, where $\C^+$, $\C^L$ are bounded (independently of $t$)
subsets of respectively  $U$ and $L$, and $\C^-$ is an arbitrary small
subset of $U^-$; in particular we can assume $\C^-\subset (K')^-$.
 Thus we have
$g_1=u_-\cdot l \cdot u$  with $u_-\in K^-$, $l\in \C^L$, $u\in \C^+$.
Then $g_2= t^{-1}g_1 t $,  and we can assume $  t^{-1}u t\in (K')^+ $.
So finally we get $(g_1, g_2)(\psi(t))= (u_-,  t^{-1}u t)(\psi(t))
\in K'\times K' (\psi(t))$. \epf

 If $t\in T_\mu$ where 
$\< \mu,\alpha \>$ is so large for $\alpha \in \Sigma^+-\Sigma_L^+$
 that $\Psi_{\C'', t}$ is defined
 for $t\in T_\mu$ and  $\C'':=\{g_1g_2^{-1}g_3\,|\,g_1,g_2,g_3\in \C\}$, then 
obviously $\Psi_{\C,t}$ and $\Psi_{\C,t'}$ agree on the intersection 
$\C t\C\cap \C t'\C$ 
whenever this intersection is nonempty. Thus we have a bijection
$\Psi:  K\backslash (\cupl_ 
{\<\lambda,  \Sigma^+-\Sigma^+_L\> \gg 0} \C\cdot T_\lambda \cdot  \C) 
/K \to K\backslash (\cupl_{\<\lambda,  \Sigma^+-\Sigma^+_L\> \gg 0}
  \C\cdot \psi(
T_\lambda) \cdot  \C) /K$.
%, where $t$ runs over the union of $T_\mu$ for $\<\mu,\alpha\>\gg 0 $ 
%for  $\alpha \in \Sigma^+-\Sigma_L^+$.

% for which $t\in T\_mu \Rightarrow  \Psi_{\C'',t}$ is defined. 

\subsubsection{} \label{Cartan}
The map $w\mapsto K_0 \cdot w\cdot K_0$ yields a canonical bijection 
$$ W_p\backslash W_{aff}/ W_p\iso K_0\backslash
G/K_0,$$ 
see \cite{BT1} 7.4.15 (notice that 
since $K_0$ and $L_0$ are in good relative position we have $K_0\supset 
L_0^c$, so the expression $ K_0 \cdot w\cdot K_0$ is meaningful).

Consider the projection $pr_L: W_p\backslash W_{aff}/ W_p
\to W_f^p\backslash W_{aff}^p/ W_f^p=\Lambda/W_f^p=\Lambda^+$.
For $x\in  W_p\backslash W_{aff}/ W_p$ we write $G_x$ for the corresponding
double coset in $G$, $\H_x$ for the space of elements in $\H$ supported
in $G_x$, and we set $G_\lambda:=\cupl_{pr_L(x)=\lambda}G_x$,
$\H_\lambda:=\oplusl_{pr_L(x)=\lambda} \H_x$.

 (Support of an element $h\in \H$
appearing here should not be confused with  support of
the corresponding element $I(h)\in \Htil$ considered before (see
\eqref{supp^spec}): the first
one is an open compact subset of $G$ while the second is a finite 
subset of $X   $.
To rule out the very possibility of such confusion we will sometimes 
refer to the first notion as to the {\it geometric support}, and to
the second one as to the {\it spectral support} and write respectively 
$\supp^{spec}$ and $\supp^{geom}$).

 Let us choose 
$\C$ to be large enough; more precisely, we require the following.
We can choose a finite number of elements $w_i\in W_{aff}$, so that for
any $w\in W_{aff}$ there exist $i_1, i_2$ such that $w_{i_1}ww_{i_2}\in X^+$.
Let us take $\C$ to be $\cupl_i (K_0\cdot w_{i})  \cup \cupl_i (w_{i} K_0)$.
Then we have
$\cupl_{\<\lambda,  \Sigma^+-\Sigma^+_L\> \gg 0}  \C\cdot T_\lambda \cdot  \C
\supset \cupl _{\mu\in \a^+,\<\mu,  \Sigma^+-\Sigma^+_L\> \gg 0}
G_\mu$.

Denote $(C_P)_\mu:=\cupl _{g\in G_\mu} \Psi(KgK)$,
$\H(C_P)_\mu=\{f\in \H (C_P) \,|\, \supp \  f \subset (C_P)_\mu\}.$

Thus we have a bijection
 $$\cupl _{\mu\in \a^+,\<\mu,  \Sigma^+-\Sigma^+_L\> \gg 0}
K \backslash G_\mu /K\iso 
\cupl _{\mu\in \a^+,\<\mu,  \Sigma^+-\Sigma^+_L\> \gg 0}
K \backslash C_\mu /K.$$

Let us  set $$\H^P_N=
 \oplusl_{\<\mu,\alpha\>\geq N, \alpha \in \Sigma^+-\Sigma^+_L} \H_\mu ,\ \
 \H(C_P)_N= \oplusl_{\<\mu,\alpha\>\geq N, \alpha \in \Sigma^+-\Sigma^+_L}
 \H(C_P)_  \mu,$$
%\{f\in C_c^\infty (C_P)^{K\times K} \,|\, \supp \  f \subset 
% \bigcup  \limits _{\<\mu,\alpha\>\geq N, \alpha \in \Sigma^+-\Sigma^+_L} 
%\Psi(G^\mu) \}$,
 and use the notation  $\Psi$ for the isomorphism $\H^P_N\iso \H(C_P)_N$
sending the $\delta$-function of a double $K$-coset $x$ on $G$ to  
the  $\delta$-function of $\Psi(x)$. (Here $N\gg 0$).

\begin{Lem}\label{odinak}
 Let $\C\subset G$ be compact. 
Let  $N$ be so large that  $\Psi_{\C\cdot \C,t}$ is
 defined for $\<\mu,\Sigma^+
-\Sigma_L^+\> \geq N$, $t\in T_\mu$ where $\C\cdot \C=\{c_1c_2|\,c_1,c_2\in \C
\}
.$
Then we have
  $\supp(h)\subset \C,\, h'\in \H^P_N\Rightarrow \Psi(h*h'), \, \Psi(h'*h)$
are defined and 
$$ \Psi(h*h')=h*\Psi(h'),$$
$$\Psi(h'*h)=\Psi(h')*h$$. 
\end{Lem}
\proof follows directly from \ref{bije}. \epf

\begin{Rem} {\it Asymptotic semigroups.}
Before going on with the argument let us explain the geometric picture standing
behind it. 

Recall from \cite{Vinberg1} that for a semisimple algebraic group 
$\underline{G}$ there exists a canonical way to degenerate
 $\underline{G}$ into one of its asymptotic cones
 $\overline{\underline C}_P$
isomorphic to the affine closure of the algebraic
variety $\underline{(G/U\times G/U^-)/L}$. 

More precisely, in \cite{Vinberg1} Vinberg constructs
the so called   {\it universal semigroup}
$\Gtil$ enjoying the following properties. (He works over an algebraically
closed field of characteristic 0; his construction carries over immediately
to the case of a split
group over a field of characteristic 0).

 Let $R$ be the set of simple 
roots of $G$; let $\Tbar \cong ({\Bbb A}^1)^R$ be the canonical partial
compactification of the abstract Cartan group of the adjoint group $G/A$
(where $A$ is the finite center of $G$). Notice that $\Tbar$ has a natural
``coordinate cross'' stratification, with the set of strata being in bijection
with subsets of $R$, i.e. with conjugacy classes of parabolics.
For a parabolic $P$ let $T_P\subset \Tbar$ be the corresponding stratum, and
 let $e_P\in T_P$ be the point all of whose coordinates are either 0 or 1
(so the coordinate of $e_P$ corresponding to a (co)root $r\in \Sigma_L$ is 0,
and the one corresponding to a (co)root in $\Sigma-\Sigma_L$ is 1).

$\Gtil $ is an algebraic
 semigroup equipped with a flat morphism $\wp: \Gtil\to \Tbar$,
 such that the preimage
 of the open stratum is isomorphic to $G\times_A \Tbar$. 
Further, for $x\in T_P$ the preimage $\wp^{-1}(x)$ is isomorphic to
 the algebraic
variety $\underline{(G/U\times G/U^-)/L}$. 

Now the canonical bijection $\Psi$ can be characterized as follows.

For any parabolic $P$ there exists an $A_L$-invariant open neighborhood 
$\U_P$ of
$\wp^{-1}(e_P)$, and a continous $A_L$-equivariant projection $\pi_P:\U_P\to 
\wp^{-1}(e_P)$, such that $\pi_P|_{\wp^{-1}(e_P)}=Id$, and the following
holds.

Whenever $KgK\subset \U_P$ and $KxK\subset \U_P$ for $g\in G$ and 
$x\in C_{P'}$, where $P'\supset P$, we have $KxK=\Psi (KgK)\Leftrightarrow
\pi_P(KxK)=\pi_P(KgK)$. In other words $KxK$ and $KgK$ are close in the 
topology of $\Gtil$. 

Notice also that the domain of definition of $\Psi: K\backslash G/K
\to K\backslash (C_P)/K$
 is a neighborhood of the stratum corresponding to 
$P$ in the DeConcini-Procesi compactification of $G$ \cite{DPr}.

Probably it is possible to work directly with this geometric definition;
we however found it more economical (though less transparent) to use the 
techniques employed here.
\end{Rem}

\begin{Lem}\label{tymnenadoel} Suppose that $l\in L_\mu$, where 
$\<\mu,\Sigma^+-\Sigma^+_L\>\gg 0$.
Then $\Psi(KlK)=K\psi(l)K$.
\end{Lem}

\proof Obvious. \epf

For a partial order $\treugol$ on $\a   $
we write $G_{\treugol ?}= \bigcup  \limits _{\mu \treugol ?} G_\mu$,
$\G_{\treugol ?}(\H)= \oplusl \limits _{\mu \treugol ?} \H_\mu$.

We are now ready for the proof of the next

\begin{Thm}\label{soglasov}
 For dominant enough $\lambda$ we have
$\H_{\preceq \mu} \cdot \G_{\preceq \lambda}(\H)\subset
\G_{\preceq \lambda+\mu}(\H)$.
\end{Thm} 

\begin{Lem}\label{tautolog}
 For $h\in \H$ and $\lambda \in \a   ^+$ the following are
equivalent: 

i) $h\in \H_{\preceq \lambda}$.

ii) For any parabolic, and any function 
 $f\in C_c^\infty (G/U)^K$ we have $\supp (f)\subset  (G/U)_\mu
\Rightarrow \supp (h*f) \subset \bigcup \limits_{\nu \preceq \lambda+\mu}
(G/U)_\nu$. 

\end{Lem}
\proof Assume that ii) holds. Recall that for any $\rho\in Cusp^K$ we have
$\Pi_\rho^K= C_c^\infty (G/U)^K\otimes _{\H(L^c)}\rho$, and the grading
on $\Pi_\rho$ comes from the grading on the first multiple in the RHS, 
given by $ C_c^\infty (G/U)^K=\oplusl_\mu C_c^\infty (G/U)_\mu^K$.
ii) $\Rightarrow$ i) is now obvious. 

Conversely, suppose that ii) does not hold, i.e.
 for some $\nu\not \preceq \lambda$ there 
exists $f\in  C_c^\infty (G/U)^K$ such that $ \supp (h*f)\cap 
(G/U)_{\lambda+\nu}\ne \emptyset$. Thus the morphism of free 
$\H(L^c)$-modules $ C_c^\infty (G/U)^K_\mu\to  C_c^\infty (G/U)^K_{\mu+
\nu}$, $f\mapsto h*f|_{ (G/U)_{\mu+
\nu}}$ is nonzero. Hence (recall that $\H(L^c)$ is finite over its center)
there exists an irreducible representation $\rho$ of $L^c$ such that 
tensoring the latter morphism over $\H(L)$ with $\rho$ we still obtain a 
nonzero morphism. Further, $\rho$ is a subquotient in a representation
of the form $i_{L_1}^L(\rho_1)|_{L^c},$ where $\rho_1$ is a cuspidal
 representation of a standard Levi $L_1\subset L$. Then there exists
$\nu_1\in \a  _{L_1}$ with $\P_{L}(\nu_1)=\nu$ such that the morphism  
$ C_c^\infty (G/U_1)^K_\mu\to  C_c^\infty (G/U)^K_{\mu+
\nu_1}$,  $f\mapsto h*f|_{ (G/U)_{\mu+
\nu_1}}$ tensored with $\rho_1$ is nonzero. This means that 
$\supp^{spec}(I_{\rho_1'}(h))\owns \nu_1$, where $\rho_1'\in Cusp$
differs from $\rho$ by twisting with an unramified character;
so $\supp^{spec}(h)\owns \nu_1$.

If i) holds then $\nu_1\preceq \lambda$. 
Since $\lambda\in \a^+$ we have $\P(\lambda)\preceq \lambda$;
since $\P_L$ preserves $\preceq$   (by \cite{BW},
Lemma 6.4 on p.139, statement 2)
  we get $\lambda\succeq \P_L(\lambda)\succeq \nu$.
We arrived at a contradiction, which proves i) $\Rightarrow$  ii). 
 \epf 

\begin{Lem}\label{nossogl} Let $\C\subset G$ be a compact set. 
Suppose that $x_1\in (C_P)_{\mu_1},$ $x_2\in (C_P)_{\mu_2},$
where $\<\mu_1,\Sigma^+-\Sigma^+_L\>\gg 0$,   
and 
$x_2=cx_1$ for $c\in \C$. Identify the left $G$-orbit on $C_P$
passing through $x_1$, $x_2$  with $G/U$,
and suppose that $x_1\in (G/U)_{\nu_1}$, $x_2\in (G/U)_{\nu_2}$.
Then $\nu_2-\nu_1=\P_L(\mu_2-\mu_1)$.
\end{Lem}

\proof 
We have a map 
$$
W_{aff}\times _{W_{aff}(L)} W_{aff}\to K_0\backslash C_P/K_0
$$
which sends $(w_1,w_2)$  to $K_0\cdot (w_1,w_2)\cdot K_0$ (it is immediate to
check that it is actually well-defined).
%%%%%%%%%%; it follows from \ref{bijection}
%%%%%%%%%%and the first paragraph of \ref{Cartan} applied to $L$). 
For $x\in W_{aff}\times_{W_{aff}(L)} W_{aff} $  denote the corresponding
double coset by $(C_P)_x$.  

As in \ref{Cartan} we pick a finite set $\{w_i\}\in W_{aff}$ so that for
any $w\in W_{aff}$ we have $w_{i}ww_{j}\in X^+$ for some $i,j$. 
For $\<\mu,  \Sigma^+-\Sigma^+_L\> \gg 0$  
Proposition \ref{bije} implies that
%%%%%%%%%%\begin{equation}
  %%%%%%%%%%\label{C_P_x}
$
(C_P)_{(w_{i}\mu, w_{j})} =\cupl_{g\in G_{w_i\mu w_j}} \Psi (KgK).
$
%%%%%%%%%%\end{equation}

In particular, $x_i\in (C_P)_{w_i \lambda_i w_i'}$ for $i=1,2$,
where $\lambda_i$ is such that $\<
\lambda_i, \Sigma^+-\Sigma^+_L\>\gg 0$, and $w_i\lambda_i w_i'\in
 W_f^p\cdot \mu_i\cdot W_f^p$. 

Since $x_1$, $x_2$ lie in one left $G$-orbit we have $(\lambda_2 w_2')L_0^c
\subset P^- (\lambda_1 w_1') L_0^c$, which is possible only if $w_2'\in
W_{aff}(L)w_2 $. Thus we can (and will) assume that $w_1'=w_2'=w'$. 

Pick $x\in (C_P)_{ w'}\cap G(x_1)$, and identify $G/U$ with the left orbit
$G(x_1)=G(x_2)$ by means of the map $g\mapsto g(x)$. Then 
we get $x_i\in (K_0\cdot (w_i\lambda_i) \cdot (K_0\cap L)U)/U
\subset (G/U)_{\overline {w_i\lambda_i}}$, where $\overline w$ is the image
of an element $w\in W_{aff}$ under the projection $W_{aff}\to
W_f^p \backslash W_{aff}/W_{aff}^L$. 

Thus  $\nu_i=\P_L(\tilde \nu_i)$, where $ w_i\lambda_i
\in W_f^p\cdot \tilde \nu_i$, while $W_f^p \mu_i  W_f^p\owns w_i\lambda_i
w'$.
Hence $\mu_i\in W(\tilde \nu _i+\mu')$, where $\mu'$ is such that $w'\in
\mu'\cdot W_f^p$. 
Notice that since $\< \mu_i,\Sigma^+-\Sigma^+_L\>\gg 0$, and $w',\mu'$
 belong
to a fixed finite set, we have $ \<\tilde\nu_i
+ \mu',\Sigma^+-\Sigma^+_L\>\gg 0$. Since $\mu_i$ is dominant, it follows
that $\mu_i\in W_L(\tilde \nu_i+\mu')$. 
In particular $\P_L(\mu_i)=\P_L(\tilde \nu_i+\mu')=\nu_i+\P_L(\mu')$.
The Lemma follows.
\epf

\proof of  Theorem \ref{soglasov}.
Consider the ``Rees algebra'' $Rees(\H):=\oplusl _{\lambda\in X   ^+}
\H_{\preceq \lambda}$. The explicit description
of the center $Z\subset \H$ (\cite{BDKV}, Theorem 2.13)
readily implies that the central subalgebra $Rees(Z)\subset Rees(\H)$ is
of finite type.  Further, from Lemma \ref{hyperplane}b) it is not hard to
deduce that there exists $\lambda_0\in X   ^+$ such that for 
$\lambda\in
\a_+$, 
 $\rho\in Cusp_L^K$ we have $h\in \H_{\preceq \lambda},\nu \in \supp^{spec}(I_\rho(h))
\Rightarrow w(\nu)\preceq \lambda+\lambda_0$ whenever $w(L)$ is a standard
 Levi.

Thus we have $Rees (\H) \imbed \oplusl_{L,\rho\in Cusp_L^K}
\oplusl_{\lambda\in\a_+} 
\capl_w 
%\oplusl _{\nu\in \cap_w w(\a_{\preceq\lambda+\lambda_0})}
% \limits_{X   ^+\times Cusp^\K}
\End(\Pi_\rho)_{w(\a_{\preceq \lambda+\lambda_0} ) }$.
 The RHS is a finite $Rees(Z)$-module. Hence the LHS is, thus $Rees(\H)$
is a Noetherian graded algebra, and in particular is finitely generated.

 Let $h_i\in \H_{\preceq \lambda_i},$
$i=(1,..,n)$ be the generators. Then it is enough to ensure that 
$f\in \G_{\preceq \mu}(\H)\Rightarrow h_i f\in \G_{\preceq \lambda_i+\mu}(\H)$.

Pick a number $N$ which satisfies the condition of \ref{odinak} with
$\C=\cup \supp (h_i)$ for all $P$. 

Also let $N'$ be such that  $h_i\in \G_{\la
N'}(\H)$ for all $i,\alpha$.

Suppose that  $\lambda \in \a   ^+$ is large enough, $\mu\preceq \lambda$,
 and $h\in \H_\mu$.

Pick a simple
root $\alpha$. We consider two cases.

I.
$\<\mu,\alpha \> > N$. Then let $P$ be the maximal proper parabolic
corresponding to $\alpha$. Applying \ref{odinak} we see that $h_i*h=
\Psi^{-1}(h_i*\Psi(h))$. 

We can decompose
 $ C_P=\bigcup \limits_{x\in \Pbar\backslash G/K} (C_P)_x$,
where  $(C_P)_x$ is the fiber of the projection $C_P\to  \Pbar\backslash G/K$,
$(g_1 U, g_2 U^-)L\mapsto \Pbar g_2K$.

Let us now write $h$ as the sum of $h_x$, where $\supp(\Psi(h_x))\in (C_P)_x$.

We have an isomorphism of $G$-modules $C_c^\infty((C_P)_x))^K\cong C_c^\infty
(G/U)^{K\cap P}$. (Here the invariants are taken with respect to the right 
action).

Hence applying to each $h_x$ \ref{tautolog} and \ref{nossogl} we obtain:
$$h_i*h \in \G_{\la \< \mu, \omega_\alpha \>+\< \lambda_i, \omega_\alpha \>}
(\H)\subset  \G_{\la \< \lambda,
 \omega_\alpha \>+\< \lambda_i, \omega_\alpha \>}
(\H).$$

II. If on the other hand $\<\mu,\alpha \> \leq N$ then 

\beq\label{rasstojanie}
\< \lambda -\mu , \omega _\alpha \> \geq { (\lambda, \alpha)-N
 \over (\alpha,\alpha) },
\end{equation}
because the RHS is the distance from $\lambda$ to the hyperplane
$\{\eta\ |\ \<\eta,\omega_\alpha\>=N\},$ while $\mu$ is separated
from $\lambda$ by this hyperplane.

If $\lambda$ is large then the RHS of 
\eqref{rasstojanie} is greater than $N'$,
so in that case
$$h_i*h \in \G_{\la \< \mu, \omega_\alpha \>  + N'} (\H)\subset
 \G_{\la \lambda}(\H)\subset  \G_{\la \lambda+ \lambda_i}(\H)
.$$

Thus in any case $h_i *\G_{\preceq \lambda}(\H)\subset \G_{\preceq
\lambda+\lambda_i}(\H)$ provided  $\lambda$ is large.

Since $h_i$ generate $Rees(\H)$ this finishes the proof. \epf

To derive an important corollary we need some information on compatibility
of our filtrations with convolution.

\begin{Lem}\label{build} We have 

a) $\G_{\preceq \nu_1}(\H)\cdot \G_{\preceq
\nu_1}(\H) \subseteq \G_{\preceq \nu_1+\nu_2}(\H)$.

b) $G_\lambda((G/U)_\mu) \subset \cupl_{\nu\preceq
\lambda+\mu}((G/U)_{\nu})$.
\end{Lem}

\proof
Recall that we have fixed a minimal Levi subgroup in good relative position
with $K_0$. This choice determines  an appartment $\Ag$ of the Bruhat-Tits
building $X$ containing the $K_0$-fixed point $p$.
Let $\Xi:X\to \Ag$ be the contraction centered at
an open polysimplex, which contains $p$ in its closure. Then an element $g\in
G$ lies in $G_\nu$ iff the coweight $\Xi (g(p))-p \in \a$ is
$W$-conjugate to $\nu$.
Now the proof is parallel to the usual proof of the fact that $\Xi$
 does not increase distances (see \cite{BT1} 7.4.20(ii)).

 Take now $g_1\in G_{\nu_1}$, $g_2\in G_{\nu_2}$; we must check that
$g_1g_2\in G_\nu$ where  $\nu$ satisfies $\nu\preceq \nu_1+\nu_2$.

Let us break the segment $ [g_1(p),g_1g_2(p))$ into the union 
$[z_0=g_1(p),z_1)$, $[z_1,z_2)$, $\dots[z_{n-1},g_1g_2(p)=z_n)$, where
$z_i$ and $z_{i+1}$  lie in the closure of one open polysimplex ({\it
loc. cit.} 7.4.21). 

According to \cite{BT1} 7.4.19, 7.4.18(i) the map $\Xi$ restricted to
the closure of any polysimplex coincides with the action of some $g\in G$,
hence 
we have $\Xi(z_i)=x_i(z_i)$, $\Xi(z_{i+1})=x_i(z_{i+1})$ for some $x_i\in
G$.

On the other hand $\Xi([p,g_2(p)])=g([p,g_2(p)])$ for some $g\in G$ is a line 
segment, and by the above remark $w(\Xi(g_2(p))-p)=\nu_2$ for some $w\in
W$.

Then obviously $\nu_2=\suml_{i=0}^{n-1} \nu^i$, where $\nu^i=
w\( \Xi(g_2(z_i))-\Xi(g_2(z_{i-1})) \)\in \a^+$.

By \cite{BT1} 7.4.8  for points $p_1,p_2,p_1',p_2'\in \Ag$ the elements
$p_1-p_2$ and $p_1'-p_2' \in \a$ are $W$-conjugate provided $p_i=g(p_i')$
for some $g\in G$. 

Thus we see that $\Xi(z_i)-\Xi(z_{i-1})$ is $W$-conjugate to 
 $\Xi(g_2(z_i))-\Xi(g_2(z_{i-1}))
=gg_1x_i^{-1}(\Xi(z_i))- gg_1x_i^{-1}( \Xi(z_{i-1}))$, and hence also to
 $\nu^i$. Let $w_i\in W$ be such that $\Xi(z_i)-\Xi(z_{i-1})=w_i(\nu^i)$.
Suppose that
$g_1g_2\in G_\nu$, i.e. $w'(\Xi(g_1g_2(p)-p)=\nu$ for some $w'\in
 W$.
Then we have 
\begin{multline*}
\nu=w'(\Xi(g_1g_2(p)-p)=w'(\nu_1+\sum w_i(\nu^i))=w'(\nu_1)+\sum
w'w_i(\nu^i) \preceq\\
 \nu+\sum \nu^i=\nu_1+\nu_2. \end{multline*}
This proves a).

To prove b) we reinterpret the condition $g U \in (G/U)_\mu$ as 
follows. Let $\Xi':X\to \Ag$ be the contraction centered at a ``vector
chamber''  corresponding to $P_0$. Then $gU\in (G/U)_\mu$ if and only if
$\P_L(p- \Xi'(g^{-1}(p)))= \nu$.

From the definition (\cite{BT1} 7.4.25) it is clear that
 $\Xi'|_{\Deltabar}=g_{\Delta}$ for any polysimplex for 
 some $g_\Delta\in G$. Hence repeating the arguments used in a) we can find
 $\lambda^i\in \a^+$ and $w_i\in W$, such that
 $\Xi(g_2^{-1}(p))-\Xi(g_2^{-1}g_1^{-1}) =\sum w_i(\lambda^i)$, while
 $\lambda =\sum \lambda^i$. Since by \cite{BW},
Lemma 6.4 on p.139 $\P_L$ preseves $\preceq$ we see that if $
g_1\in G_\lambda$, $g_2 U\in (G/U)_\mu$, $g_1g_2\in (G/U)_{\mu'}$ 
then 
\begin{multline*}
\mu'=\P_L(p-\Xi'(g_2^{-1}(p))) + \P_L(
\Xi(g_2^{-1}(p))-\Xi(g_2^{-1}g_1^{-1})
)=  \\
\P_L(\mu)+\P_L(\sum w_i(\lambda^i))\preceq 
 \mu +\sum
 \lambda_i=\mu+\lambda.
\end{multline*}
This proves the Lemma. \epf

\begin{Cor}\label{easydir} For any $\lambda \in X   ^+$ we have
$\H_{\preceq \lambda}  
\supset  \G_{\preceq \lambda}(\H)$.
\end{Cor}

\proof Compare  Lemma \ref{tautolog} with Lemma \ref{build}b). \epf

The latter statement together with Theorem \ref{soglasov} yield the
following 

\begin{Prop}\label{F} Fix $\lambda_0$ in $X   ^+$, and define filtration
$F_{\leq \lambda}(\H):= \suml_{\mu \in \a^+_{\leq\lambda}} \G_{\preceq
\lambda_0+\mu}(\H)$. Then $F_{\leq}$ is a good filtration on the free
$\H$-module, provided $\lambda_0$ is large enough.
\end{Prop}

\proof From Theorem \ref{soglasov} we see that the filtrations on the
algebra and on the module are compatible if $\lambda_0$ is large. 

It remains to show that for all but finitely many $\nu$ we have 
\begin{equation}\label{nado1}
F_{\leq \lambda}\subset \suml_{\lambda_2\in \a^+_{<\lambda}} 
\suml_{\lambda_1\in \a_{\leq \lambda-\lambda_2}}
\H_{\leq
\lambda_1} \cdot F_{\leq \lambda_2}.
\end{equation}
We will in fact show that    
\begin{equation}\label{nado2}
F_{\leq \lambda}\subset  \suml_{\lambda_2\in \a^+_{<\lambda}} 
\suml_{\lambda_1\in \a_{\leq \lambda-\lambda_2}}
 \G_{\preceq
\lambda_1}(\H) \cdot F_{\leq \lambda_2}
\end{equation}
 for almost all
$\lambda$; \eqref{nado2} is stronger than \eqref{nado1} by 
\ref{easydir}.

Pick some $\lambda$ and $\nu $ such that $\H_\nu\subset F_{\leq
\lambda}(\H)$. We want to check that $\H_\nu$ lies in the RHS of
\eqref{nado2}.

 By the definition   $\nu \preceq
\lambda'+\lambda_0$ for some $\lambda'\in \a^+_{\leq \lambda}$. 
If $\lambda'< \lambda$ then $\H_\nu$ does lie in the RHS of \eqref{nado2},
 for it lies in the summand corresponding to  $\lambda_2=\lambda'$,
$\lambda_1=0$. So assume that $\lambda'=\lambda\succeq\nu$.

For any $w\in W_f^p$ consider the intersection $w^{-1}W_{aff}w\cap
\Lambda$. It is a sublattice of finite index in $\Lambda$;  in
particular there exists some $n\in \Zet^{>0}$ such that $n\omega_\alpha\in
 w^{-1}W_{aff}w\cap \Lambda$ for any $w\in W_f^p$ and any simple coroot
$\alpha$. 

For all but finitely many $\lambda$ there exists a simple coroot $\alpha$
such that 
$(\lambda,\omega_\alpha)\geq n+1$; thus we can suppose that 
 $\lambda,\alpha$ are like that. If $\nu \la n$ then $\H_\nu\subset F_{\leq
\lambda -\omega_\alpha}$, so we can take $\lambda_1=0$,
$\lambda_2=\lambda-\omega_\alpha$. 

It remains to consider the case  $(\nu, \omega_\alpha) >   n$.
 
We have $\H_\nu=\oplus\H_x$, where $x$ runs over double cosets
 $W_p w W_p\in W_p\backslash
W_{aff}/W_p =K_0\backslash G/K_0$ such that  $w\in 
 W_{aff} \cap (W_f^p\nu W_f^p)$.

Fix one such $w$; let us write it down as
$w=w_1\nu 
w_2$, $w_i\in W_f^p$, and consider the element
 $y=W_p \cdot ( w_1(n\omega_i)w_1^{-1})\cdot W_p\in W_p\backslash
W_{aff}/W_p$. 

It is obvious  that for any three double cosets $z_1$, $z_2$, $z_3$ where
 $z_1=
W_p \cdot w'W_p$,
 $z_2=W_p \cdot w''W_p$,
$z_3=W_p \cdot w' w''W_p\in  W_p\backslash
W_{aff}/W_p$ we have  $\H_{z_1}\cdot \H_{z_2}\supset \H_{z_3}$.

Hence we have $\H_x\subset \H_y\cdot \H_z$, where $z=W_f^p
w_1 (\nu-n\omega_\alpha)w_2 $. 

Since $\H_z\subset \H_{\nu-n\omega_\alpha}\subset
\H_{\preceq(\lambda-n\omega_\alpha) +\lambda_0}\subset F_{\leq
\lambda-n\omega_\alpha }$ we see that
in this case $\H_\nu$ lies in the summand of the RHS of \eqref{nado2}
corresponding to  $\lambda_1=n\omega_i$,
$\lambda_2=\lambda-n\omega_i$. The Proposition is
proved. 
\epf

\begin{Rem} We will also use the following equivalent description of 
$F_{\leq}$. 

For a coweight $\lambda\in \a$ let $\wp(\lambda)$ denote the point
of $\a_+$ closest to $\lambda$.  Then % $\wp(\lambda)\succeq \lambda$;
                                      % moreover 
\begin{equation}\label{Langllem}
\left.
\begin{array}{l}
\wp(\lambda)=\suml_\alpha a_\alpha\omega_\alpha\\
\wp(\lambda)-\lambda=\suml_\beta b_\beta \beta;
\end{array}
\right.
\end{equation}
where $\alpha$ and $\beta$ run over nonintersecting sets of simple roots, and
$a_\alpha>0$, $b_\beta\geq 0$; moreover for any coweight $\lambda$
there exists a unique such decomposition (the latter statement is  the 
so-called  Langlands combinatorial Lemma, see e.g. 
\cite{BW}, Lemma 6.11 on p.143).

It is not hard to see that $\wp(\mu)$ is the minimal element in the set
$\{ \lambda\in \a_+\,|\, \lambda\succeq \mu\}$. It follows that 
\begin{equation}\label{Langlequiv}
%h\in F_{\leq \lambda}(\H)
%\Longleftrightarrow \supp^{geom}(h)\in \cupl_{\wp(\mu-\lambda_0)\leq \lambda }
% G_\mu
 F_{\leq \lambda}(\H)=\oplusl_{\wp(\mu-\lambda_0)\leq \lambda } \H_\mu
\end{equation}
\end{Rem}

\subsubsection{} Another result we will need is the following

\begin{Thm}\label{upada}  Fix $\lambda_1\in\a^+$.
 Make the following assumptions:

\  $(e_1,r_j)>0$ for all  $r_j\in \Sigma^+$;

\  $h_x^\lambda$ (notation of \ref{hxlam}) satisfies the condition
of \ref{hxlam} with large enough a (the bound on $a$ depends on $\lambda_1$); 

\  $\mu\in\a_+$
is very dominant. 

Then we have:

a) If $\nu\in \a_+$ is such that $\<\nu,\Sigma^+-\Sigma^+_L\> \gg 0$
then
$$
%%%\H_\nu\subset  F_{\leq \mu+\lambda_1}(\H),
\wp(\nu-\lambda_0)\leq \mu+\lambda_1,
\,h\in 
%F_{\leq \mu+\lambda_1}( )
\H_\nu
 \Longrightarrow \Psi(h_x^\lambda h)=\iota(\lambda)m_x \Psi(h)+\Psi(h'),$$
where $h'\in F_{\leq \lambda+\mu}(\H).$

b) Otherwise $ h_x^\lambda h\in F_{\leq \lambda+\mu}(\H).$
\end{Thm}

\begin{Lem}\label{tautolog1}  For $L'\subset L$ and
 $f\in C_c^\infty(G/U')_\mu^K$
we have $h_x^\lambda f= \iota(\lambda)m_x (f)+ f'$, where $f'\in \suml_{r\in 
\Sigma
^+-\Sigma^+_L}   C_c^\infty(G/U')_{\preceq \mu+\lambda -a\cdot \P_{L'}(r)}$.
\end{Lem}

\proof is similar to the proof of Lemma \ref{tautolog}. It is enough to
prove that for any irreducible representation $\rho_L$ of $L$ the action of 
$h_x^\lambda -\iota(\lambda)$ on $\Pi_{\rho_L}^K=C_c^\infty (G/U)^K\otimes
_{\H(L^c)}\rho_L$ belongs to $\suml_{\alpha\in\Sigma^+-\Sigma^+_L}\End(\Pi_
{\rho_L})_{\preceq \lambda-ar}.$ The irreducible representation
 $\rho_L$ is a subquotient  
in the induced representation $i_{L'}^L(\rho_{L'})$, where $\rho_{L'}$
is a cuspidal representation of a standard Levi.
 Then $\Pi_{\rho_L}=\Pi_{\rho_{L'}}$, and the statement 
reduces to $I_{\rho_{L'}}(h_x^\lambda)-[\lambda]_{\rho_{L'}}\in
\suml_{\alpha\in \Sigma^+-\Sigma^+_L}
\End(\Pi_{\rho_{L'}}^K)_{\preceq \lambda - a r}$, which follows from 
Proposition \ref{image}.
\epf

\begin{Lem}\label{weight} 
For a weight $\nu\in\a_+$ such that    $\<\nu,\Sigma^+-\Sigma^+
_L\>\gg 0$ we have 
$\wp(\nu+\lambda_1 - a\cdot r - \lambda_0)\leq
\wp(\nu-\lambda_0)$ for $r\in \Sigma^+-\Sigma^+_L$, provided
$a>   \frac{ (e_1,\lambda_1)}{ (e_1,\alpha)}$
for any simple root $\alpha \in   \Sigma^+-\Sigma^+_L$.
\end{Lem}

\proof Let $\alpha\in \Sigma^+-\Sigma^+_L$ be a simple root such 
$\alpha\preceq r$. 

If $\nu$ satisfies  $\<\nu,\Sigma^+-\Sigma^+_L\>\gg 0$ then also
 $\wp(\nu-\lambda_0)+\lambda_1-a\alpha,\Sigma^+-\Sigma^+_L\>\gg 0$.
Since $\<\alpha,\beta\> \leq 0$ for $\beta\in \Sigma^+_L$
it follows that $ \wp(\nu-\lambda_0)+\lambda_1-a\alpha\in \a_+$.
From the condition of Lemma we get $\left( e_1,
 \wp(\nu-\lambda_0)+\lambda_1-a\alpha \right)< \left( e_1,
 \wp(\nu-\lambda_0)\right)$, hence $ \wp(\nu-\lambda_0)+\lambda_1-a\alpha
 \leq  \wp(\nu-\lambda_0)$. The Lemma now follows from \eqref{Langlequiv}.   
\epf

\proof of Theorem \ref{upada} 
follows the same scheme as the proof of \ref{soglasov}.

Let $L'$ be minimal among the standard Levi subgroups for which $\Psi_{L'}|_{G_
\nu}$ is defined. Under the conditions of a) we have $L'\subset L$.

In the notations of the proof of \ref{soglasov}, we can assume that 
$\Psi(h)\in \H(C_{P'})_x$ for some $x\in  \Pbar\backslash G/K$. 
Then using Lemma \ref{tautolog1} together with 
 \ref{tautolog}, \ref{nossogl}  we get 
 $$\Psi(h_x^\lambda h)-
m_x\iota(\lambda) \Psi(h) \in
\suml_{r\in\Sigma^+-\Sigma^+_L}
 \suml_{\eta\preceq \nu+\lambda-a r}\H(C_{P'})_\eta,$$
hence 
 $$\supp^{geom}(h_x^\lambda h-
\Psi^{-1}(m_x\iota(\lambda) \Psi(h)) )\subset
\cupl_{r\in\Sigma^+-\Sigma^+_L}
 \cupl_{\P_{L'}(\eta)\preceq \nu+\lambda-a r}G_\eta.$$

Let, first, $\eta$ be such that $G_\eta\cap \supp^{geom}(h_x^\lambda h-
\Psi^{-1}(m_x\iota(\lambda) \Psi(h)) )\ne\emptyset$. As we have just seen
there exists $r\in\Sigma^+-\Sigma^+_L$ such that $\eta\la \nu+\lambda-ar$
for any simple coroot $\alpha\in \Sigma^+-\Sigma^+_{L'}$. By Lemma \ref{weight}
$\nu+\lambda-ar\preceq \lambda_0+\lambda'$, where $\lambda'\leq \mu+\lambda$.

Thus for a  simple coroot $\alpha\in \Sigma^+-\Sigma^+_{L'}$ we have
$\eta\la \nu+\lambda-ar\la  \lambda_0+\lambda'$.

On the other hand, if $\alpha$ is a simple coroot  of $L'$
then $\<\omega_\alpha,\nu\>$ is bounded. Since $\lambda_0$ is very dominant
we have $\eta\la \lambda_0\la \lambda_0+\lambda'$.

So $\eta\la \lambda_0+\lambda'$ for any simple coroot $\alpha$, i.e.
 $\eta\preceq \lambda_0+\lambda'$, which implies $h_x^\lambda h-
\Psi^{-1}(m_x\iota(\lambda) \Psi(h))\in F_{\leq\mu+\lambda}(\H)$. a) is proved.

Assume now we are in the situation of b). Then for some simple root $\alpha\in 
\Sigma_L^+$ the number $\<\nu,\omega_\alpha\>$ is bounded.
 
The condition $(e_1,\alpha)>0$ implies that we
 can find $a\in \RE_{>0}$ such that $\lambda_1-a\alpha\leq 0$; we can then 
assume $\lambda\succeq a\alpha$. 
Also we can assume (after possibly increasing $\lambda_0$) that 
$\H_{\preceq \lambda}\cdot \G_{\preceq \mu}(\H)\subset \G_{\preceq\lambda+\mu}$
holds already for $\mu\succeq \lambda_0-a\alpha$. By the definition
$\nu \preceq \lambda_0+\zeta$ for some $\zeta\leq \mu+\lambda_1$. But
since $\<\omega_\alpha,\nu\>$ is bounded we also have
 $\nu \preceq \lambda_0+\zeta-a\alpha$.
 Then we have
$h_x^\lambda h\in \G_{\preceq  \lambda_0+\zeta-a\alpha+\lambda}\subset
F_{\leq \lambda-a\alpha+\zeta}= F_{\leq \lambda+\mu+(\lambda_1-a\alpha)}
\subset  F_{\leq \lambda+\mu}$.
\epf

\subsection{End of the proof: orbital integrals}
\begin{Thm}\label{sumEii}
The functional $Tr(g,M)$ defined above satisfies \ref{iii}iii).
\end{Thm}

 We start with recalling the following standard computation.

 We will say that a  subset $\cS \subset \Lambda   ^+$ is large  if
it contains $(\Lambda^+)_{\preceq \lambda}$ for large $\lambda   $.

For any  
$\cS \subset \Lambda^+$ denote $G_\cS = \bigcup _{\lambda   \in \cS}
G_\lambda$,
 $\H_\cS=\bigoplus\limits
_{\lambda\in \cS} \H_\lambda$
and let $pr_\cS:\H^{\oplus m}\to \H^{\oplus m}_\cS \subset \H^{\oplus m}$ be
 the
projector along $\H^{\oplus m}_{\Lambda^+-\cS}$.

\begin{Lem}\label{standcomputa} Keep the assumptions of \ref{sumEii}. Then
 for any large enough finite set $\cS \subset \Lambda^+$
we have
$$O_{g^{-1}}(\sum E_{ii})=tr(g\circ E\circ pr_\cS)$$
(the RHS is the well-defined trace of a finite-rank endomorphism of
$\H^{\oplus m}$.)
\end{Lem}
\proof of the Lemma.
We have $
tr(g\circ E\circ pr_\cS, \H^{\oplus m})=
tr(\delta_{g\cdot K}\circ E\circ pr_\cS, k[G/K]^{\oplus m})=
\int_{g'\in \K\cdot g}tr( g'\circ E\circ  pr_\cS^{\oplus m},
k[G/K]^{\oplus m}),$
where we reused the symbol $pr_\cS$ to denote also the
projector
  $k[G/\K]^{\oplus m}\to k[G_\cS/K]^{\oplus m}\subset k[G/K]^{\oplus m}$
 along 
 $k[G_{\Lambda^+-\cS}/K]^{\oplus m}$.

In $k[G/\K]$ we have the standard basis consisting of delta-functions;
 it gives also a basis in $k[G/\K]^{\oplus m}$.
 We use this basis to
``compute'' the trace of $ g'\circ E\circ  pr_\cS^{\oplus m}$.
 Precisely, let $dg$ be the Haar 
measure on $G$ such that $\int_{\K} dg=1$; if $x\in  G/\K$,
 and $x_{(i)}= (\underbrace{0,\dots,0}_{i-1},\delta_x,
\underbrace{0,\dots,0}_{m-i})$ is a basis element
of $\H^{\oplus m}$, then the corresponding diagonal entry of 
$g'\circ E\circ pr_\cS^{\oplus m}$ is   equal to

\begin{equation} \label{intTr}
\left\{   
\begin{array}{l}
   0 \ \ \operatorname{if}
 \ \  x\not \in G_{\cS}/\K \\
  {E_{ii} \over dg}\left( x^{-1}(g')^{-1}x \right)\ \ \operatorname{if} 
\  \
 x\in G_{\cS}/\K 
\end{array}
\right.
\end{equation}
If  $\K\cdot g\owns g'$  is contained in the regular elliptic set,
then for   $\cS$ large enough we have 
 ${E_{ii}\over dg} ( x^{-1}(g')^{-1}x) \not = 0
\Rightarrow x \in  G_{ \cS}/\K $. 
For such 
$\cS$ summation of \eqref{intTr} over $x$ gives
$O_{(g')^{-1}}(E_{ii})$.

Altogether we get
$$tr(g\circ E\circ pr_\cS, \H^{\oplus m})=tr(\int_{g'\in\K\cdot g}
g'\circ E\circ pr_\cS, k[G/\K]^{\oplus m})dg'=
$$
$$=\int_{g'\in\K\cdot g}
O_{(g')^{-1}}(\sum E_{ii})=O_{g^{-1}}(\sum E_{ii})$$
where  the last equality follows from the assumption that
$O_{(g')^{-1}}(\sum E_{ii})$ is
constant for $g'\in g\cdot K$.  \epf

We proceed now to the proof of the Theorem.

Fix very dominant $\lambda_0$, and
set $\cS_\lambda= \cupl  _{\mu \leq \lambda} X   ^+_{\preceq \mu
+ \lambda_0}$.
Then $ \H^{\oplus m}_{\cS_\lambda}=F_{\leq \lambda}(\H)^{\oplus m}$
is a good filtration on the free module of rank $m$ by \ref{F}.

Thus the image of $F_{\leq \lambda}(\H)^{\oplus m}$ under $E$ is a good
filtration on $M$, hence we have $Tr(g,M)= tr\(g, E(F_{\leq
\lambda}(\H)^{\oplus m} ) \)$ for large $\lambda$.

So to prove the Theorem we need only to show that
\begin{equation} \label{twosides}
tr(g,E(\H^{\oplus m}_{\cS_\lambda}))=
tr(g\circ E\circ pr_{\cS_\lambda},
\H^{\oplus m}).
\end{equation}
It is convenient to consider the vector space
 $E(\H^{\oplus m}_{\cS_\lambda})+
\H^{\oplus m}_{\cS_\lambda}$
 which is invariant under the action of $E,\,
pr_{\cS_\lambda},\, g$. Obviously
                                     $$
tr(g,E(\H^{\oplus m}_{\cS_\lambda}))=tr(g\circ E,
E(\H^{\oplus m}_{\cS_\lambda})+
\H^{\oplus m}_{\cS_\lambda}),$$
$$ tr(g\circ E\circ pr_{\cS_\lambda},
\H^{\oplus m})=
 tr(g\circ E\circ pr_{\cS_\lambda},
E(\H^{\oplus m}_{\cS_\lambda})+
\H^{\oplus m}_{\cS_\lambda}).$$

Thus the difference of the two sides of \eqref{twosides}
coincides with 
\begin{multline*}
tr(g\circ E\circ (Id-pr_{\cS_\lambda}),
E(\H^{\oplus m}_{\cS_\lambda})+
\H^{\oplus m}_{\cS_\lambda})= tr(g\circ E\circ pr_{X   _+-\cS_\lambda},
E(\H^{\oplus m}_{\cS_\lambda})+
\H^{\oplus m}_{\cS_\lambda})=\\
 tr(g\circ  pr_{X   _+-\cS_\lambda}\circ E,
\(E(\H^{\oplus m}_{\cS_\lambda})+
\H^{\oplus m}_{\cS_\lambda}\)\cap\H^{\oplus m}_{X   ^+-\cS _\lambda} )
, \end{multline*}
 so the Theorem
 is  reduced to the following assertion.

\begin{Prop}\label{traceless}
 The endomorphism $g\circ pr_{X   ^+-\cS _\lambda} \circ E$
of the vector space
 $\( E(\H^{\oplus m}_{\cS_\lambda})+\H^{\oplus m}_{\cS_\lambda} \) \cap
\H^{\oplus m}_{X   ^+-\cS _\lambda}$
is traceless for large $\lambda$.
\end{Prop}
 Denote  $$
\( E(\H^{\oplus m}_{\cS_\lambda})+\H^{\oplus m}_{\cS_\lambda} \) \cap 
\H^{\oplus m}_{X   ^+-\cS _\lambda}=pr_{X   ^+-\cS
_\lambda}\(E(\H^{\oplus m}_{\cS_\lambda})\)\iso\(
E(\H^{\oplus m}_{\cS_\lambda})+\H^{\oplus m}_{\cS_\lambda} \)/
\H^{\oplus m}_{\cS_\lambda} $$
 by $W_\lambda$.  For $\nu \leq \lambda$ let us
denote by $(W_\lambda)_{\leq \nu}\subset W_\lambda$ the subspace
$\(E\(F_{\leq \nu}(\H^{\oplus m})\) +F_{\leq \lambda}(\H^{\oplus m})\)
/F_{\leq
\lambda}(\H^{\oplus m})= \(E\(F_{\leq \nu}(\H^{\oplus m})\) +
\H^{\oplus m}_{\cS_\lambda} \)/
\H^{\oplus m}_{\cS_\lambda}
$. 

For a number $N$ let $\lambda_N$ be the maximal element in the (finite)
set $X   _N:=\{\lambda \in X   ^+\ |\ ( \lambda, e_1)\leq N
\}$. Then  
$ X   ^+_{\leq \lambda_N}=X   _N$; since $X   _N$ is arbitrary
large for $N\gg 0$ it is enough (in view of \ref{standcomputa})
to prove the Proposition for $\lambda=\lambda_N$. We write $W_N, \cS_N$
 instead of $W_{\lambda_N}, \cS_{\lambda_N}$.

Consider the associated graded space $grW:=\oplusl_N \oplusl _{ \lambda \in
X   _ N} (W_N)_{\leq \lambda}/(W_N)_{< \lambda}$. It carries a natural
 action of $gr\H$ provided by the obvious maps $gr\H_\mu \otimes  (W_N)_{\leq
\lambda}/(W_N)_{< 
\lambda} \to   (W_{N+(e_1, \mu)} )_{\leq\mu+ \lambda}/(W_{N+(e_1,
\mu)} )_{<\mu+
\lambda}$ which are induced by multiplication in $\H$.

\begin{Lem}
$grW$ is a finitely generated $gr\H$-module.
\end{Lem}

\proof It is obvious that $gr(W_N)_\lambda \ne 0 \Rightarrow N-n\leq
\<\lambda ,e_1\>\leq N$ where $n$ is some constant (depending on
$\supp^{geom}(E_{ij})$). Thus $grW$ is a finite sum of modules
$grW^{(i)}:= \oplusl_{N-(\lambda,e_1)=i} gr(W_N)_\lambda$, 
each one of which is a quotient of $gr_F(\H^{\oplus m})$, so
 is finitely generated.
\epf

Now by the proof of \ref{tr0}b)
 we can break $\Lambda_0^+$ into a finite union of
cosets 
\begin{equation}\label{decompos}
\Lambda_0^+=\bigcup \mu_i+ X   _{L_i}^+
\end{equation}
 so that $\lambda \in \mu_i+
X   _{L_i}^+ \Rightarrow grW_\lambda= M_i$ where $M_i$ is a finite
$\A_{L_i}$-module (notations of \eqref{AL}). 
  
Notice that for $\lambda\in \a_+$ the set $\cS_\lambda-\cS_{<\lambda}$
looks as follows. Suppose that $\lambda$ lies strictly inside $\a_L^+$.
Then $\cS_\lambda-\cS_{<\lambda}=(\lambda_0+\wp^{-1}(\lambda-\lambda_0))
\cap \Lambda^+=\Lambda^+\cap
\{\lambda-\sum a_ir_i\,|\,0\leq a_i\leq \<\omega_i, \lambda_0\>\}$
 where $r_i$ runs over the set of simple 
roots of $L$.

In particular we have $\< \lambda,\Sigma^+-\Sigma^+_L\>\gg 0,\,
\nu\in \cS_\lambda-\cS_{<\lambda}\Rightarrow \<\nu,\Sigma^+-\Sigma^+_L\>\gg 0$.

Proceeding by induction in $\dim(\a_L)$   %%??
 we can choose decomposition
\eqref{decompos} so that $\<\nu,\Sigma^+-\Sigma^+_{L_i}\>$ is arbitrary
large if $\nu\in\cS_\lambda-\cS_{<\lambda}$ where $\lambda\in \mu_i+X_{L_i}^+$.
In particular we can guarantee that $G_{\cS_\eta-\cS_{<
\eta}}$   lies in the domain of $\Psi:K\backslash G/K\to K\backslash (C_{P_i})
/K$
for $\eta \in \mu_i+X   
_{P_i}^+$.  We shall make this
assumption from now on. 

Then  we get a direct sum decomposition:
\begin{equation}\label{decompo}\H^{\oplus m}_{\cS_\eta-\cS_{<
\eta}}= \oplusl _{x\in  K\backslash G/P} \Psi^{-1}
C^\infty (C_P^{\cS_\eta-\cS_{< \eta}})_x^{\oplus m}
\end{equation}

%% say that X   _G=0, A_G=k etc.

So take $\eta \in \mu_i+ X   _{P_i}^+$, and
consider the  surjection $\pi:\H^{\oplus m}_{\cS_\eta-\cS_{< \eta}}
\to gr (W_N)_\eta\cong E(\H^{\oplus m}_{\cS_\eta-\cS_{< \eta}})
/  \left( E(\H^{\oplus m}_{\cS_\eta-\cS_{< \eta}})\cap
( E(\H^{\oplus m}_{\cS_{< \eta}})+ \H^{\oplus m}_{\cS_N})
\right)$.

Since its  target is an $\A_{L_i}$-module, it also decomposes as a direct
sum: 
$ gr (W_N)_\eta= \oplusl _{x\in  K\backslash G/P}  gr (W_N)_\eta^x$.
%%Napish, gad!!??
It is not hard to see that $\pi$ commutes with the $\A_{L_i}$-action, and hence
sends respective summands to respective
summands (i.e. $\pi\Psi^{-1}\(C^\infty
(C_P^{\cS_\eta-\cS_{< \eta}})_x\)\subset
 gr (W_N)_\eta^x$).

Hence for each $N,\eta$ we can choose a subspace $(V_N)_\eta\subset
\H^{\oplus m}_{\cS_\eta-\cS_{<
\eta}}$ in such a way that:

a) $\pi$ maps $(V_N)_\eta$ to $gr(W_N)_\eta$ isomorphically.

b)  $(V_N)_\eta =\oplusl _{x\in  K\backslash G/P} \left(\Psi^{-1}
C^\infty(C_P^{\cS_\eta-\cS_{< \eta}})_x^{\oplus m}\right) \cap
(V_N)_\eta $.

c) $(V_N)_\eta$ is $K_0$ invariant.

d) For $\mu \in X   _{L_i}^+$ the isomorphism $\Psi_{P_i}^{-1}\circ
 \iota(\mu)\circ \Psi_{P_i}: 
\H^{\oplus
m}_{\cS_\eta-\cS_{< \eta}} \cong
\H^{\oplus m}_{\cS_{\mu+\eta}-\cS_{<
\mu+\eta}}$ sends  $(V_N)_\eta$ to  $(V_{(e_1,\mu)+N})_{\mu+\eta}$.

Set now $(W_N)_\eta=E\((V_N)_\eta\)$. By property a) we have a
direct sum decomposition $W_N =
\oplusl_\nu (W_N)_\nu$.

 Let $E_N^\nu\in End\((W_N)_\nu\)$
 denote the corresponding diagonal block of $pr_{X   ^+-\cS_N}\circ E$.

To finish the proof we need the following (key)

\begin{Lem}\label{transla}
For  $N$ large %(depending on $m$, $\lambda$, $E$??)
and any  $x\in P\backslash G/K$ we have:

a)  If $\<\nu,\Sigma^+-\Sigma^+_L\>\gg 0$ then
$h_{x,\lambda}    \cdot (W_N)_\nu \subset
(W_{N+(e_1,\lambda)})_{\nu+\lambda}+\H_{\cS_{N+(e_1,\lambda)}}$.

b) $h_{x,\lambda}   \cdot
(W_N)_\nu \subset \H_{\cS_{N+(e_1,\lambda)}}$ otherwise.
\end{Lem}
\proof By \ref{odinak} property d) implies that $\Psi_P(W_N)_{\nu+\lambda}
=\iota(\lambda)\Psi_P(W_N)_{\nu}$ if $\<\nu,\Sigma^+-\Sigma^+_L\> \gg 0$.
Also from the definition it is clear that under the same condition
on $\nu$ the space $\Psi((V_N)_\nu)$, and hence $\Psi(
(W_N)_\nu)$
is invariant under the action of $m_x$. So a) follows from Theorem
\ref{upada}a) (take $\lambda_1$ to be any coweight satisfying $E_{ij}
\in \H_{\leq \lambda_1}$). 

Statement b) follows Theorem \ref{upada}b).

\begin{Cor}\label{commut}  
$E_N^\nu$ commutes with the action of 
 $\A_L$ on $(W_N)_\nu\iso gr(W_N)_\nu$ provided $\<\nu,\Sigma^+-\Sigma_L^+\>
\gg 0$.
\end{Cor}
\proof Fix $x\in K\backslash G/P$ 
and consider the element $h_{x}^\lambda   \in \H$ (see \ref{hxlam}),
where $\lambda\in X_{L}^+$ is dominant enough as in \ref{hxlam}.

 We have $h'\in (W_N)_\nu \Rightarrow h'\cdot E= E_N^\nu (h')+\suml
_{\nu'\ne \nu}h'_{\nu'}+h_0$ where $h'_{\nu'}\in (W_N)_{\nu'}$
and $h_0 \in \H_{\cS_N}$. Multiplying the last equality 
by $h_x^\lambda$ on the left and applying \ref{transla} we get the statement.
\epf

\proof of \ref{traceless}. From \ref{commut} it follows that for large $N$
and any $\nu$ the endomorphism $E_N^\nu$
preserves the direct sum decomposition $(W_N)_\nu=\oplusl_{x\in K\backslash
G/P}
\((W_N)_\nu \)_x$. Since $g$ permutes these summands without fixed points,
the endomorphism $g\circ E_N^\nu$ is traceless. Hence 
 $tr(g\circ pr_{X   ^+-\cS _N} \circ E)=\suml_\nu tr(g\circ
E_N^\nu)=0$.

Proposition \ref{traceless}, and thus also Theorem \ref{sumEii} is proved.
\epf 

\medskip

We can now summarize our results in the following

\begin{Thm}\label{cocentr}  Let $\Mmod$ be a projective $G$-module, 
and let $\K\subset K_0$ be a normal open  subgroup,  so  small that:

$\ \alpha)$ $g$ is $K$-elliptic, so $Tr(g,\Mmod^K)$ is defined 
(see \ref{Kellipt}, \ref{iii}).

$\ \beta)$  $\Mmod$ is generated by its $K$-fixed vectors.
 
$\ \gamma)$ Conditions of \ref{iii}iii) are satisfied for $M=\Mmod^K$.

Then 
\beq\label{TrOg}
Tr(g,\Mmod^\K)=O_{g^{-1}}\(\< \Mmod)\> \).
\end{equation}
\end{Thm}

\proof There exists an idempotent $E\in Mat_m(\H)$, such that
$\Mmod$ is isomorphic to the image of $E$ acting on the right on $k[G/K]
^{\oplus m}$. Since $K$ is nice by the assumption (see \ref{iii}iii), we have
 $\<\Mmod\>= \suml_i E_{ii} \mod [\H(G),\H(G)]$, so
the statement follows from Theorem \ref{sumEii}. \epf

%VSE PEREPISAT'

\subsubsection{Proof of Theorem \ref{Main}}\label{proofMain} 
Let $0\to P_n \to \dots \to P_0\to \rho \to 0$ be a projective resolution
of $\rho$. Pick a nice open compact $K$, such that $g$ is $K$-elliptic,
properties $\beta)$ and $\gamma)$ of \ref{cocentr} hold for $\Mmod=P_i$,
and the character $\chi_\rho$ is constant on the coset $g\cdot K$.
Then by additivity of $Tr(g,M)$ (\ref{iii}(i)=\ref{Tr}(b)) 
and \ref{cocentr} we get
$Tr(g,\rho^K)=O_{g^{-1}}(\<\rho\> )$. On the other hand, by  
\ref{iii}(ii)=\ref{Tr}(c)
we have $Tr(g,\rho^K)=\chi_\rho(g)$. The proof is finished. \epf

\subsection{Elliptic pairing}\label{ellsect} It  remains to
deduce Theorem \ref{Kazhdan} from
 \ref{Main}.
This argument  can also be found in \cite{Schn1} and is included here for
 the sake of completeness.

We set $k=\Ce$, $\char(F)=0$ till the end of the argument.

Recall that the measure $d\mu$ on the set $Ell$ of regular semisimple elliptic
conjugacy classes is characterized by the equality
\begin{equation}\label{Weyl} \int\limits_{Ell} O_g(h)d\mu(g)=\int\limits_G h
\end{equation}
being true for all $h\in \H$ supported inside the set of regular elliptic
elements.

\begin{Lem}\label{Mrho} For any $\Mmod \in \M$ and $\rho\in \R$ we have
$$
 \sum (-1)^i \dim
Ext^i(\Mmod,\rho)=
\int_{G} \chi_{\rho}\cdot \< \Mmod \>.
$$
\end{Lem}
\prf Since both sides are additive on short exact sequences in $\Mmod$,
it is enough to consider the case of  projective $\Mmod$. 
As usual we write $\Mmod$ as the image of the right action of an idempotent
$E\in Mat_m(\H)$ on $\H^{\oplus m}$.
 We have
\begin{equation}\label{Hom_G}
\Hom_G(\H,\rho)
 \cong \hat{\rho},
\end{equation}
 where $\hat {\rho}\overset{\fed}{=} \varprojlim_K (\rho_K)$. (For $\rho\in \R$
we have also $\hat{\rho}
=*(\rho\galochka)$).
%=\varprojlim_{V}
%\subset \rho,\codim(V)< \infty}
%\rho/V,$ where $V$ runs over the set of subspaces of finite codimension
%in the space of $\rho$.
The right action of $Mat_n(\H)$ on $\H^{\oplus n}$ induces an action
of $Mat_n(\H)$ on $\Hom_G(\H^{\oplus n } , \rho),$
$m(f)(x):=f(m(x))$. Evidently, this action agrees under \eqref{Hom_G}
 with the action of
$Mat_n(\H)$ on $\hat{\rho}^{\oplus n}$ inherited from the action of $\H$
on $\rho$. 
Thus we have \begin{multline}
\dim (\Hom_G(\Mmod, \rho))=\dim (\Hom_G(\Im(E,\H^{\oplus n}),\rho))=
\dim(\Im(E, \hat{\rho}^{\oplus n}))=\\
\dim(\Im(E, {\rho}^{\oplus n}))=
tr(E,\rho^{\oplus n})= \int_G \chi_\rho\cdot \sum_i E_{ii}.
\end{multline}
The Lemma is proved.

\begin{Lem}\label{nonell}
Let $\Phi$ be an invariant generalized function on $G$ supported on the set
of regular non-elliptic elements. Then for any admissible representation
$\rho$ we have:
$$\int_G  \<\rho\>\cdot \Phi =0.$$
\end{Lem}
\prf  By Theorem 0 in \cite{K} it is enough to show that 
$$\int\limits_G \<\rho\> \cdot \chi_{i_L^G(\rho_L)}=0$$
for any admissible representation $\rho_L$ of a proper standard
 Levi subgroup $L$.
By Lemma \ref{Mrho} it we have to show that
$$\suml_i (-1)^i \dim Ext^i_G(\rho, i_P^G(\rho_L))=0.$$
By Frobenius adjointness we have $Ext^i_G(\rho, i_P^G(\rho_L))\cong 
Ext^i_L(r^G_L(\rho),\rho_L)$. The Lemma now follows from the next

\begin{Claim}\label{noncomp} \cite{Bunp} If the center of $G$ is not compact,
then for any admissible representations $\rho_1,\rho_2\in \R(G)$ we have
$$\suml_i (-1)^i \dim Ext^i_G(\rho_1, \rho_2)=0.$$
\end{Claim}
\proof For any $\Mmod\in Sm(G)$ the space $Ext^i(\Mmod,\rho_2)$ is finite
dimensional. Thus the map $[\Mmod]\mapsto  
\suml_i (-1)^i \dim Ext^i_G(\Mmod, \rho_2)$ is a well-defined homomorphism
$K^0(\M(G))\to \Zet$.

So we will be done if we show that the class  of $\rho_1$, 
$[\rho_1]\in K^0(\M(G))$ vanishes. 

$G/G^c$ is a free abelian group of rank equal to the split rank of the center
of $G$. Hence, if the center of $G$ is not finite, $G/G^c$ is nontrivial,
and there exists a nontrivial homomorphism $G\to \Zet$ with open kernel.
Let us  endow $k[\Zet]$  with the structure of a smooth $G$-module
by means of this homomorphism. 
 Consider the short exact sequence
of $G$-modules: $0\to k[\Zet]\overset{[n]\to[n+1]}{\longrightarrow}
 k[\Zet]\to k\to 0$. Tensoring it with $\rho_1$ we get
 \begin{equation}\label{0}
0\to \rho_1\otimes k[\Zet]\to
\rho_1\otimes k[\Zet]\to \rho_1\to 0.
\end{equation}

It is easy to see that  $\rho_1\otimes k[\Zet]$ is finitely generated provided
$\rho_1$ is admissible. Thus \eqref{0} implies that $[\rho_1]=0$. \epf

\subsubsection{Proof of \ref{Kazhdan}} By Lemma \ref{Mrho} $$
 \sum (-1)^i \dim
Ext^i(\rho_1,\rho_2)=
\int_{G} \chi_{\rho_2}\cdot \< \rho_1 \>.$$
By Harish-Chandra's Theorem on
integrability of characters we can rewrite $$
\int_{G} \chi_{\rho_2}\cdot \< \rho_1 \>
=\lim\limits_U
\int_{U} \chi_{\rho_2}\cdot \< \rho_1 \>,
$$
where $U$ runs over an exhausting family of conjugation-invariant subsets of
 regular elements.
Each $U$ in this family can be written as a disjoint union $U=U^{ell}\cup
U^{nonell}$ of a subset of the elliptic (respectively, non-elliptic)
set. We have
$\int_{U^{nonell}} \chi_{\rho_2}\cdot \< \rho_1 \>=0,$
for all $U$ by \ref{nonell}, while $$ \int_{U^{ell}}
 \chi_{\rho_2}\cdot \< \rho_1 \>=\int_{U^{ell}/Ad}\chi_{\rho_2}(g)
O_g(\<\rho_1\>)d\mu=\int_{U^{ell}/Ad}\chi_{\rho_2}(g)\chi_{\rho_1}(g^{-1})d\mu,$$
where the first equality is \eqref{Weyl}, and the second follows from 
\ref{Main}.
As $U^{ell}$ increases, the latter expression tends to $\int_{Ell} \chi_{\rho_1}
(g^{-1}) \chi_{\rho_2}(g)d\mu(g).$ The proof is finished.

\end{section}

\end{document}